\documentclass[11pt]{article}
\usepackage{amsmath,amsthm,amssymb,amsfonts,amscd}
\usepackage{amsxtra}
\usepackage{eucal}
\usepackage{mathrsfs}
\usepackage{color}
\usepackage{xcolor}
\usepackage{hyperref}
\usepackage{enumerate}
\usepackage{blkarray}
\usepackage{enumitem}
\usepackage{tcolorbox}
\usepackage{tikz,tikz-cd}
\usetikzlibrary{3d}
\usepackage{tikz-3dplot}
\usetikzlibrary{positioning}
\usetikzlibrary{decorations.pathreplacing,decorations.markings,decorations.pathmorphing}

\tikzset{
  on each segment/.style={
    decorate,
    decoration={
      show path construction,
      moveto code={},
      lineto code={
        \path [#1]
        (\tikzinputsegmentfirst) -- (\tikzinputsegmentlast);
      },
      curveto code={
        \path [#1] (\tikzinputsegmentfirst)
        .. controls
        (\tikzinputsegmentsupporta) and (\tikzinputsegmentsupportb)
        ..
        (\tikzinputsegmentlast);
      },
      closepath code={
        \path [#1]
        (\tikzinputsegmentfirst) -- (\tikzinputsegmentlast);
      },
    },
  },
  mid arrow/.style={postaction={decorate,decoration={
        markings,
        mark=at position .5 with {\arrow[#1]{stealth}}
      }}},
}

\tikzstyle{v}=[circle, draw, solid, fill=black, inner sep=0pt, minimum width=3pt]
\tikzstyle{vo}=[circle, draw, solid, fill=white, inner sep=0pt, minimum width=3pt]

\tikzset{snake it/.style={decorate, decoration=snake}}

\tikzset{box1/.style={draw=black, fill=gray!30, thick, rectangle, rounded corners, minimum height=0.63cm, minimum width=0.63cm}}
\tikzset{box2/.style={draw=black, dotted, rectangle, rounded corners, minimum height=0.63cm, minimum width=0.63cm}}

\tikzset{box1s/.style={draw=black, fill=gray!30, thick, rectangle, rounded corners=0.07cm, minimum height=0.3cm, minimum width=0.3cm}}
\tikzset{box2s/.style={draw=black, dotted, rectangle, rounded corners=0.07cm, minimum height=0.3cm, minimum width=0.3cm}}

\usepackage{algorithm}
\usepackage{algpseudocode}
\usepackage[a4paper,margin=2cm]{geometry}
\definecolor{forestgreen(traditional)}{rgb}{0.0, 0.27, 0.13}
\definecolor{forestgreen(web)}{rgb}{0.13, 0.55, 0.13}
\definecolor{airforceblue}{rgb}{0.36, 0.54, 0.66}

\newcommand{\sgn}{\mathrm{sign}}

\newcommand{\symdiff}{\triangle}
\newcommand{\supp}{\mathrm{supp}}
\newcommand{\pf}{\mathrm{pf}}

\newcommand{\antlift}{\mathrm{ant}}

\newcommand{\cA}{\mathcal{A}}
\newcommand{\cB}{\mathcal{B}}
\newcommand{\cC}{\mathcal{C}}
\newcommand{\cD}{\mathcal{D}}

\newcommand{\cG}{\mathcal{G}}

\newcommand{\cP}{\mathcal{P}}

\newcommand{\cS}{\mathcal{S}}
\newcommand{\cT}{\mathcal{T}}

\newcommand{\bF}{\mathbb{F}}

\newcommand{\bI}{\mathbb{I}}

\newcommand{\bK}{\mathbb{K}}

\newcommand{\bN}{\mathbb{N}}

\newcommand{\bP}{\mathbb{P}}

\newcommand{\bR}{\mathbb{R}}
\newcommand{\bS}{\mathbb{S}}
\newcommand{\bT}{\mathbb{T}}
\newcommand{\bU}{\mathbb{U}}

\newcommand{\bZ}{\mathbb{Z}}

\newcommand{\lift}{\mathrm{lift}}
\newcommand{\skewpair}[1]{\{#1,#1^*\}}

\newcommand{\sympl}[2]{\omega(#1,#2)}
\newcommand{\ortho}[2]{\alpha(#1,#2)}
\newcommand{\ground}[1]{[#1]\cup[#1]^*}

\newcommand{\ul}[1]{\underline{#1}}

\newcommand{\dist}{\mathrm{dist}}

\newcommand{\smaller}[2]{|#1<#2|}
\newcommand{\smallereq}[2]{|#1\le #2|}
\newcommand{\gr}[3]{\mathrm{Gr}_{#1}(#2,#3)}
\newcommand{\lag}[2]{\mathrm{SpGr}_{#1}(#2,2#2)}
\newcommand{\dressian}[2]{\mathrm{Dr}(#1,#2)}

\hypersetup{
  colorlinks,
  citecolor=forestgreen(web),
  linkcolor=airforceblue,
  urlcolor=magenta,
  pdftitle={Baker--Bowler theory for Lagrangian Grassmannians},
  pdfauthor={Donggyu Kim}}

\newtheorem{theorem}{Theorem}[section]
\newtheorem{lemma}[theorem]{Lemma}

\newtheorem{claim}{Claim}[theorem]
\newtheorem{corollary}[theorem]{Corollary}
\newtheorem{proposition}[theorem]{Proposition}

\newtheorem*{theorem-bij}{Theorem~\ref{thm: bij}}

\theoremstyle{definition}
\newtheorem{definition}[theorem]{Definition}
\newtheorem{example}[theorem]{Example}

\theoremstyle{remark}
\newtheorem{remark}[theorem]{Remark}

\usepackage{authblk}
\usepackage{lineno}
\title{Baker--Bowler theory for Lagrangian Grassmannians}
\author[ ]{Donggyu Kim}
\affil[ ]{School of Mathematics, Georgia Institute of Technology, Atlanta, USA}
\affil[ ]{E-mail: \texttt{donggyu@gatech.edu}}
\date{\today}	

\begin{document}
\maketitle
\begin{abstract}
Baker and Bowler (2019) showed that the Grassmannian can be defined over a tract, a field-like structure generalizing both partial fields and hyperfields.
This notion unifies theories of matroids over partial fields, valuated matroids, and oriented matroids.
We extend Baker--Bowler theory to the Lagrangian Grassmannian which is the set of maximal isotropic subspaces in a $2n$-dimensional symplectic vector space.
By Boege et al. (2019), the Lagrangian Grassmannian is parameterized as a subset of the projective space of dimension $2^{n-2}(4+\binom{n}{2})-1$ and its image is cut out by certain quadrics.
We simplify a list of quadrics so that these are apparently induced by the Laplace expansions only concerning principal and almost-principal minors of a symmetric matrix.
From the idea that the strong basis exchange axiom of matroids captures the combinatorial essence of the Grassmann--Pl\"{u}cker relations, we define matroid-like objects, called antisymmetric matroids, derived from the quadrics for the Lagrangian Grassmannian.
We also provide a cryptomorphic definition in terms of circuits capturing the orthogonality and maximality of a Lagrangian subspace.
We define antisymmetric matroids over tracts in two equivalent ways, which generalize both Baker--Bowler theory and the parameterization of the Lagrangian Grassmannian.
It provides a new perspective on the Lagrangian Grassmannian over hyperfields such as the tropical hyperfield and the sign hyperfield.
Our proof involves a homotopy theorem for graphs associated with antisymmetric matroids, which generalizes Maurer's homotopy theorem for matroids.
We also prove that if a point in the projective space satisfies the $3$-/$4$-term quadratic relations for the Lagrangian Grassmannian and its supports form the bases of an antisymmetric matroid, then it satisfies all quadratic relations, a result motivated by the earlier work of Tutte (1958) for matroids and the Grassmannian.
\end{abstract}

\section{Introduction}\label{sec: intro}

For a field $k$ and integers $0\le r\le n$, the Grassmannian $\gr{k}{r}{n}$ is the set of $r$-dimensional linear subspaces in the $n$-dimensional vector space $k^n$, which can be parameterized as a subset of the projective space of dimension $\binom{n}{r}-1$ by the Grassmann--Pl\"{u}cker embedding $p$.
The image of $p$ is exactly cut out by the Grassmann--Pl\"{u}cker relations, which are homogeneous quadrics.
For every linear subspace $V\in \gr{k}{r}{n}$, the support of the Pl\"{u}cker vector $p(V)$ forms the bases of a matroid $M$, and the set of the minimal supports of nonzero vectors in $V$ forms the circuits of the dual matroid $M^\perp$.
Therefore, matroids are regarded as the combinatorial essence of linear spaces.
A similar combinatorial abstraction exists for the Lagrangian orthogonal Grassmannian $\mathrm{OGr}_k(n,2n)$, which is the set of maximal isotropic subspaces of $k^{2n}$ equipped with the standard symmetric bilinear form; see~\cite{Bouchet1988repre,Wenzel1993b,Wenzel1996b,BJ2023,JK2023}.
In this article, we explore a combinatorial structure established on the Lagrangian symplectic Grassmannian.

Let $E := \ground{n} = \{1,2,\ldots,n\} \cup \{1^*,2^*,\ldots,n^*\}$, and let $\omega$ be the standard symplectic form on $k^{2n} = k^{E}$, i.e., $\sympl{X}{Y} := \sum_{i=1}^n (X(i)Y(i^*) - X(i^*)Y(i))$.
The Lagrangian (symplectic) Grassmannian $\lag{k}{n}$ is the set of maximal isotropic subspaces in~$k^{2n}$, which are $n$-dimensional and are also called Lagrangian subspaces.
De Concini~\cite{DeConcini1979} showed that $\lag{k}{n}$ is parameterized as a subset of the projective space of dimension $\binom{2n}{n}-1$ by Grassmann--Pl\"{u}cker embedding, and the image is cut out by the Grassmann--Pl\"{u}cker relations together with certain linear relations.
Such linear relations capture a property that if $\Lambda = \begin{bmatrix} A \,|\, B \end{bmatrix}$ is an $n$-by-$E$ matrix such that its row-space is Lagrangian, %
then $AB^t$ is a symmetric matrix.
Boege, D'Al\`i, Kahle, and Sturmfels~\cite{BDKS2019} showed that the Lagrangian Grassmannian $\lag{k}{n}$ is parameterized as a subset of much lower dimensional projective space $\bP^{2^{n-2}\left(4+\binom{n}{2}\right)-1}(k)$.
Their parameterization still uses the Grassmann--Pl\"{u}cker embedding but only considers the Pl\"{u}cker coordinates corresponding to principal and almost-principal minors of the symmetric matrix $AB^t$.
The defining ideal $J_n$ of its image is scheme-theoretically generated by its quadratic part $(J_n)_2$, and the quadratic relations in $(J_n)_2$ can be classified into four types~{\cite[Theorem~5]{BDKS2019}}. 
The first type is nothing but a $3$-term Grassmann--Pl\"{u}cker relation, but the remaining three types are less clear to be understood as determinantal identities coming from a symmetric matrix.
We resolve this problem by refining quadratic relations cutting out the Lagrangian Grassmannian into a single type, which corresponds to the Laplace expansions of principal and almost-principal minors of symmetric matrices.

\begin{theorem}\label{thm: parameterization}
    The Lagrangian Grassmannian $\lag{k}{n}$ is parameterized as a subset of the projective space of dimension $2^{n-2}(4+\binom{n}{2})-1$, which is set-theoretically cut out by the following quadrics:
    \begin{align*}
        \sum_{e\in S\setminus T}
        (-1)^{\smaller{S}{e} + \smaller{T}{e}}
        x_{S-e}
        x_{T+e} = 0
        \tag{\ref{eq: r G--P}}
    \end{align*}
    where $S$ and $T$ are subsets of $E$ of sizes $n+1$ and $n-1$, respectively, and $S$ contains exactly one $\{i,i^*\}$ for some $i\in[n]$ and $T$ contains no $\{j,j^*\}$.
\end{theorem}

We call such refined quadrics~\eqref{eq: r G--P} the \emph{restricted Grassmann--Pl\"{u}cker relations}. %
Undefined notations will be clarified in Section~\ref{sec: Preliminaries}.
The parameterization of the Lagrangian Grassmannian and the proof of Theorem~\ref{thm: parameterization} will be provided in the same section.

In Section~\ref{sec: antisymmetric matroids}, we define \emph{antisymmetric matroids}, in terms of bases, by extracting zero and nonzero patterns of the restricted Grassmann--Pl\"{u}cker relations, and we provide an equivalent definition, in terms of circuits, capturing the orthogonality and maximality of a Lagrangian subspace (Theorem~\ref{thm: cryptomorphism}).
Antisymmetric matroids generalize matroids as the Grassmannian $\gr{k}{r}{n}$ %
is embedded into the Lagrangian Grassmannian $\lag{k}{n}$ by mapping $V$ to $V\oplus V^\perp$, and moreover, we show that a matroid is representable over a field $k$ in the usual sense if and only if it is representable over $k$ as an antisymmetric matroid (Proposition~\ref{prop: embedding of matroids to antisym mats}).

There were several approaches for understanding type C analogues of matroids, but we believe that none of them gave satisfactory extensions compared to the type D analogue~\cite{Bouchet1988repre,Wenzel1993b,Wenzel1996b,BJ2023,JK2023}. 
Section~\ref{sec: type C} is devoted to explaining relations between antisymmetric matroids and other combinatorial structures.
The most well-known type C counterpart is delta-matroids~\cite{Bouchet1987sym}, which have many different names and equivalent concepts such as symmetric matroids~\cite{Bouchet1987sym}, pseudomatroids~\cite{CK1988}, metroids~\cite{DH1986}, $2$-matroids~\cite{Bouchet1997}, and Lagrangian (symplectic) matroids~\cite{BGW1998}.
A crucial disadvantage of delta-matroids is the lack of the `strong' basis exchange property; see Example~\ref{eg: no strong basis exchange}.
Note that the support of a Pl\"{u}cker vector is the set of bases of a matroid, and the strong basis exchange axiom for matroids is exactly the Grassmann--Pl\"{u}cker relations over the Krasner hyperfield~$\bK$~\cite{BB2019}.
Hence, to utilize delta-matroids as a combinatorial abstraction of the Lagrangian Grassmannian, 
it is natural to demand a representable delta-matroid as the support of a point parameterizing a Lagrangian subspace into the projective space.
However, it is impossible because the defining ideal of the Lagrangian Grassmannian in the projective space is homogeneous, and some representable delta-matroids do not satisfy the strong basis exchange property.
Another disadvantage of delta-matroids is that fundamental circuits do not necessarily exist.
We point out that the fundamental cocircuits of a linear matroid~$M$ with respect to a fixed basis~$B$ are the supports of the rows of a standard matrix representation of~$M$ with respect to~$B$.
We highlight that antisymmetric matroids resolve these problems.
We also explain, in Sections~\ref{sec: gaussoids} and~\ref{sec: oriented gaussoids}, that antisymmetric matroids are closely related to gaussoids which were introduced to understand behaviors of almost-principal minors of a positive definite symmetric matrix~\cite{LM2007}.

\emph{Matroids with coefficient} were founded by 
Dress and Wenzel~\cite{Dress1986,DW1991,DW1992} offering a unified approach to theories of representable matroids, oriented matroids~\cite{BVSWZ1999}, valuated matroids~\cite{DW1992b}, and ordinary matroids.
It was based on comprehending Grassmann--Pl\"{u}cker relations over \emph{fuzzy ring} and was culminated by Baker and Bowler~\cite{BB2019} introducing more tractable and extensive field-like structures, called \emph{tracts}.
Their approach encompasses partial field representations of matroids~\cite{SW1996} previously not covered by fuzzy rings, and it was generalized to flag matroids~\cite{JL2024} and even delta-matroids~\cite{JK2023}.
We extend Baker--Bowler theory for antisymmetric matroids in Sections~\ref{sec: antisymmetric matroids over tracts} and~\ref{sec: cryptomorphism} by introducing \emph{antisymmetric matroids with coefficients in tracts} in two equivalent ways.

\begin{theorem}\label{thm: bij}
    For a tract $F$, there is a natural bijection between antisymmetric $F$-matroids and antisymmetric $F$-circuit sets.
\end{theorem}

The precise definition of \emph{antisymmetric $F$-matroids} and \emph{antisymmetric $F$-circuit sets} will be presented in Section~\ref{sec: antisymmetric matroids over tracts}.
We briefly note that antisymmetric $F$-matroids generalize points in the projective space which satisfy the restricted Grassmann--Pl\"{u}cker relations and antisymmetric $F$-circuit sets generalize Lagrangian subspaces in the standard symplectic vector space, and thus their equivalence implies Theorem~\ref{thm: parameterization}. %
Furthermore, antisymmetric matroids with coefficients in tracts $F$ encompasses
\begin{enumerate}[label=\rm(\arabic*)]
    \item matroids with coefficients in tracts~\cite{BB2019}, and 
    \item the Lagrangian orthogonal matroids (Coxeter matroids of type D~\cite{BGW2003}) with coefficients in tracts~\cite{JK2023} if $-1=1$ in $F$.
\end{enumerate}
They are also compatible with several other concepts such as 
\begin{enumerate}[label=\rm(\arabic*)]
    \setcounter{enumi}{2}
    \item the symplectic Dressian and isotropic tropical linear spaces~\cite{Rincon2012,BO2023} if $F = \bT$ is the tropical hyperfield, and
    \item oriented gaussoids~\cite{BDKS2019} if $F = \bS$ is the sign hyperfield.
\end{enumerate}
All these connections are addressed through Sections~\ref{sec: type C} and~\ref{sec: examples}--\ref{sec: oriented gaussoids}.
Theorem~\ref{thm: bij} is proved in Section~\ref{sec: cryptomorphism}, which involves a homotopy theorem for graphs associated with antisymmetric matroids (Theorem~\ref{thm: homotopy}).
The homotopy theorem informally states that all cycles in such a graph are generated by short cycles.
We underline that our homotopy theorem implies Maurer's homotopy theorem for matroids~\cite{Maurer1973} and Wenzel's homotopy theorem for even delta-matroids~\cite{Wenzel1995}.

We finally show that antisymmetric matroids extend one of the fundamental results on the Grassmannian and matroids, that is, a point in the $\left( \binom{n}{r}-1 \right)$-dimensional projective space is a solution of the all Grassmann--Pl\"{u}cker relations if and only if it satisfies the $3$-term Grassmann--Pl\"{u}cker relations and its supports form the bases of a rank-$r$ matroid on $n$ elements.
This was implicitly shown by Tutte~\cite{Tutte1958} in terms of chain-group representations of matroids and was generalized to matroids with coefficients over perfect tracts by Baker and Bowler~{\cite[Theorem~3.46]{BB2019}}.
We prove a counterpart of Tutte's theorem for the Lagrangian Grassmannian using antisymmetric matroids.
\begin{theorem}\label{thm: tutte for Lag}
    A point in the projective space of dimension $2^{n-2}(4+\binom{n}{2})-1$ satisfies all restricted Grassmann--Pl\"{u}cker relations if and only if it satisfies the $3$-/$4$-term restricted Grassmann--Pl\"{u}cker relations and its supports form the bases of an antisymmetric matroid.
\end{theorem}

This paper is organized as follows.
Section~\ref{sec: Preliminaries} provides basic notations and recalls fundamental relations between matroids and linear spaces.
Moreover, it explicitly describes the parameterization of the Lagrangian Grassmannian~\cite{BDKS2019} and proves Theorem~\ref{thm: parameterization}.
In Section~\ref{sec: antisymmetric matroids}, we define antisymmetric matroids in two equivalent ways and explore their properties such as fundamental circuits, minors, and representability.
We also show how we understand matroids as a special case of antisymmetric matroids.
Section~\ref{sec: type C} presents connections between antisymmetric matroids and delta-matroids.
In Section~\ref{sec: homotopy}, we prove a homotopy theorem for graphs associated with antisymmetric matroids, which is an extension of Maurer's homotopy theorem for matroids.
In Section~\ref{sec: antisymmetric matroids over tracts}, we define antisymmetric matroids with coefficients in tracts, which notion is defined in two tantamount ways as a generalization of the parameterization of the Lagrangian Grassmannian and Baker--Bowler theory for matroids.
We also handle their relations with even delta-matroids with coefficients, symplectic Dressian, and oriented gaussids.
Theorem~\ref{thm: bij}, the equivalence between two notions of antisymmetric matroids with coefficients, is proved in Section~\ref{sec: cryptomorphism} using the homotopy theorem.
In Section~\ref{sec: tutte}, we prove Theorem~\ref{thm: tutte for Lag}, an analogue of Tutte's theorem for the Lagrangian Grassmannian.
In Section~\ref{sec: concluding remarks}, we finish the paper with several open problems on antisymmetric matroids motivated by matroid theory.

\section{Preliminaries}\label{sec: Preliminaries}

For two sets $S$ and $T$, we often write $S\cup T$ and $S\setminus T$ as $S+T$ and $S-T$, respectively.
If $T=\{x\}$, we abuse the notation and write $S+x$ and $S-x$ rather than $S+\{x\}$ and $S-\{x\}$.
When the symbols `$+$' and `$-$' are used more than once, we read them from left to right, such as $S-x+y = (S-x)+y$.
We often omit brackets and commas while denoting a set, such as $abc = \{a,b,c\}$.

For a set $E$ and an integer $r$, we denote by $\binom{E}{r}$ the set of all $r$-element subsets of $E$.
Suppose $E$ is equipped with a linear ordering $<$.
Then for $S\subseteq E$ and $x\in E$, let $\smaller{S}{x}$ be the number of elements $y\in S$ smaller than $x$.
We similarly define $\smallereq{S}{x}$ as the number of elements $y\in S$ smaller than or equal to $x$.
We usually denote a field by $k$.
The~\emph{support} of a vector $X\in k^E$, denoted by $\supp(X)$ or $\ul{X}$, is the set of $i\in E$ such that $X(i)\ne 0$.

\paragraph{Matroids} %
A \emph{matroid} is a pair $M = (E,\cB)$ of a finite set $E$ and a nonempty set $\cB$ of subsets of $E$ satisfying the \emph{(strong) basis exchange axiom}:
\begin{itemize}
    \item For all $B,B'\in\cB$ and $e\in B\setminus B'$, there is $f\in B'\setminus B$ such that $B-e+f\in\cB$ and $B'+e-f\in\cB$.
\end{itemize}
Each element in $\cB$ is called a \emph{basis} of $M$, and the \emph{rank} of $M$ is the size of a basis.

A matroid is a combinatorial abstraction of linear independence in a vector space, which is revealed clearly from the Grassmann--Pl\"{u}cker relations.
For a field $k$ and integers $0\le r\le n$, the \emph{Grassmannian} $\gr{k}{r}{n}$ is the set of $r$-dimensional vector spaces in the $n$-dimensional vector space $k^n$.
It is parameterized as a subset of the projective space of dimension $\binom{n}{r}-1$ by the Grassmann--Pl\"{u}cker embedding and is exactly the solution of the \emph{Grassmann--Pl\"{u}cker relations}:
\begin{align*}
    \sum_{x\in S\setminus T} (-1)^{\smaller{S}{x} + \smaller{T}{x}} p_{S-x}p_{T+x} = 0
    \text{ for all $S \in \binom{[n]}{r+1}$ and $T \in \binom{[n]}{r-1}$.}
    \tag{$\dagger$}
    \label{eq: G--P}
\end{align*}
For a given point $p$ satisfying all Grassmann--Pl\"{u}cker relations, let $\cB$ be the set of $B\in \binom{[n]}{r}$ such that $p_B \ne 0$.
Then for any $B,B'\in \cB$ and $e\in B\setminus B'$, we have
\[
    \sum_{x\in (B'+e)\setminus (B-e)} (-1)^{\smaller{(B'+e)}{x} + \smaller{(B-e)}{x}} p_{B'+e-x} p_{B-e+x} = 0.
\]
Thus, there is $f\in B'\setminus B$ such that $B'+e-f$ and $B-e+f$ are in $\cB$, implying that $\cB$ is the set of bases of a rank-$r$ matroid on $[n]$.
A matroid is \emph{representable over $k$} if it is obtainable by the previous construction up to isomorphism.

There are more ways to understand a matroid as an underlying combinatorial structure of a linear vector space.
We first recall two cryptomorphic definitions of a matroid.
A \emph{circuit} of a matroid $M=(E,\cB)$ is a minimal subset of $E$ contained in no basis.
The \emph{dual} of $M$ is a matroid $M^\perp := (E,\cB^\perp)$ where $\cB^\perp := \{E\setminus B : B\in \cB\}$.
A \emph{cocircuit} of $M$ is a circuit of the dual $M^\perp$.
A \emph{clutter} is a set $\cC$ of subsets of a finite set such that no element in $\cC$ is a proper subset of another element in $\cC$.
A clutter is \emph{nontrivial} if it does not contain the empty set.

\begin{lemma}[see~\cite{Oxley2011matroid}]\label{lem: mat circuits}
    Let $\cC$ be a set of subsets of $E$. 
    Then $\cC$ is the set of circuits of a matroid if and only if $\cC$ is a nontrivial clutter and satisfies the \emph{circuit elimination axiom}:
    \begin{itemize}
        \item For distinct $C,C'\in\cC$ and $e\in C\cap C'$, there is $C''\in \cC$ such that $C'' \subseteq (C\cup C') - e$.
    \end{itemize}
\end{lemma}

\begin{lemma}[Minty's Painting Axiom~\cite{Minty1966}]\label{lem: mat Minty} %
    Let $\cC$ and $\cD$ be sets of subsets of $E$. 
    Then $\cC$ is the set of circuits of a matroid and $\cD$ is the set of cocircuits of the same matroid if and only if $\cC$ and~$\cD$ satisfy the following:
    \begin{enumerate}[label=\rm(\roman*)]
        \item $\cC$ and $\cD$ are nontrivial clutters.
        \item\label{item: mat orth} $|C\cap D| \ne 1$ for all $C\in\cC$ and $D\in\cD$.
        \item\label{item: mat max} For every tripartition $(P,Q,\{e\})$ of $E$, either there is $C\in\cC$ such that $e\in C\subseteq P+e$ or there is $D\in\cD$ such that $e\in D\subseteq Q+e$.
    \end{enumerate}
\end{lemma}

We often call Lemma~\ref{lem: mat Minty}\ref{item: mat orth} the \emph{orthogonality} of matroids.

Let $V$ be an $r$-dimensional linear space in $k^n = k^{[n]}$, and let $\cC$ be the set of minimal supports of nonzero vectors in~$V$.
Then $\cC$ is the set of circuits of a matroid $M$, because for any $X,Y\in \cC$ with $X(e)=Y(e)\ne 0$, $\supp(X-Y) \subseteq (\supp(X) + \supp(Y)) - e$.
We note that the rank of $M$ is $n-r$ and $\cC$ satisfies the \emph{strong circuit elimination axiom}: For all $C,C'\in\cC$, $e\in C\cap C'$, and $f\in C\setminus C'$, there is $C''\in \cC$ such that $f\in C'' \subseteq (C\cup C') - e$.
We further remark that if we let $N$ be the matroid induced by the supports of $p(V)$, where $p$ is the Grassmann--Pl\"{u}cker embedding, then $N=M^\perp$.
The linear space $V$ also induces the sets of circuits and cocircuits of a matroid as follows.

\begin{lemma}[folklore]
    Let $V$ be a linear space in $k^n$.
    Let $\cC$ be the set of minimal supports of nonzero vectors in $V$, and 
    let $\cD$ be the set of minimal supports of nonzero vectors in the orthogonal complement~$V^\perp$.
    Then $\cC$ and $\cD$ are the sets of circuits and cocircuits, respectively, of a matroid.
\end{lemma}
\begin{proof}
    It suffices to show that $\cC$ and $\cD$ fulfill~\ref{item: mat orth} and~\ref{item: mat max} in Lemma~\ref{lem: mat Minty}.
    Because $\sum_{i=1}^n X(i)Y(i) = 0$ for each $X\in V$ and $Y\in V^\perp$, two sets $\cC$ and $\cD$ satisfies~\ref{item: mat orth}, the orthogonality. 

    Suppose to the contrary that $(P,Q,\{e\})$ is a tripartition of $[n]$ violating~\ref{item: mat max}, i.e., there is no $C\in \cC$ such that $e\in C\subseteq P+e$ and there is no $D\in \cD$ such that $e\in D\subseteq Q+e$.
    Let $P'$ a minimal subset of $P$ such that $P\setminus P'$ does not contain any $C\in\cC$, and let $Q'$ be a minimal subset of $Q$ such that $Q\setminus Q'$ does not contain any $D\in \cD$.
    We denote by $P'' := P-P'+Q'$ and $Q'' := Q-Q'+P'$.
    Then $(P'',Q'',\{e\})$ is a tripartition of $[n]$.
    Note that for each $x\in Q'$, there is $D_x\in\cD$ such that $D_x\subseteq Q$ and $D_x\cap Q' = \{x\}$.
    If $P''+e$ contains some $C\in \cC$, then $C\cap Q' \ne \emptyset$ and thus $C\cap D_x = \{x\}$ for $x\in C\cap Q'$, contradicting the orthogonality.
    Hence $P''+e$ does not contain any $C\in \cC$, and thus $\dim V \le |[n]-(P''+e)|$.
    Similarly, $Q''+e$ does not contain any $D\in \cD$, and so $\dim V^\perp \le |[n]-(Q''+e)|$.
    Therefore, $\dim V + \dim V^\perp \le n-1$, a contradiction.
\end{proof}

\paragraph{Lagrangian Grassmannian}
Let $E:= \ground{n} = \{1,\ldots,n\} \cup \{1^*,\ldots,n^*\}$ with a linear ordering $1<\cdots<n<1^*<\cdots<n^*$ and let $*$ be the natural involution on $E$ mapping $i\in[n]$ to $i^*$.
A {\em skew pair} is a $2$-element subset of $E$ of the form $\{i,i^*\}$.
Let $\chi:E \to \{0,1\}$ be a map such that $\chi(i) = 0$ and $\chi(i^*)=1$ for $i\in[n]$.
Then the \emph{standard symplectic form} $\omega$ on~$k^E$ is represented as $\sympl{X}{Y} = \sum_{i\in E} (-1)^{\chi(i)}X(i)Y(i^*)$.
Let $\cT_n := \{T \in \binom{E}{n} : |T\cap T^*| = 0\}$ and
$\cA_n := \{A \in \binom{E}{n} : |A\cap A^*| = 1\}$.
Then $|\cT_n| = 2^n$ and $|\cA_n| = n(n-1)2^{n-2}$.
We call each element in $\cT_n$ %
(resp. $\cA_n$) %
a {\em transversal} (resp. an {\em almost-transversal}).
A \emph{subtransversal} is a subset of a transversal.

The \emph{Lagrangian Grassmannian} $\lag{k}{n}$ is the set of all maximal isotropic subspaces of the standard symplectic space~$k^{E}$.
In 1979, De Concini~\cite{DeConcini1979} showed that $\lag{k}{n}$ is parameterized as a subset of the projective space of dimension $\binom{2n}{n}-1$, which is exactly the solution of the Grassmann--Pl\"{u}cker relations together with the linear relations $\sum x_{S + \skewpair{i}} = 0$ for all $S\in \binom{E}{n-2}$, where the sum is over all skew pairs $\skewpair{i}$ not intersecting with $S$; see~{\cite[Proposition~2.1]{BO2023}}.\footnote{These linear relations in~\cite{BO2023} seem to miss some signs because $x_{121^*} = x_{233^*}$ for $n=3$.}
Boege, D'Al\`i, Kahle, and Sturmfels~{\cite{BDKS2019}} proved that $\lag{k}{n}$ is also parameterized as a subset of the projective space of dimension $2^{n-2}(4+\binom{n}{2})-1$ whose image is cut out by the quadratic relations of its defining ideal. 
These quadratic relations are classified into four types according to Theorem~5 in~\cite{BDKS2019}.

Now we review the latter parameterization of $\lag{k}{n}$ by~\cite{BDKS2019} and prove Theorem~\ref{thm: parameterization}.
We regard the coordinates of the projective space of dimension $2^{n-2}(4+\binom{n}{2})-1$ as elements in $\cT_n$ and $\cA_n$ along with the identification $(-1)^{i} x_{S+\skewpair{i}} = (-1)^{j} x_{S+\skewpair{j}}$ for each subtransversal $S$ of size $n-2$ and distinct elements $i,j \in [n] \setminus (S\cup S^*)$.
Let $\Lambda = \begin{bmatrix} \Lambda_1 \,|\, \Lambda_2 \end{bmatrix}$ be an $n\times E$ matrix whose row-space is a Lagrangian subspace $W$ in $k^E$ or, equivalently, $\Lambda_1 \Lambda_2^t$ is a symmetric matrix.
We parameterize $W$ into the projective space as $\Phi(W) := ( \det(\Lambda[n,B] ) )_{B \in \cT_n \cup \cA_n}$.
Then the point $x=\Phi(W)$ evidently satisfies the following \emph{restricted Grassmann--Pl\"{u}cker relations} (in short, \emph{restricted G--P relations}):
\begin{align*}
    \sum_{b \in S_1 \setminus S_2} (-1)^{\smaller{S_1}{b} + \smaller{S_2}{b}} x_{S_1-b}x_{S_2+b} = 0
    \tag{$\ddagger$}
    \label{eq: r G--P}
\end{align*}
for all $S_1 := T_1+a_1$ and $S_2 := T_2-a_2$ such that $T_1,T_2\in \cT_n$ and $a_1,a_2 \in T_2\setminus T_1$ (possibly, $a_1=a_2$).
The simplest restricted G--P relations consist of three terms.
The first example is
\[
    x_{Sab}x_{Sa^*b^*} + x_{Sab^*}x_{Sa^*b} - x_{Saa^*}x_{Sbb^*} = 0
\]
where $Sab = S + \{a,b\}$ is a transversal with $a,b\in [n]$, which is a restricted G--P relation~\eqref{eq: r G--P} applied to $S_1 = S + \{a,a^*,b^*\}$ and $S_2 = S + \{b\}$. 
It contains a square term $x_{Saa^*}x_{Sbb^*} = (-1)^{a+b} x_{Saa^*}^2$ and thus we call such relations \emph{square relations}.
The second example is 
\[
    (-1)^{\smaller{L}{a}} x_{S abc} x_{S bb^*c^*}
    +
    (-1)^{\smaller{L}{b^*}}
    x_{S abc^*} x_{S bb^*c}
    +
    (-1)^{\smaller{L}{c^*}}
    x_{S abb^*} x_{S bcc^*}
    =0 
\]
where $Sabc = S + \{a,b,c\}$ is a transversal and $L = \{a,c,b^*,c^*\}$, which is a restricted G--P relation applied to $S_1 = S + \{a,b,b^*,c^*\}$ and $S_2 = S + \{b,c\}$.
We call these relations \emph{edge relations}.\footnote{The corresponding trinomials of square and edge relations are respectively called \emph{square} and \emph{edge} trinomials in~\cite{BDKS2019}.}

The restricted G--P relations can be identified with the Laplace expansion of the symmetric matrix $\Sigma := \Lambda_1 \Lambda_2^t$ as follows.
For simplicity, we assume that $\Lambda_1$ is the identity matrix.
Then for every $X,Y\subseteq [n]$ such that $|X|=|Y|$ and $|X\setminus Y| \le 1$, we have
\begin{align*}
    \det(\Sigma[X,Y])
    =
    (-1)^m
    \cdot 
    x_{[n] - X + Y^*}
\end{align*}
where $m:= \sum_{i\in[n]-X} \big( i+\smallereq{([n]-X)}{i} \big)$.
Thus, 
for $U\subseteq V \subseteq [n]$ and $i,j \in [n]\setminus V$, 
the restricted G--P relation~\eqref{eq: r G--P} applied to $S_1 = [n]-V+(V+j)^*$ and $S_2 = [n]-(U+i)+U^*$ is identified with the Laplace expansion 
\begin{align*}
    &\det(\Sigma[U,U])
    \det(\Sigma[V+i,V+j]) \\
    &\hspace{3cm}=
    \sum_{j'\in (V+j)\setminus U} 
    (-1)^{\smaller{V\setminus U}{i}+\smaller{(V+j)\setminus U}{j'}}
    \det(\Sigma[U+i,U+j'])
    \det(\Sigma[V,V+j-j']).
\end{align*}

We finally prove Theorem~\ref{thm: parameterization} stating that $\Phi$ is a parameterization of the Lagrangian Grassmannian whose image is set-theoretically cut out by the restricted G--P relations.
For $i,j\in \ground{n}$, let $1_{i<j}$ be the indicator for the inequality $i<j$, i.e., it is $1$ if $i<j$ and $0$ otherwise.
Then $\chi(i) = 1_{i^*<i}$.

\begin{proof}[\mbox{\bf Proof of Theorem~\ref{thm: parameterization}}]
    It suffices to show that for each point $x$ satisfying all restricted G--P relations, there is a Lagrangian subspace $W$ such that $\Phi(W) = x$.
    By the square relations, there is a transversal $T \in \cT_n$ such that $x_{T} \ne 0$.
    For each $i\in T$, let $X_i$ be a vector in $k^E$ such that $\supp(X_i) \subseteq T^* + i$ and 
    $X_i(j) = (-1)^{\smaller{T}{i}+\smaller{(T-i)}{j}} \frac{x_{T-i+j}}{x_{T}}$ 
    for each $j\in T^*+i$.
    As $X_i(i) = 1$, the $n$ vectors $X_i$'s are linearly independent.
    
    We claim that the span of $\{X_i : i\in T\}$, say $W$, is Lagrangian.
    It suffices to check that $\sympl{X_i}{X_j} = 0$ for all $i,j\in T$.
    Clearly, $\omega(X_i,X_i) = (-1)^{\chi(i)}X_i(i)X_i(i^*) + (-1)^{\chi(i^*)}X_i(i^*)X_i(i) = 0$ for each $i\in T$.
    For distinct $i,j\in T$, let $U = T-ij$.
    Then we have $\smaller{T}{i} = \smaller{U}{i} + 1_{j<i}$ and $\smaller{(T-i)}{j^*} = \smaller{U}{j^*} + \chi(j^*)$.
    Let $\ul{i} = \begin{cases}
        i & \text{if } i\in[n], \\
        i^* & \text{otherwise}.
    \end{cases}$
    We similarly define $\ul{j}$.
    Since $U$ is a subtransversal of size $n-2$ non-intersecting with $\skewpair{i}$ and $\skewpair{j}$, we deduce that $\sum_{e\in\{i,i^*,j,j^*\}} \smaller{U}{e} \equiv 1+\ul{i}+\ul{j} \pmod{2}$.
    Hence $\smaller{T}{i} + \smaller{(T-i)}{j^*} + \smaller{T}{j} + \smaller{(T-j)}{i^*} 
    \equiv
    \chi(i^*)+\chi(j^*)+ 
    \ul{i}+\ul{j} \pmod{2}$.
    Note that $x_{T-i+j} = (-1)^{\ul{i}+\ul{j}}x_{T-j+i}$.
    Then we have
    \begin{align*}
        X_i(j^*) 
        &= 
        (-1)^{\smaller{T}{i}+\smaller{(T-i)}{j}} \frac{x_{T-i+j}}{x_{T}}
        =
        (-1)^{\smaller{T}{j}+\smaller{(T-j)}{i} + \chi(i^*)+\chi(j^*)} \frac{x_{T-j+i}}{x_{T}} 
        =
        (-1)^{\chi(i^*)+\chi(j^*)} X_j(i^*)
    \end{align*}
    and so $\sympl{X_i}{X_j} = (-1)^{\chi(i)} X_i(i) X_j(i^*) + (-1)^{\chi(j^*)} X_i(j^*) X_j(j) =  0$.
    Thus, $W$ is Lagrangian.

    We finally show that $\Phi(W) = x$.
    Let $\Lambda$ be an $n\times E$-matrix whose rows are $X_i$'s ordered with respect to the linear ordering $1<\cdots<n<1^*<\cdots<n^*$.
    Then $\Lambda[n,T]$ is the identity matrix and for each $i\in T$ and $j\in T^*$, we have
    \begin{align*}
        \det(\Lambda[n,T-i+j])
        = (-1)^{\smaller{T}{i}+\smaller{(T-i)}{j}}
        X_i(j)
        = 
        \frac{x_{T-i+j}}{x_{T}}.
    \end{align*}
    Therefore, $\det(\Lambda[n,B]) = \frac{x_{B}}{x_{T}}$ for every $B\in \cT_n \cup \cA_n$ and thus $\Phi(W) = x$.
\end{proof}

\section{Antisymmetric matroids}\label{sec: antisymmetric matroids}

We introduce a combinatorial abstraction, called an \emph{antisymmetric matroid}, of a Lagrangian subspace in two equivalent ways. %
First, we define this notion in terms of bases so that its typical example is the support of a point satisfying all restricted G--P relations.
It is reminiscent of the relation between the basis exchange axiom of a matroid and the Grassmann--Pl\"{u}cker relations.
Second, in Section~\ref{sec: antisymmetric matroids circuits}, we provide a tantamount definition in terms of circuits, which extends the Minty's Painting Axiom (Lemma~\ref{lem: mat Minty}).
It captures properties of certain minimal supports of nonzero vectors in a Lagrangian subspace.

We define the representability and minors of antisymmetric matroid in Sections~\ref{sec: representability} and~\ref{sec: minors}, respectively.
We also show how antisymmetric matroids encompass matroids through Sections~\ref{sec: antisymmetric matroids circuits}--\ref{sec: representability}; see also Section~\ref{sec: matroids with coefficients}.

Let $E:=\ground{n}$ through this section.

\begin{definition}[Antisymmetric Matroids]\label{def: antisymmetric matroids}
    A pair
    $M = (\ground{n},\cB)$ is an {\em antisymmetric matroid} if $\cB \subseteq \cT_n \cup \cA_n$ and the following hold:
    \begin{enumerate}[label=\rm(B\arabic*)]
        \item\label{item: sB1} $\cB \ne \emptyset$.
        \item\label{item: sB2} 
        For $T \in \cT_n$ and distinct skew pairs $p$ and $q$, $T+p-q \in \cB \cap \cA_n$ if and only if $T-p+q \in \cB \cap \cA_n$.
    \end{enumerate}
    \begin{enumerate}[label=\rm(Exch)]
        \item\label{item: exchange} 
        For $B, B' \in \cB$ and $e\in B \setminus B'$,
        if $B-e$ has no skew pair and $B'+e$ has exactly one skew pair, then there is $f\in B' \setminus B$ such that both $B-e+f$ and $B'+e-f$ are in $\cB$.
    \end{enumerate}
    We call each element in $\cB(M):= \cB$ a {\em basis} of~$M$.
\end{definition}

One can rewrite~\ref{item: exchange} as follows, which captures the zero and nonzero patterns of a point in the projective space of dimension $2^{n-2}(4+\binom{n}{2})$ satisfying all restricted Grassmann--Pl\"{u}cker relations~\eqref{eq: r G--P}.

\begin{enumerate}[label=\rm(Exch$'$)]
    \item\label{item: exchange'}
    For arbitrary transversals $T,T'$ and $e,f \in T' \setminus T$ (possibly, $e=f$)
    there are no or at least two elements $g \in (T+e) \setminus (T'-f)$ such that  $\{T+e-g,T'-f+g\} \subseteq \cB$.
\end{enumerate}

\begin{lemma}\label{lem: 3-term basis exchange}
    Let $M = (E,\cB)$ be an antisymmetric matroid.
    Let $T$ be a transversal and $p$, $q$ be distinct skew pairs.
    Then none or at least two of 
    \[
        \{T+p-q, T-p+q\}, \; \{T, T\symdiff(p + q)\}, \;\text{and}\; \{T\symdiff p, T\symdiff q\}
    \]
    are contained in~$\cB$.
    In particular, if $T+p-q \in \cB \cap \cA_n$, then $\{T, T\symdiff(p + q)\} \subseteq \cB$ or $\{T \symdiff p, T \symdiff q\} \subseteq \cB$.
\end{lemma}
\begin{proof}
    We denote by $\{x\} = T\cap p$ and $\{y\} = T\cap q$.
    By taking $T$, $T'$, $e$, and $f$ in~\ref{item: exchange'} as $T-p+q$, $T+p-q$, $x$, and $x$,
    there is no or at least two $g \in \{x,y,y^*\} = (T-p+q+x) \setminus (T+p-q-x)$ such that $T+p-q-x+g$ and $T-p+q+x-g$ are bases of $M$.
    Note that
    \begin{align*}
        T+p-q-x+g = 
        \begin{cases}
            T+p-q & \text{if } g=x, \\
            T\symdiff(p + q) & \text{if } g=y^*, \\
            T\symdiff p & \text{otherwise},
        \end{cases}
        \;\text{ and }\;
        T-p+q+x-g = 
        \begin{cases}
            T-p+q & \text{if } g=x, \\
            T & \text{if } g=y^*, \\
            T\symdiff q & \text{otherwise}.
        \end{cases}
    \end{align*}
    By~\ref{item: sB2}, $T+p-q\in \cB$ if and only if $T-p+q\in \cB$.
    Hence the proof is completed.
\end{proof}

\subsection{Circuits}\label{sec: antisymmetric matroids circuits}

A {\em circuit} of an antisymmetric matroid $M$ on $E=\ground{n}$ is a minimal subset $C$ of $E$ such that $C$ contains at most one skew pair and $C$ is not a subset of any basis of~$M$.
We denote by $\cC(M)$ the set of circuits of~$M$.
Note that every circuit is nonempty by~\ref{item: sB1}, and $\cB(M)$ is equal to the set of $B\in \cT_n\cup \cA_n$ such that $B$ is not a superset of any circuit of~$M$.
We present a cryptomorphic definition of antisymmetric matroids in terms of circuits, which is reminiscent of Minty's Painting Axiom (Lemma~\ref{lem: mat Minty}).

\begin{theorem}\label{thm: cryptomorphism}
    Let $\cC$ be a set of subsets $C$ of $E$ such that $C$ contains at most one skew pair.
    Then $\cC$ is the set of circuits of an antisymmetric matroid on $E$ if and only if it satisfies the following:
    \begin{enumerate}[label=\rm(C\arabic*)]
        \item\label{item: sC1} $\emptyset \notin \cC$.
        \item\label{item: sC2} If $C_1,C_2 \in \cC$ and $C_1 \subseteq C_2$, then $C_1 = C_2$.
    \end{enumerate}
    \begin{enumerate}[label=\rm(Orth)]
        \item\label{item: isotropic} $|C_1 \cap C_2^*| \ne 1$ for all $C_1,C_2 \in \cC$.
    \end{enumerate}
    \begin{enumerate}[label=\rm(Max)]
        \item\label{item: maximal} For every transversal $T \in \cT_n$ and element $e\in T^*$, there is $C\in\cC$ such that $C\subseteq T \cup \{e\}$.
    \end{enumerate}
\end{theorem}

We prove Theorem~\ref{thm: cryptomorphism} at the end of this subsection.
We note that~\ref{item: maximal} can be replaced with the following stronger condition:
\begin{enumerate}[label=\rm(Max$'$)]
    \item\label{item: maximal'} For every transversal $T\in \cT_n$ and element $e\in T^*$, there is $C\in \cC$ such that $C\subseteq T\cup\{e\}$ and $C\cap \skewpair{e} \ne \emptyset$.
\end{enumerate}
It is due to the next lemma.
\begin{lemma}
    Let $\cC$ be a set of subsets $C$ of $E$ such that $C$ contains at most one skew pair.
    If $\cC$ satisfies~\ref{item: isotropic} and~\ref{item: maximal}, then it satisfies~\ref{item: maximal'}.
\end{lemma}
\begin{proof}
    Let $T\in \cT_n$ and $e\in T^*$.
    Let $S$ be a minimal subset of $T-e^*$ such that $T-e^*-S$ does not contain any $C\in\cC$.
    Then for each $f\in S$, there is $C_f\in \cC$ such that $C_f\subseteq T-e^*$ and $C_f\cap S = \{f\}$.
    Let $T' = T-S+S^*$.
    By~\ref{item: maximal}, there is $D\in \cC$ such that $D\subseteq T'+e$.
    If $D\cap S^* \ne \emptyset$, then $D\cap C_{f}^* = \{f^*\}$ for $f^* \in D \cap S^*$, contradicting~\ref{item: isotropic}.
    Thus, $D\subseteq T'+e-S^* = T+e-S$ and $D\cap \{e,e^*\} \ne \emptyset$.
\end{proof}

The set $\bigsqcup_{0\le r\le n} \gr{k}{r}{n}$ of linear spaces can be embedded into the Lagrangian Grassmannian $\lag{k}{n}$ by mapping a linear space $V$ %
to $V\oplus V^\perp$. %
Analogously, there is a natural injection from the set of matroids on $[n]$ to the set of antisymmetric matroids on $\ground{n}$.
For a matroid $M$ on $[n]$, and let $\cD$ be the set $\cC(M) \oplus \cC(M^\perp) := \cC(M) \cup \{C^*\subseteq [n]^* : C\in \cC(M^\perp)\}$.
Then $\cD$ satisfies \ref{item: sC1},~\ref{item: sC2}, \ref{item: isotropic}, and~\ref{item: maximal} as follows.
As $\cC(M)$ and $\cC(M^\perp)$ are the sets of circuits and cocircuits of $M$, it obviously satisfies~\ref{item: sC1} and~\ref{item: sC2}.
To check~\ref{item: isotropic}, we can assume that $C_1\in \cC(M)$ and $C_2^* \in \cC(M^\perp)$.
Then $|C_1 \cap C_2^*| \ne 1$ is exactly the orthogonality of matroids.
We finally examine~\ref{item: maximal} for a transversal $T$ and an element $e\in T^*$.
We may assume that $e\in [n]$.
Applying Lemma~\ref{lem: mat Minty}(iii) to a tripartition $(P,Q,\{e\}) = \big( T \cap [n], (T^*-e)\cap[n], \{e\} \big)$, the matroid $M$ has a circuit contained in $(T\cap[n])\cup\{e\}$ or has a cocircuit contained in $T^*\cap[n]$.
Then there is $C\in \cD$ that is a subset of $T+e$.
Therefore, by Theorem~\ref{thm: cryptomorphism}, we conclude that $\cD$ is the set of circuits of an antisymmetric matroid, denoted by $\antlift(M)$, on $\ground{n}$.
Moreover, we deduce the following commutative diagram
\begin{center}
    \begin{tikzcd}
        \gr{k}{r}{n} \ar[r] \ar[d] & \lag{k}{n} \ar[d] \\
        \mathrm{Mat}_{r,n} \ar[r] & \mathrm{AntMat}_{n} \\
    \end{tikzcd}
    \hspace{1cm}
    \begin{tikzcd}
        V \ar[r,mapsto] \ar[d,mapsto] & V\oplus V^\perp \ar[d,mapsto] \\
        M \ar[r,mapsto] & \antlift(M) \\
    \end{tikzcd}
\end{center}
\vspace{-0.7cm}
where $M$ is the rank-$r$ matroid on $[n]$ whose set of cocircuits is the set of minimal supports of vectors in $V \setminus \{\mathbf{0}\}$.
Here, $\mathrm{Mat}_{r,n}$ denotes the set of rank-$r$ matroids on $[n]$ and $\mathrm{AntMat}_{n}$ denotes the set of antisymmetric matroids on $\ground{n}$.
The map $\lag{k}{n} \to \mathrm{AntMat}_n$ is defined by sending a Lagrangian subspace $W$ to an antisymmetric matroid $M$ such that $\cC(M)^* = \{C^* : C\in \cC(M)\}$ equals the set of minimal supports $D$ of nonzero vectors in $W$ such that $D$ contains at most one skew pair, of which well-definedness will be checked in Proposition~\ref{prop: Lag subsp and antisym mat}.
In Section~\ref{sec: even delta-matroids}, we further develop this diagram by extending the domain $\mathrm{Mat}_{r,n}$ to the set of even delta-matroids whenever the field $k$ has characteristic two. %

\begin{example}
    Let $V$ be the row-space of a $2$-by-$3$ matrix $A = \begin{bmatrix}
        1 & 0 & 1 \\
        0 & 1 & 1
    \end{bmatrix}$ over an arbitrary field~$k$.
    Then $V$ is $2$-dimenaional and $\{X \subseteq [3] : \det(A[2,X]) \ne 0\} = \{12,13,23\}$ forms the set of bases of a uniform matroid $M=U_{2,3}$.
    The orthogonal complement $V^\perp$ is the span of $\begin{bmatrix} 1 & 1 & -1 \end{bmatrix}$, which induces the dual matroid $M^\perp = U_{1,3}$.
    Note that $\cC(M) = \{123\}$ and $\cC(M^\perp) = \{12,13,23\}$, which are the sets of minimal supports of nonzero vectors in $V^\perp$ and $V$, respectively.
    Then the set of circuits of an antisymmetric matroid $\antlift(M)$ is $\cC(M) \oplus \cC(M^\perp) = \{123, 1^*2^*, 1^*3^*, 2^*3^*\}$ by definition, which equals to the set of minimal supports of nonzero vectors in $V^\perp \oplus V$.
    The bases of an antisymmetric matroid $\antlift(M)$ are $123^*$, $12^*3$, and $1^*23$.
\end{example}

We show several properties of circuits to prove Theorem~\ref{thm: cryptomorphism}.

\begin{lemma}\label{lem: basis from circuit}
    Let $C$ be a circuit of an antisymmetric matroid on $E$ and let $e\in C$.
    If $C-e$ has no skew pair, then there is a basis $B$ such that $B$ is a transversal and $C\setminus B = \{e\}$.
\end{lemma}
\begin{proof}
    Since $C-e$ is not a circuit, there is a basis $B$ such that $C-e \subseteq B$ and $e\notin B$.
    We may assume that $B \in \cA_n$ and let $p$ and $q$ be skew pairs such that $p\subseteq B$ and $q\cap B = \emptyset$.
    Then $p \ne \{e,e^*\}$.
    Let $f \in p \setminus C$.
    By Lemma~\ref{lem: 3-term basis exchange}, for some $g\in q$, $B' := B-f+g$ is a basis.
    Then $B'$ is a transversal and $C\setminus B' = \{e\}$.
\end{proof}

\begin{lemma}\label{lem: fundamental circuit}
    Let $M = (E,\cB)$ be an antisymmetric matroid and let $S$ be a subset of $E$ such that $|S| = n+1$, $S$ has exactly one skew pair, and $S-e\in \cB$ for some $e\in S$.
    Then there is a unique circuit $C$ contained in $S$.
    Moreover, $C = \{e \in S: S-e \in \cB\}$.
\end{lemma}
\begin{proof}
    Let $C := \{e \in S: S-e \in \cB\}$ and $\{z,z^*\} \subseteq S$.
    By the assumption, $C \ne \emptyset$.
    If $C \setminus \{z,z^*\} \ne \emptyset$, then $S-x$ is a basis for $x\in C\setminus \{z,z^*\}$ and by Lemma~\ref{lem: 3-term basis exchange}, $S-z$ or $S-z^*$ is a basis.
    Thus, $C\cap \{z,z^*\} \ne \emptyset$.
    By relabelling we can assume that a transversal $B := S-z$ is a basis.

    We claim that $C$ is a circuit of $M$.
    For every $x\in C$, we have $C-x \subseteq S-x \in \cB$.
    Hence every proper subset of $C$ is not a circuit.
    Therefore, it suffices to check that $C$ is not a subset of any basis.
    Suppose to the contrary that there is a basis~$B'$ containing~$C$.
    Then $z\in C\setminus B \subseteq B'\setminus B$.
    By Lemma~\ref{lem: 3-term basis exchange}, we may assume that $B'-z$ has no skew pair.
    By~\ref{item: exchange}, there is $y \in B\setminus B'$ such that $S-y = B+z-y \in \cB$.
    Then $y\in C$ and thus $C\not\subseteq B'$, a contradiction.
    Therefore, $C$ is a circuit.

    Let $D$ be a circuit of $M$ such that $D\subseteq S$.
    If $e \in S\setminus D$, then $S-e \notin \cB$ and hence $e \notin C$.
    Then $C\subseteq D$.
    This implies that $C$ is a unique circuit contained in $S$.
\end{proof}

For an antisymmetric matroid $M$, a transversal basis $B$, and an element $e\in B^*$, 
the unique circuit contained in $B+e$ is called the {\em fundamental circuit} of $M$ with respect to $B$ and $e$.
Obviously, such a fundamental circuit contains $e$.

\begin{lemma}\label{lem: orthogonality}
    If $C_1$ and $C_2$ are circuits of an antisymmetric matroid, then $|C_1\cap C_2^*| \ne 1$.
\end{lemma}
\begin{proof}
    By Lemma~\ref{lem: basis from circuit}, there are transversal bases $B_1$ and $B_2$ such that $|C_i \setminus B_i| = 1$. 
    For each $i\in\{1,2\}$, let $S_i = B_i \cup C_i$ and
    let $q_i$ be the skew pair contained in $S_i$.
    For each $e \in S_1 \setminus (S_2-q_2)$, by Lemma~\ref{lem: fundamental circuit},
    $e\in C_1$ if and only if $S_1-e\in\cB$.
    Similarly, 
    $e^*\in C_2$ if and only if $S_2-e^*\in\cB$, and the latter condition is equivalent to $S_2-q_2+e \in \cB$ by~\ref{item: sB2}.
    By~\ref{item: exchange'} applied to $S_1$ and $S_2-q_2$, we deduce that $|C_1 \cap C_2^*| \ne 1$.
\end{proof}

\begin{lemma}\label{lem: circuits with two intersection}
    Let $C$ be a circuit of an antisymmetric matroid $M$ and let $e,f\in C$ be distinct elements such that $C-e$ is a subtransversal.
    Then there is a circuit $D$ such that $C\cap D^* = \{e,f\}$.
\end{lemma}
\begin{proof}
    By Lemma~\ref{lem: basis from circuit}, $M$ has a transversal basis $B$ such that $C \setminus B = \{e\}$.
    Let $D$ be the fundamental circuit with respect to $B$ and $f^*$.
    Then $f \in C\cap D^* \subseteq \{e,f\}$ and thus by Lemma~\ref{lem: orthogonality}, $C\cap D^* = \{e,f\}$.
\end{proof}

\begin{lemma}\label{lem: addition}
    Let $\cC$ be a set of subsets $C$ of $E$ such that $C$ contains at most one skew pair.
    If $\cC$ satisfies~\ref{item: isotropic} and~\ref{item: maximal}, then it satisfies
    the following:
    \begin{enumerate}[label=\rm(Add)]
        \item\label{item: addition} For distinct $C_1,C_2 \in \cC$ and $e\in C_1\cap C_2$, if $(C_1 \cup C_2) -e$ contains at most one skew pair, then there is $C_3 \in \cC$ such that $C_3 \subseteq (C_1 \cup C_2) -e$.
    \end{enumerate}
\end{lemma}
\begin{proof}
    We say a subtransversal $I\subseteq E$ is \emph{$\cC$-independent} if $I$ contains no element in $\cC$.
    If a subtransversal $I$ is $\cC$-independent and a skew pair $\skewpair{e}$ does not intersect with $I$, then $I+e$ or $I+e^*$ is $\cC$-independent by~\ref{item: isotropic}.

    Let $C_1,C_2$ be distinct elements in $\cC$ and let $e\in C_1\cap C_2$ such that $J:=(C_1 \cup C_2) \setminus \{e\}$ contains at most one skew pair.
    Suppose to contrary that $J$ does not contain any element in $\cC$.

    We first assume that $J$ is a subtransversal.
    Then $J$ is $\cC$-independent and thus there is a $\cC$-independent transversal $J'$ containing $J$.
    Then $e^* \in J'$.
    As $C_1$ and $C_2$ are distinct, there is $f\in C_2\setminus C_1$.
    By~\ref{item: maximal}, there is $D\in \cC$ such that $f^* \in D \subseteq J'+f^*$.
    By~\ref{item: isotropic} applied to $D$ and $C_2$, we have $e^*\in D$.
    Then $D\cap C_1^* = \{e^*\}$ contradicting~\ref{item: isotropic}.

    Now we may assume that $J$ has a skew pair, say $\skewpair{f}$. Note that $\skewpair{e} \ne \skewpair{f}$ because $e\notin J$. We consider two cases depending on whether $e^*$ is in $J$ or not.
    
    \smallskip
    
    \noindent\textbf{Case I.} Suppose that $e^*\in J$. By interchanging $C_1$ and $C_2$ if necessary, we may assume that $\skewpair{e} \subseteq C_1$. By interchanging $f$ and $f^*$ if necessary, we may assume that $f^*\notin C_1$. Then $f^*\in C_2$. Because $e\in C_2$, we deduce that $f\in C_1$ by~\ref{item: isotropic}. In short, $\{e,e^*,f\} \subseteq C_1$ and $\{e,f^*\} \subseteq C_2$. Let $K := J-f^*$. Then $K$ is a $\cC$-independent subtransversal and thus there is a $\cC$-independent transversal $K'$ containing $K$. By~\ref{item: maximal}, we have an element $D\in \cC$ such that $f^* \in D \subseteq K'+f^*$. By the assumption that $J$ contains no element in $\cC$, there is $g\in D\setminus J$. By~\ref{item: maximal}, $\cC$ has an element $D'$ such that $g^* \in D' \subseteq K'+g^*$. As $g \in D\cap (D')^* \subseteq \{g,f^*\}$, we deduce that $f^*\in (D')^*$ by~\ref{item: isotropic}. Then $f^* \in C_2 \cap (D')^* \subseteq \{f^*,e\}$ and by~\ref{item: isotropic}, $e\in (D')^*$. Then $C_1 \cap (D')^* = \{e\}$ that contradicts~\ref{item: isotropic}.
    
    \smallskip

    \noindent\textbf{Case II.} Suppose that $e^*\notin J$. We subdivide this case into two subcases.

    \smallskip

    \noindent\textbf{Subcase II-1.} $C_1$ or $C_2$ is a subtransversal. By symmetry, we may assume that $C_1$ is a subtransversal. Then $\skewpair{f} \cap C_1 = \emptyset$ and $\skewpair{f} \subseteq C_2$ by~\ref{item: isotropic}. The remainder is basically the same as the Case~I. Let $K := J-f^*$. Then $K$ is a $\cC$-independent subtransversal and thus there is a $\cC$-independent transversal $K'$ containing $K$. Then $e^* \in K'$. By~\ref{item: maximal}, we have an element $D\in \cC$ such that $f^* \in D \subseteq K'+f^*$. Note that $e^* \notin D$ by~\ref{item: isotropic} applied to $D$ and $C_1$. Because $J$ contains no element in $\cC$, there is $g\in D\setminus J$. By~\ref{item: maximal}, $\cC$ has an element $D'$ such that $g^* \in D' \subseteq K'+g^*$. As $g \in D\cap (D')^* \subseteq \{g,f^*\}$, we deduce that $f^*\in (D')^*$ by~\ref{item: isotropic}. Then $f^* \in C_2 \cap (D')^* \subseteq \{f^*,e\}$ and by~\ref{item: isotropic}, $e\in (D')^*$. Then $C_1 \cap (D')^* = \{e\}$ that contradicts~\ref{item: isotropic}.

    \smallskip

    \noindent\textbf{Subcase II-2.} Neither of $C_1$ and $C_2$ is a subtransversal. Then $\skewpair{f}$ is a subset of $C_1$ and $C_2$. Hence, $(C_1\cup C_2) \setminus \{f\}$ is a subtransversal, which implies that there is $C_3 \in \cC$ such that $C_3 \subseteq (C_1\cup C_2) \setminus \{f\}$ by a part of~\ref{item: addition} covered already. If $e \notin C_3$, then $C_3 \subseteq J$ and we are done. Thus, we may assume that $e\in C_3$. Then by applying Subcase II-1 to $C_1$ and $C_3$, there is $C_4 \in \cC$ such that $C_4 \subseteq (C_1\cup C_3) \setminus \{e\} \subseteq J$.
\end{proof}

We remark that~\ref{item: addition} extends the circuit elimination axiom for matroids.

\begin{proof}[\mbox{\bf Proof of Theorem~\ref{thm: cryptomorphism}}]
    The forward direction is done by Lemmas~\ref{lem: fundamental circuit} and~\ref{lem: orthogonality}.

    Now to show the converse we assume that $\cC$ satisfies the four clauses~\ref{item: sC1},~\ref{item: sC2},~\ref{item: isotropic}, and~\ref{item: maximal}.
    By Lemma~\ref{lem: addition}, $\cC$ also satisfies~\ref{item: addition}.
    Let $\cB$ be the set of transversals and almost-transversals that do not contain any $C\in \cC$.
    It is enough to prove that $\cB$ satisfies \ref{item: sB1},~\ref{item: sB2}, and~\ref{item: exchange'}.

    We first show~\ref{item: sB1}, i.e., $\cB\ne \emptyset$.
    Let $C \in \cC$.
    By~\ref{item: sC1}, $C \ne \emptyset$.
    We choose an element $e\in C$, and
    we additionally assume that $\{e,e^*\}\subseteq C$ if $C$ contains a skew pair.
    Let $I_0 = C - e + e^*$.
    Then by~\ref{item: isotropic}, there is no $D\in\cC$ contained in $I_0$.
    By~\ref{item: isotropic}, if $I$ is a subtransversal containing no set in $\cC$, then for each $\{f,f^*\} \subseteq E-I$, at least one of two sets $I+f$ and $I+f^*$ contains no set in $\cC$.
    Hence we can obtain $B\in \cB \cap \cT_n$ such that $B \supseteq I_0$.

    Second we claim~\ref{item: sB2}.
    Let $T$ be a transversal and $p,q$ be distinct skew pairs such that $T+p-q \in \cB$.
    Suppose to the contrary that there is $C\in \cC$ such that $C \subseteq T-p+q$.
    Then $C \cap q \ne \emptyset$ and let $x \in C\cap q$.
    Replacing $T$ with $T\symdiff q$ if necessary, we can assume that $x\in T$.
    There is $D \in \cC$ such that $D \subseteq (T\symdiff q) + p$ by~\ref{item: maximal}.
    Then $x^* \in D$ because otherwise $D \subseteq T+p-q$.
    Then $C\cap D^* = \{x\}$, contradicting~\ref{item: isotropic}.

    Finally, we show~\ref{item: exchange'}.
    Let $T$, $T'$ be transversals and $e$, $f$ be elements in $T'\setminus T$.
    Let $S := T+e$, $S' := T'-f$, and $q := \{f,f^*\}$.
    We claim that there are no or at least two elements $g\in S\setminus S'$ such that both $S-g$ and $S'+g$ are in $\cB$.
    We may assume that $S'$ does not contain any set in $\cC$, since otherwise $S'+x \notin \cB$ for any $x\in E$.
    Then by~\ref{item: addition} and~\ref{item: isotropic}, 
    $\cC$ has a unique element $D$ such that $D\subseteq S'+q$. %
    By~\ref{item: maximal}, there is $C\in\cC$ such that $C \subseteq S$.
    If there is another $C_2\in\cC$ such that $C_2\subseteq S$, then by~\ref{item: addition}, we deduce that $S-x \notin \cB$ for every $x\in S$.
    Hence we may assume that $C$ is the unique element in $\cC$ such that $C \subseteq S$.
    Then
    \begin{align*}
        C = \{x\in S: S-x \in\cB\}
        \;\text{ and }\;
        D = \{x\in S'+q: S'+q-x \in\cB\}.
    \end{align*}
    For each $x\in S \setminus S'$, we have that $S'+q-x^* \in \cB$ if and only if $S'+x \in \cB$ by~\ref{item: sB2}. 
    Therefore, for each $x\in S\setminus S' = S\cap (S'+q)^*$, $x\in C\cap D^*$ if and only if $S-x \in \cB$ and $S'+x \in \cB$.
    Then the claim follows~\ref{item: isotropic}.
\end{proof}

\subsection{Representability}\label{sec: representability}

Recall that for an $r$-dimensional linear space $V$ in $k^n$, the set of minimal supports of vectors in $V\setminus\{\mathbf{0}\}$ is the set of cocircuits of a rank-$r$ matroid $M$ on $[n]$.
In addition, the set of bases of $M$ is equal to the set of $r$-element subsets $B$ of $[n]$ such that $p(V)_B \ne 0$, where $p$ is the Grassmann--Pl\"{u}cker embedding.
We show an analogous result for Lagrangian subspaces and antisymmetric matroids.
Recall that the parameterization $\Phi$ of the Lagrangian Grassmannian $\lag{k}{n}$ into the projective space of dimension $2^{n-2}(4+\binom{n}{2})-1$ was defined in Section~\ref{sec: Preliminaries}.

\begin{proposition}\label{prop: Lag subsp and antisym mat}
    Let $W$ be a Lagrangian subspace in $k^E$.
    Let $\cB := \{ B\in \cT_n\cup \cA_n : \Phi(W)_B \ne 0 \}$ and
    let $\cC$ be the set of minimal supports $C$ of vectors in $W \setminus \{\mathbf{0}\}$ such that $C$ contains at most one skew pair.
    Then $\cB$ is the set of bases of an antisymmetric matroid, and $\cC$ is the set of circuits of an antisymmetric matroid $(E,\cB^*)$, where $\cB^* := \{B^* : B\in\cB\}$.
\end{proposition}

\begin{proof}%
    The set $\cB$ satisfies~\ref{item: sB1} and~\ref{item: sB2} trivially, and~\ref{item: exchange'} is deduced by the restricted Grassmann--Pl\"{u}cker relations~\eqref{eq: r G--P}.
    Thus, $M_1 = (E,\cB)$ is an antisymmetric matroid.

    We now prove that $\cC$ is the set of circuits of an antisymmetric matroid.
    By definition, \ref{item: sC1} and~\ref{item: sC2} hold.
    As $W$ is isotropic, $\cC$ satisfies~\ref{item: isotropic}.
    Let $T \in \cT_n$ and $e\in T^*$.
    Suppose that there is no $C\in\cC$ such that $C\subseteq T+e$.
    Then $\dim W \le |E - (T+e)| = n-1$, a contradiction.
    Thus,~\ref{item: maximal} holds.
    Then by Theorem~\ref{thm: cryptomorphism}, $\cC$ is the set of circuits of an antisymmetric matroid, say $M_2$, on $E$.

    Finally, we show that $M_2 = (E,\cB^*)$. %
    Let $\Lambda$ be an $n\times E$ matrix such that its row-space is~$W$.
    Then $\Phi(W)_B = \det(\Lambda[n,B])$ for each $B\in \cT_n \cup \cA_n$ by definition.

    \begin{claim}
        For each $B\in \cB$, there is no $C\in \cC$ such that $C\subseteq B^*$.
    \end{claim}
    \begin{proof}
        We denote by $B=\{b_1,\ldots,b_n\}$ and $E-B = \{a_1,\ldots,a_n\}$.
        Since $B\in\cB$, $\Phi(W)_B = \det(\Lambda[n,B])$ is nonzero.
        Hence $\Lambda$ is row-equivalent to a matrix $\Lambda'$ such that $\Lambda'[n,B]$ is an identity matrix.
        So $W$ has $n$ independent vectors $X_1,\ldots,X_n$ such that $\supp(X_i) \cap B = \{b_i\}$.

        Suppose that $B$ is a transversal. %
        As $X_1,\ldots,X_n$ span $W$, every element $C\in \cC$ intersects with $B$.
        Thus, $\cC$ has no element contained in $B^* = E-B$.

        Hence we can assume that $B$ is an almost-transversal.
        By relabelling, we may assume that $b_2 = b_1^*$ and $a_2 = a_1^*$.
        Then $B\cap \{a_1,a_1^*\} = \emptyset = (E-B)\cap \{b_1,b_1^*\}$.
        Hence $0= \sympl{X_1}{X_2} = \sum_{e\in \{b_1,a_1,a_2\}} (-1)^{\chi(e)}X_1(e) X_2(e^*)$.
        Since $X_1(b_1) = 1 = X_2(b_1^*)$, we deduce that $(X_1(a_1),X_1(a_2))$ and $(X_2(a_1),X_2(a_2))$ are not a scalar multiple of each other.
        Then the support of each nonzero linear combination of $X_1,\ldots,X_n$ intersects with $E-B^* = \{a_1,a_2\} \cup \{b_3,\ldots,b_n\}$.
        Therefore, $\cC$ has no element contained in $B^*$.
    \end{proof}
    
    \begin{claim}
        For each pair $(C,e)$ such that $e \in C\in \cC$, there is $B\in \cB$ such that $C-e\subseteq B^*$.
    \end{claim}
    \begin{proof}
        First, suppose that $C-e$ has no skew pair.
        Then $M_2$ has a transversal basis $B = \{b_1,\ldots,b_n\}$ such that $C-e \subseteq B$ by Lemma~\ref{lem: basis from circuit}.
        Let $X_1,\ldots,X_n$ be vectors in $W$ such that each $\supp(X_i)$ is the fundamental circuit of $M_2$ with respect to $B$ and $b_i^*$.
        Because $b_i^* \in \supp(X_i) \subseteq B+b_i^*$, $n$ vectors $X_1,\ldots,X_n$ are independent.
        Hence $\Lambda$ is row-equivalent to a matrix $\Lambda'$ consisting of $X_1,\ldots,X_n$ as rows.
        A submatrix $\Lambda'[n,B^*]$ only has nonzero entries for $(i,b_i^*)$ with $1\le i\le n$.
        Thus, $\Lambda[n,B^*]$ is nonsigular and $B^*\in \cB$.

        Now we assume that $C-e$ has a skew pair, say $\skewpair{x}$.
        Then $M_2$ has a transversal basis $B$ such that $C-x \subseteq B$ by Lemma~\ref{lem: basis from circuit}.
        Then $B_2:=B+x-e$ is also a basis of $M_2$ by Lemma~\ref{lem: fundamental circuit}.
        Note that $B_2 \supseteq C-e$.
        We denote by $B = \{b_1,\ldots,b_n\}$ such that $b_1=x^*$ and $b_2 = e$.
        Let $X_1,\ldots,X_n$ be vectors in $W$ such that each $\supp(X_i)$ is the fundamental circuit of $M_2$ with respect to~$B$ and~$b_i^*$.
        Then $\supp(X_1) = C \ni e = b_2$.
        Then by the orthogonality~\ref{item: isotropic}, $x^* = b_1\in \supp(X_2)$.
        Let $Y_2 := X_2$ and for $i\in [n] \setminus \{2\}$, let $Y_i := X_i - \frac{X_i(b_1)}{X_2(b_1)}$.
        Then an $n\times E$ matrix $\Lambda'$ consisting of $Y_1,\ldots,Y_n$ is row-equivalent to~$\Lambda$ and its square submatrix $\Lambda'[n,\{b_1^*,b_1\}\cup \{b_3^*,\ldots,b_n^*\}] = \Lambda'[n,B_2^*]$ is a nonsingular diagonal matrix.
        Thus, $\Phi(B_2^*) = \det(\Lambda[n,B_2^*]) \ne 0$ and $B_2^*\in \cB$.
        As $C-e\subseteq B_2$, the claim follows.
    \end{proof}

    By the above two claims, $\cC$ is the set of circuits of $(E,\cB^*)$.
\end{proof}

\begin{example}
    Let $\Lambda = 
        \begin{bmatrix}
            1 & 0 & 1 & 1 \\
            0 & 1 & 1 & 1 \\
        \end{bmatrix}$
    be a $2\times (\ground{2})$ matrix over a field $k$.
    Then its row-space, say $W$, is in $\mathrm{Lag}_k(2,4)$.
    Then $\cC := \{12, 11^*2^*, 21^*2^*\}$ is the set of minimal supports $C$ of nonzero vectors in $W$ such that $C$ contains at most one skew pair.
    Let $\cB := \{B\in \cT_2\cup \cA_2: \det(\Lambda[2,B]) \ne 0\} = \{12,11^*,12^*,21^*,22^*\}$.
    Then the set of circuits of an antisymmetric matroid $(E,\cB)$ is $\{1^*2^*, 121^*, 122^*\}= \cC^*$.
\end{example}
\begin{definition}
    An antisymmetric matroid $M$ on $E=\ground{n}$ is \emph{representable} over a field $k$ if there is an $n\times E$ matrix $\Lambda$ over $k$ such that its row-space is Lagrangian in $k^E$ and $\cB(M) = \{B\in \cT_n \cup \cA_n : \det(\Lambda[n,B]) \ne 0 \}$.
\end{definition}

For a Lagrangian subspace $W$ in $k^E$, let $M(W)$ be an antisymmetric matroid on $E$ such that $\cC(M(W))$ is the set of minimal supports $C$ of vectors in $W\setminus\{\mathbf{0}\}$ such that $C$ contains at most one skew pair.
Then an antisymmetric $M$ is representable over $k$ if and only if $M = M(W)$ for some Lagrangian subspace $W$.

In Section~\ref{sec: antisymmetric matroids circuits}, we observed that a linear matroid produces a representable antisymmetric matroid.
More strongly, we show that the representability of antisymmetric matroids extends that of matroids.

\begin{proposition}\label{prop: embedding of matroids to antisym mats}
    A matroid $M$ on $[n]$ is representable over a field $k$ if and only if an antisymmetric matroid $\antlift(M)$ is representable over $k$.
\end{proposition}
\begin{proof}
    Recall that, in Section~\ref{sec: antisymmetric matroids circuits}, $\antlift(M)$ is defined as an antisymmetric matroid whose set of circuits is $\cC(M) \oplus \cC(M^\perp) = \cC(M) \cup \{C^*: C\in \cC(M^\perp)\}$.
    We denote the rank of $M$ by $r$.

    Suppose that $M$ is representable over $k$, and let $V$ be an $r$-dimensional linear subspace in $k^n$ representing $M$ so that the set $\cC(M^\perp)$ of cocircuits of $M$ is exactly the set of minimal supports of nonzero vectors of $V$.
    Then $\cC(\antlift(M)) = \cC(M) \oplus \cC(M^\perp)$ equals the set of minimal supports of nonzero vectors in $V^\perp \oplus V$.
    Thus, $\antlift(M) = M(V^\perp \oplus V)$.

    Suppose that $\antlift(M)$ is representable over $k$. 
    Then there is a Lagrangian subspace $W$ in $k^{\ground{n}}$ such that $\cC(\antlift(M))^*$ equals the set of minimal supports $C$ of nonzero vectors of $W$ such that $C$ contains at most one skew pair.
    Let $V$ be the projection of $W$ into the space $k^{[n]}$ regarding the first $n$ coordinates.
    Since $\cC(\antlift(M))^* = \cC(M^\perp) \oplus \cC(M)$, the projection $V$ is a linear space representing~$M$.
\end{proof}

\subsection{Minors}\label{sec: minors} %

We define minors of an antisymmetric matroid and show that the representability is closed under taking minors.
For a set $\cS$ of subsets, let $\mathrm{Min}(\cS)$ be the set of inclusion-minimal elements of $\cS$.

\begin{proposition}\label{prop: minor}
    Let $M$ be an antisymmetric matroid on $E$ and let $i\in E$.
    We define $\cB(M)|i := \{B-i \subseteq E\setminus \skewpair{i} : B \in \cB(M), \; B\cap\skewpair{i}=\{i\}\}$ if $M$ has a basis containing $i$, and $\cB(M)|i := \{B-i^* \subseteq E\setminus \skewpair{i} : B \in \cB(M) \}$ otherwise.
    We define $\cC(M)|i := \mathrm{Min} \{C - i : i^* \notin C \in \cC(M) \text{ and } C \ne \{i\}\}$.
    Then 
    \begin{enumerate}[label=\rm(\roman*)]
        \item $\cB(M)|i$ is the set of bases of an antisymmetric matroid on $E\setminus\skewpair{i}$, and 
        \item $\cC(M)|i$ is the set of circuits of the same antisymmetric matroid.
    \end{enumerate}
\end{proposition}

We denote the resulting antisymmetric matroid by $M|i$ and call it an {\em elementary minor} of $M$.
An antisymmetric matroid $N$ is a {\em minor} of another antisymmetric matroid $M$ if $N = M|i_1|i_2 \cdots|i_k$ for some $i_1,\ldots,i_k$.
For a matroid $N$ on $[n]$ and $i\in[n]$, we note that $\cC(\antlift(N))|i = \cC(N/i) \oplus \cC(N^\perp \setminus i)$.
Equivalently, $\antlift(N)|i = \antlift(N/i)$.
We similarly have that $\cC(\antlift(N))|i^* = \cC(N\setminus i) \oplus \cC(N^\perp / i)$ and $\antlift(N)|i^* = \antlift(N\setminus i)$.

\begin{example}
    Let $N = U_{3,4}$ be the uniform matroid on $[4]$ of rank $3$.
    Then $\cC(\antlift(N)) = \{1234\} \cup \{i^*j^*: ij\in \binom{[4]}{2}\}$. 
    Note that $N/4 = U_{2,3}$ and $N\setminus 4 = U_{3,3}$.
    Thus, $\cC(\antlift(N))|4 = \{123\} \cup \{i^*j^*: ij\in \binom{[3]}{2}\} = \cC(N/4) \oplus \cC(N^\perp \setminus 4)$
    and 
    $\cC(\antlift(N))|4^* = \{1^*,2^*,3^*\} = \cC(N\setminus 4) \oplus \cC(N^\perp / 4)$.
\end{example}

\begin{proof}[\mbox{\bf Proof of Proposition~\ref{prop: minor}}]
    By definition, $\cB(M)|i$ trivially satisfies~\ref{item: sB2} and~\ref{item: exchange}.
    Now we check~\ref{item: sB1}.
    Suppose that $M$ has a basis $B$ containing $i$.
    If $B\cap\skewpair{i} = \{i\}$, then $\cB(M)|i \ne \emptyset$.
    Therefore, we may assume that $B$ is an almost-transversal containing $\skewpair{i}$.
    Let $\skewpair{j}$ be a skew pair non-intersecting with $B$.
    Then by Lemma~\ref{lem: 3-term basis exchange}, $B-i^*+j$ or $B-i^*+j^*$ is a basis of $M$.
    Thus, $\cB(M)|i \ne \emptyset$.
    Hence we can assume that $M$ has no basis containing $i$.
    If $M$ has a basis which does not include $i^*$, then by~\ref{item: sB2}, $M$ has a basis containing $\skewpair{i}$, a contradiction.
    Thus, every basis of $M$ contains $i^*$, implying that $\cB(M)|i \ne \emptyset$.
    Therefore, $\cB(M)|i$ is the set of bases of an antisymmetric matroid on $E\setminus \skewpair{i}$.

    By definition, $\cC(M)|i$ satisfies~\ref{item: sC1} and~\ref{item: sC2}.
    Let $C_1,C_2\in \cC(M)|i$ and let $D_1,D_2\in \cC(M)$ such that $C_a = D_a - i$ with $a\in\{1,2\}$.
    Then $|C_1 \cap C_2^*| = |D_1 \cap D_2^*| \ne 1$, so $\cC(M)|i$ satisfies~\ref{item: isotropic}.
    Let $T$ be a transversal in $E\setminus\skewpair{i}$ and let $j\in T^*$.
    Suppose that $\{i\}$ is not a circuit of $M$.
    By~\ref{item: maximal}, $M$ has a circuit $D$ contained in $(T+i)+j$.
    Then $D \ne \{i\}$.
    Hence there is $C\in \cC(M)|i$ such that $C \subseteq D-i \subseteq T+j$.
    So we may assume that $\{i\}$ is a circuit of $M$.
    Then every circuit of $M$ does not contain $i^*$ by~\ref{item: isotropic}.
    By~\ref{item: maximal}, $M$ has a circuit $D'$ contained in $(T+i^*)+j$.
    Then $D'\subseteq T+j$.
    Thus, $\cC(M)|i$ satisfies $\ref{item: maximal}$.
    By Theorem~\ref{thm: cryptomorphism}, $\cC(M)|i$ is the set of circuits of an antisymmetric matroid on $E\setminus \skewpair{i}$.

    It remains to show that $\cB(M)|i$ and $\cC(M)|i$ are associated with the same antisymmetric matroid.
    We first assume that $M$ has a basis containing $i$.
    Then $\{i\}$ is not a circuit of $M$.
    Let $\cA$ be the set of subsets $C$ of $E\setminus\skewpair{i}$ such that $C$ contains at most one skew pair.
    We note that if $D\in\cA$ and $D+i$ is a subset of some almost-transversal basis $B$ of $M$ such that $B\supseteq \skewpair{i}$, then by Lemma~\ref{lem: 3-term basis exchange}, $D+i$ is a subset of a transversal basis of $M$.
    Then
    \begin{align*}
        \cC(M|i)
        &=
        \mathrm{Min}\{ C-i : i^* \notin C \in \cC(M) \} \\
        &=
        \mathrm{Min}\{ D\in \cA : D \text{ or } D+i \in \cC(M) \} \\
        &=
        \mathrm{Min}\{ D\in \cA : D+i \not\subseteq B \text{ for all $B\in\cB(M)$ with } i\in B \} \\
        &=
        \mathrm{Min}\{ D\in \cA : D+i \not\subseteq B \text{ for all $B\in\cB(M)$ with } B\cap\skewpair{i} = \{i\} \} \\
        &=
        \mathrm{Min}\{ D\in \cA : D \not\subseteq B' \text{ for all $B'\in \cB(M)|i$} \}. 
    \end{align*}
    Thus, $C(M|i)$ is the set of circuits of the antisymmetric matroid $(E\setminus\skewpair{i},\cB(M)|i)$.

    Next, we suppose that $M$ has no basis containing $i$.
    Then every basis of $M$ contains $i^*$.
    Hence $\{i\}$ is a circuit of $M$, and every circuit of $M$ does not contain $i^*$ by~\ref{item: isotropic}.
    Thus, $\cC(M|i)=\{ C : \{i\} \ne C \in \cC(M)\}$ and it is equal to the set of circuits of the antisymmetric matroid with the basis set $\cB(M)|i = \{B-i^*:B\in \cB(M)\}$.
\end{proof}

The class of antisymmetric matroids representable over a given field is closed under taking minors by the following proposition.

\begin{proposition}\label{prop: repre and minors}
    Let $W \in \lag{k}{n}$ and let $\pi : k^{\ground{n}} \to k^{\ground{n}\setminus\skewpair{i}}$ be the natural projection.
    Then
    $M(\pi(W\cap \mathbf{e}_i^\perp)) = M(W)|i$.
\end{proposition}

\begin{lemma}%
    Let $W$ be in $\lag{k}{n}$ and $v$ be a nonzero vector in $k^{\ground{n}}$ such that $\supp(v) \subseteq \skewpair{n}$.
    Then the natural projection of $W \cap v^\perp$ into $k^{\ground{n-1}}$ is in $\mathrm{SpGr}_{k}(n-1)(2n-2)$.
\end{lemma}
\begin{proof}
    We can assume that $v\notin W$.
    Let $w_1,\ldots,w_n$ be a basis of $W$.
    By relabelling, we may assume that $w_n \notin v^\perp$.
    Then for each $1\le i\le n-1$, there is $c_i\in k$ such that $u_i := w_i - c_i w_n \in v^\perp$.
    We denote by $\pi$ the natural projection from $k^{\ground{n}}$ to $k^{\ground{n-1}}$.
    Then $\pi(u_1),\ldots,\pi(u_n)$ is a basis of $\pi(W\cap v^\perp)$ and $\sympl{\pi(u_i)}{\pi(u_j)} = 0$ for all distinct $i,j\in[n-1]$.
    Therefore, $\pi(W\cap v^\perp) \in \mathrm{SpGr}_{k}(n-1)(2n-2)$.
\end{proof}

\begin{proof}[\mbox{\bf Proof of Proposition~\ref{prop: repre and minors}}]
    Straightforwardly, the set of circuits of $M(\pi(W\cap \mathbf{e}_i^\perp))$ equals $\mathrm{Min} \{C \setminus \{i\} : i^* \notin C \in \cC(M(W)) \text{ and } C \ne \{i\}\}$.
\end{proof}

\section{Other combinatorial structures of type C}\label{sec: type C}

We review several combinatorial concepts derived from symmetric matrices and investigate their relations to antisymmetric matroids.
In Section~\ref{sec: delta-matroids}, we discuss delta-matroids which have been understood as the best combinatorial abstractions for symmetric matrices and Lagrangian symplectic subspaces.
We spotlight that delta-matroids cannot be attained from homogeneous quadratic relations cutting out the Lagrangian Grassmannian $\lag{k}{n}$, in contrast to that antisymmetric matroids are acquired from the restricted Grassmann--Pl\"{u}cker relations~\eqref{eq: r G--P}.
In addition, the circuit axiom for delta-matroids does not possess the existence of fundamental circuits.

In Section~\ref{sec: even delta-matroids}, we handle even delta-matroids which are delta-matroids with a parity condition on the bases.
Even delta-matroids are Coxeter matroids of type D~\cite{BGW2003}, which successfully generalize numerous results on matroids and capture combinatorial properties of skew-symmetric matrices and Lagrangian orthogonal subspaces~\cite{Bouchet1988repre,Wenzel1993b,Wenzel1996b,BJ2023,JK2023}.
We show that the set of even delta-matroids can be regarded as a subclass of antisymmetric matroids, and
further discussions will be continued in Section~\ref{sec: orthgonal matroids with coefficients}.
In Section~\ref{sec: gaussoids}, we show that certain antisymmetric matroids induce gaussoids which arise for studying singular almost-principal submatrices of positive definite real symmetric matrices.

\subsection{Delta-matroids and symmetric matroids}\label{sec: delta-matroids}

A \emph{delta-matroid} is a pair $(E,\cB)$ of a finite set $E$ and a nonempty set $\cB$ of subsets of $E$ satisfying the \emph{symmetric exchange axiom}~\cite{Bouchet1987sym}:
\begin{enumerate}[label=\rm(SEA)]
    \item\label{item: sea} For every $B_1,B_2 \in \cB$ and $e \in B_1 \symdiff B_2$, there is $f \in B_1 \symdiff B_2$ (possibly, $e=f$) such that $B_1 \symdiff\{e,f\} \in \cB$.
\end{enumerate}
We call each element of $\cB$ a \emph{basis}.
Obviously, every matroid is exactly a delta-matroid whose bases have the same cardinality.
A delta-matroid and its equivalent notions were extensively studied by independent researchers under various names such as \emph{symmetric matroids}~\cite{Bouchet1987sym}, \emph{$2$-matroids}~\cite{Bouchet1997}, \emph{metroid}~\cite{DH1986}, \emph{pseudomatroids}~\cite{CK1988}, and \emph{Lagrangian (symplectic) matroids}~\cite{BGW1998}.
For compatibility with antisymmetric matroids, we will work on the following equivalent concept of delta-matroids.

\begin{definition}
    A {\em symmetric matroid} is a pair $M = (\ground{n},\cB)$ such that $\emptyset \ne \cB \subseteq \cT_n$ and $\cB$ satisfies the \emph{symmetric exchange axiom}:
    \begin{enumerate}[label=\rm(SEA$'$)]
        \item\label{item: sea'} For all $B_1,B_2 \in \cB$ and $x \in B_1 - B_2$, there is $y \in B_1 - B_2$ such that $B_1 \symdiff\{x,x^*,y,y^*\} \in \cB$.
    \end{enumerate}
    We call each member of $\cB$ a {\em basis}, and we denote by $\cB(M) = \cB$ the set of bases.
    It is straightforward that $N = ([n],\cB)$ is a delta-matroid if and only if $\lift(N) := (\ground{n}, \{B \cup ([n]\setminus B)^* \in \cT_n : B \in \cB\})$, called the \emph{lift} of $N$, is a symmetric matroid.
\end{definition}

A subset $C \subseteq \ground{n}$ is a {\em circuit} of a symmetric matroid $M$ if $C$ is a subtransversal and it is not contained in any bases of $M$.
We denote by $\cC(M)$ the set of circuits of $M$.
Booth, Moreira, and Pinto presented the circuit axiom for symmetric matroids.

\begin{theorem}[{\cite[Theorem~8]{BMP2003}}]\label{thm: circuits of symmetric}
    Let $\cC$ be a set of subsets $C$ of $E = \ground{n}$ such that $C$ is a subtransversal.
    Then $\cC$ is the set of circuits of a symmetric matroid if and only if it satisfies  \ref{item: sC1},~\ref{item: sC2},~\ref{item: isotropic}, and the following:
    \begin{enumerate}[label=\rm(Add$'$)]%
        \item\label{item: addition'} If $C_1,C_2 \in \cC$ such that $C_1 \cup C_2$ is a subtransversal and $e \in C_1 \cap C_2$, then there is $C_3 \in \cC$ such that $C_3 \subseteq (C_1\cup C_2) - e$.
    \end{enumerate}
\end{theorem}

We note that the condition~\ref{item: addition'} is equivalent to~\ref{item: addition} under the assumption that every element in $\cC$ is a subtransversal.

Similar to antisymmetric matroids, a typical example of delta-matroids (or symmetric matroids) is derived from symmetric matrices.
Let $\Sigma$ be an $n\times n$ symmetric matrix.
Then $M = ([n],\{X \subseteq [n]: \det(\Sigma[X,X]) \ne 0\})$ is a delta-matroid~{\cite[(4.1)]{Bouchet1987sym}}.
By Proposition~\ref{prop: Lag subsp and antisym mat}, $M' = (\ground{n}, \{B\in \cT_n \cup \cA_n : \det(\Lambda[n,B])\})$ is an antisymmetric matroid, where $\Lambda := \begin{bmatrix} \Sigma \, | \, I \end{bmatrix}$ of which columns are indexed by $1,\ldots,n,1^*,\ldots,n^*$ in order.
Note that if $X,Y\subseteq [n]$ with $|X|=|Y|$, then $\det(\Sigma[X,Y]) = (-1)^{\sigma} \det(\Lambda[n,[n]-X+Y^*])$ for some $\sigma\in\{0,1\}$.
Thus, the set of bases of $\lift(M)$ is equal to the set of transversal bases of $M'$, i.e., $\cB(\lift(M)) = \cB(M') \cap \cT_n$.
We extend this observation for general antisymmetric matroids as follows.

\begin{proposition}\label{prop: antisym to sym}
    Let $M = (\ground{n}, \cB)$ is an antisymmetric matroid.
    Then $(\ground{n}, \cB \cap \cT_n)$ is a symmetric matroid.
\end{proposition}
\begin{proof}
    By~\ref{item: sB1} and Lemma~\ref{lem: 3-term basis exchange}, $\cB \cap \cT_n \ne \emptyset$.
    Let $B_1, B_2 \in \cB \cap \cT_n$ and let $x\in B_1-B_2$. Then $x^* \in B_2 - B_1$. %
    We may assume that $B_1 - x + x^* \notin \cB\cap \cT_n$.
    Then by~\ref{item: exchange},
    there is $y \in (B_1-B_2)-x$ %
    such that $B_1 + x^* - y \in \cB \cap \cA_n$. 
    By Lemma~\ref{lem: 3-term basis exchange}, $(B_1+x^*-y) -x + y^* = B_1 \symdiff\{x,x^*,y,y^*\} \in \cB\cap \cT_n$.
\end{proof}

A symmetric matroid is \emph{even} if all intersections of its bases with $[n]$ have the same parity.
The converse of Proposition~\ref{prop: antisym to sym} holds for even symmetric matroids (Theorem~\ref{thm: even to antisymmetric}), which is unknown in general.
We remark that two distinct antisymmetric matroids can induce the same non-even symmetric matroid.

\begin{example}\label{eg: two sym mat with same antisym mat}
    Two antisymmetric matroids $M_1 = (\ground{2},\cT_2)$ and $M_2 = (\ground{2},\cT_2\cup\cA_2)$ induce the same symmetric matroid that is identified with $M_1$.
    It is easily seen that the antisymmetric matroid $M_1$ is representable over the binary field, but $M_2$ is not.
    If $A_1 = 
    \begin{bmatrix}
        1 & 0 \\ 0 & 1
    \end{bmatrix}$
    and
    $A_2 =
    \begin{bmatrix}
        1 & 1 \\ 1 & -1
    \end{bmatrix}$ are ternary matrices, then $M_i$ is represented by $\begin{bmatrix} I \,|\, A_i \end{bmatrix}$ and thus both $M_1$ and $M_2$ are representable over the ternary field.
\end{example} 

We point out two disadvantages of symmetric matroids.
\begin{enumerate}[label=\rm(\roman*)]
    \item The `strong' basis exchange does not hold; see Example~\ref{eg: no strong basis exchange}.
    \item Fundamental circuits do not necessarily exist; see Example~\ref{eg: no fundamental circuit}.
\end{enumerate}
The lack of strong basis exchange implies that symmetric matroids are not appropriate for understanding a parameterization of the Lagrangian Grassmannian into a projective space.
For any parameterization of~$\lag{k}{n}$, the quadratic relations cutting out its image, if exist, are homogeneous, and thus if a symmetric matroid is the support of a point parameterizing a Lagrangian subspace, then it should exchange two bases simultaneously.
However, this fails even for representable symmetric matroids; see Example~\ref{eg: no strong basis exchange}.

The deficiency of fundamental circuits makes it difficult to regard the circuits of a symmetric matroid as certain minimal supports of nonzero vectors in a Lagrangian subspace.
To be more precise, we review the role of fundamental (co)circuits of linear matroids.
Let $M$ be a rank-$r$ linear matroid on $[n]$, let $B$ be its basis, and let $A$ be an $r$-by-$n$ matrix representation of $M$.
By relabelling columns, we may assume that the first $r$ columns of $A$ are indexed by elements of $B$, and let $A'$ be the reduced row-echelon form of $A$.
Then the $r$ fundamental cocircuits of $M$ with respect to $B$ are exactly the supports of $r$ rows in $A'$.
One can similarly observe the same phenomenon for a Lagrangian subspace and antisymmetric matroids.
Let $N$ be an antisymmetric matroid and $\Lambda$ be an $n$-by-$\ground{n}$ matrix such that its row-space is Lagrangian and $\cB(N) = \{B\in\cT_n\cup \cA_n : \det(\Lambda[n,B]) \ne 0\}$.
Let $B_0$ be a transversal basis of $N$, and we relabel the columns of $\Lambda$ in order $b_1,\ldots,b_n,b_1^*,\ldots,b_n^*$, where $b_1 < \ldots < b_n$ are elements of $B_0$.
We denote by the supports $C_1,\ldots,C_n$ of the rows of the reduced row-echelon form of $\Lambda$.
Then each $C_i^*$ is the fundamental circuit of $N$ with respect to $B_0$ and $b_i^*$.
We note that $C_i^*$ is a circuit of a symmetric matroid $N':=(\ground{n},\cB(N)\cap \cT_n)$ if $C_i^*$ is a subtransversal.
However, it is possible that none of $C_i^*$ is a subtransversal, and hence the basis $B$ of $N$ may not have any fundamental circuits; see Example~\ref{eg: no fundamental circuit}.

Thus, symmetric matroids are not appropriate for the study of the Lagrangian Grassmannian.
In contrast, the basis axiom of antisymmetric matroids already possesses the `strong' basis exchange property by definition, and the fundamental circuits of antisymmetric matroids exist always by~\ref{item: maximal} (or Lemma~\ref{lem: fundamental circuit}).

\begin{example}\label{eg: no strong basis exchange}
    Let $M = (\ground{3},\{1^*2^*3^*, 12^*3^*, 1^*23^*, 1^*2^*3, 123\})$ be a symmetric matroid, and let $B=1^*2^*3^*$ and $B'=123$ be its bases.
    Then for each $e\in B'\setminus B = [3]$, exactly one of $B\symdiff 1 1^* e e^*$ and $B'\symdiff 1 1^* e e^*$ is a basis.
    Let
    \[
        \Sigma = \begin{bmatrix} 1 & -1 & -1 \\ -1 & 1 & -1 \\ -1 & -1 & 1 \end{bmatrix}
    \] 
    be a symmetric matrix over a field of characteristic not two.
    Then the delta-matroid associated with $\Sigma$ is $([3], \{\emptyset, 1, 2, 3, 123\})$, and $M$ is its lift.
    There is a unique antisymmetric matroid $N$ on $\ground{3}$ such that $\cB(M) = \cB(N) \cap \cT_n$.
    It has twelve almost-transversal bases $1^*22^*$, $1^*33^*$, $2^*11^*$, $2^*33^*$, $3^*11^*$, $3^*22^*$, $122^*$, $133^*$, $211^*$, $233^*$, $311^*$, and $322^*$.
\end{example}

\begin{example}\label{eg: no fundamental circuit}
    Let $M = (\ground{n},\cT_n)$, which is a symmetric matroid and an antisymmetric matroid, simultaneously.
    It is representable over all fields, witnessed by $\Lambda = \begin{bmatrix}
        I_n \, | \, I_n
    \end{bmatrix}.$
    As a symmetric matroid, $M$ has no circuits.
    Thus, for any basis $B$ and an element $e\in B^*$, there is no circuit contained in $B+e$.
    In contrast, as an antisymmetric matroid, $M$ has the fundamental circuit $\skewpair{e}$ with respect to each $B$ and $e$, which coincide with the $n$ supports of the rows of $\Lambda$. 
\end{example}

We remark that symmetric matroids have a convex polytope characterization.
First, we recall a celebrated theorem on matroid basis polytopes by Gelfand, Goresky, MacPherson, and Serganova~\cite{GGMS1987}.
We denote by $\mathbf{e}_i\in \bR^n$ the $i$-th standard basis vector.

\begin{theorem}[\cite{GGMS1987}]
    Let $\cB$ be a set of subsets of $[n]$ and let $P$ be the convex hull of $\mathbf{e}_B = \sum_{i\in B} \mathbf{e}_i$ with $B\in \cB$.
    Then $([n],\cB)$ is a matroid if and only if every $1$-dimensional face of $P$ is a translation of $\mathbf{e}_i-\mathbf{e}_j$ for some $i,j\in [n]$.
\end{theorem}

We set $e_{i^*} := -e_i$ for each $i\in[n]$ and $\mathbf{e}_S := \sum_{i\in S} \mathbf{e}_i = \sum_{i\in S\cap[n]} \mathbf{e}_i - \sum_{j\in S^*\cap[n]} \mathbf{e}_j$ for each subset $S\subseteq \ground{n}$.

\begin{theorem}[Borovik, Gelfand, and White~\cite{BGW1998}]\label{thm: Lagrangian symplectic matroid polytope}
    Let $\cB \subseteq \cT_n$ and let $P$ be the convex hull of $\mathbf{e}_B$ with $B\in \cB$.
    Then $(\ground{n},\cB)$ is a symmetric matroid if and only if every $1$-dimensional face of $P$ is a translation of $2\mathbf{e}_i$, $2(\mathbf{e}_i + \mathbf{e}_j)$, or $2(\mathbf{e}_i - \mathbf{e}_j)$ for some $i,j \in [n]$ with $i\ne j$.
\end{theorem}

We define a convex polytope $P_M$ associated with an antisymmetric matroid $M=(E,\cB)$, by taking $P_M$ as the convex hull of $\mathbf{e}_B$ with $B\in \cB$.
Then it is the same with the polytope associated with a symmetric matroid $(E,\cB\cap \cT_n)$ by the following lemma.
Figure~\ref{fig: no strong basis exchange} illustrates the polytope defined by the symmetric matroid and the antisymmetric matroid in Example~\ref{eg: no strong basis exchange}.

\begin{corollary}
    For each basis $B$ of an antisymmetric matroid $M=(\ground{n},\cB)$, the vector $\mathbf{e}_B$ is a vertex of $P_M$ if and only if $B$ is a transversal.
    In particular, $P_M$ is equal to the polytope associated with a symmetric matroid $(\ground{n},\cB\cap \cT_n)$.
\end{corollary}
\begin{proof}     
    Let $A$ be an almost-transversal basis, let $p$ and $q$ be skew pairs such that $p\subseteq A$ and $q\cap A = \emptyset$, and let $S:= A-p$.
    By Lemma~\ref{lem: 3-term basis exchange}, $M$ has two transversal bases $T_1$ and $T_2$ such that $T_1 = S+ij$ and $T_2 = S+i^*j^*$ for some elements $i\in p$ and $j\in q$.
    Then $\mathbf{e}_A = \mathbf{e}_S$ is lying on the line segment between $\mathbf{e}_{T_1} = \mathbf{e}_S + \mathbf{e}_{ij}$ and $\mathbf{e}_{T_2} = \mathbf{e}_S - \mathbf{e}_{ij}$.
    Therefore, $P_M$ is indeed equal to the convex hull of the vectors $\mathbf{e}_B$ with transversal bases $B \in \cB \cap \cT_n$, i.e., the polytope associated with the symmetric matroid $(\ground{n}, \cB \cap \cT_n)$.
    Then by Theorem~\ref{thm: Lagrangian symplectic matroid polytope}, the vertices of $P_M$ are exactly the vectors $\mathbf{e}_B$ with transversal bases $B$ of $M$.
\end{proof}

\begin{figure}
    \centering
    \tdplotsetmaincoords{70}{40}
    \begin{tikzpicture}
            [tdplot_main_coords,
                front/.style={thick,black},
                back/.style={thick,dashed,black},
                subdivfront/.style={very thin,blue},
                subdivback/.style={very thin,dashed,blue},
                axis/.style={->,very thin,gray},
                ]

        \begin{scope}[scale=1.5]
        \coordinate (aaa) at (-1,-1,-1);
        \coordinate (baa) at (1,-1,-1);
        \coordinate (aba) at (-1,1,-1);
        \coordinate (aab) at (-1,-1,1);
        \coordinate (bbb) at (1,1,1);
        
        \coordinate (b00) at (1,0,0);
        \coordinate (0b0) at (0,1,0);
        \coordinate (00b) at (0,0,1);

        \coordinate (a00) at (-1,0,0);
        \coordinate (0a0) at (0,-1,0);
        \coordinate (00a) at (0,0,-1);

        \fill[color=blue!25] (aaa) -- (baa) -- (aab) -- cycle;
        \fill[color=blue!12] (bbb) -- (baa) -- (aab) -- cycle;

        \draw[front] (aaa) -- (baa) -- (bbb);
        \draw[back] (aaa) -- (aba) -- (bbb);
        \draw[front] (aaa) -- (aab) -- (bbb);
        \draw[back] (baa) -- (aba) -- (aab);
        \draw[front] (aab) -- (baa);

        \node[vo] () at (aaa) {};
        \node[vo] () at (baa) {};
        \node[vo] () at (aba) {};
        \node[vo] () at (aab) {};
        \node[vo] () at (bbb) {};
        \node[vo] () at (b00) {};
        \node[vo] () at (0b0) {};
        \node[vo] () at (00b) {};
        \node[vo] () at (a00) {};
        \node[vo] () at (0a0) {};
        \node[vo] () at (00a) {};

        \node[left] () at (aaa) {\small{$-e_{123}$}};
        \node[below] () at (baa) {\small{$+e_1-e_{23}$}};
        \node[above left] () at (aab) {\small{$+e_3-e_{12}$}};
        \node[above right] () at (bbb) {\small{$+e_{123}$}};
        \node[right] () at (b00) {\small{$+e_1$}};
        \node[above] () at (00b) {\small{$+e_3$}};
        \node[left] () at (0a0) {\small{$-e_2$}};
        \end{scope}
    \end{tikzpicture}
    \caption{The basis polytope of a symmetric matroid $M=(\ground{3},\{1^*2^*3^*, 12^*3^*, 1^*23^*, 1^*2^*3, 123\})$. 
    The six mid-points $\pm \mathbf{e}_i$ on $1$-dimensional faces represent almost-transversals of the antisymmetric matroid $N$ such that $\cB(M) = \cB(N) \cap \cT_n$.}
    \label{fig: no strong basis exchange}
\end{figure}

\begin{remark}
    Matroids naturally arise from not only linear spaces but also graphs.
    Similarly, one important subclass of delta-matroids is defined from \emph{ribbon graphs}, which are identified with graphs $2$-cell embedded in closed surfaces.
    We refer the readers to a comprehensive survey on delta-matroids by Moffatt~\cite{Moffatt2019}.
    As every delta-matroid associated with a ribbon graph is binary, we can define an antisymmetric matroid from a ribbon graph as well in terms of both bases and circuits.
\end{remark}

\subsection{Even symmetric matroids}\label{sec: even delta-matroids} %

The most important subclass of delta-matroids is even delta-matroids.
A delta-matroid is \emph{even} if its every basis has the same parity. 
Accordingly, a symmetric matroid is \emph{even} if each intersection $B\cap[n]$ of a basis $B$ and $[n]$ has the same parity.
It is also called a \emph{Lagrangian orthogonal matroid}~\cite{BGW2003}.
Even symmetric matroids satisfy the strong basis exchange, and the fundamental circuit always exists for a given basis and element outside the basis.

\begin{theorem}[Wenzel~\cite{Wenzel1993}]
    Let $\cB \subseteq \cT_n$.
    Then $(\ground{n},\cB)$ is an even symmetric matroid if and only if for every $B_1,B_2 \in \cB$ and $x\in B_1 - B_2$, there is $y\in (B_1 - B_2) - x$ such that  
    $B_1 \symdiff \{x,x^*,y,y^*\}$, $B_2 \symdiff \{x,x^*,y,y^*\} \in \cB$.
\end{theorem}

\begin{theorem}[\cite{BMP2003}]
    Let $M$ be an even symmetric matroid.
    For each basis $B$ and element $e\in B^*$, there is a unique circuit $C$ of $M$ such that $e\in C \subseteq B\symdiff \skewpair{e}$.
\end{theorem}

Even symmetric matroids are nice combinatorial abstractions for skew-symmetric matrices and  Lagrangian subspaces in $k^{2n} = k^{\ground{n}}$ equipped with the standard symmetric bilinear form $\ortho{X}{Y} = \sum_{i=1}^n (X(i)Y(i^*) + X(i^*)Y(i))$.
Suppose that our base field $k$ is not characteristic two.
For an $n\times n$ skew-symmetric matrix $\Sigma$, a pair $([n], \{X\subseteq [n] : \pf(\Sigma[X,X]) \ne 0\})$ is an even delta-matroid~\cite{Bouchet1988repre}, where $\pf(\Sigma)$ means the pfaffian of $\Sigma$.
We note that $\det(\Sigma) = \pf(\Sigma)^2$ and the pfaffian of any $(2m+1)\times (2m+1)$ skew-symmetric is zero.
The skew-symmetric matrix $\Sigma$ can be parameterized as a subset of the projective space of dimension $2^{n}-1$ by sending it to a point $(\pf(\Sigma[X,X]): X\subseteq [n])$.
More generally, every maximal isotropic subspace of $k^{2n}$ with respect to $\ortho{\cdot}{\cdot}$ is parameterized as a subset of the projective space, of which image is cut out by homogeneous quadratic relations, called \emph{Wick relations}.
These relations stand for the pfaffian identities analogous to the Laplace expansions of determinants.

Even symmetric matroids with coefficients were successfully generalized by~\cite{Wenzel1996b,Rincon2012,JK2023}, which will be discussed more in Section~\ref{sec: orthgonal matroids with coefficients} by comparing with antisymmetric matroids with coefficients in tracts $F$ with $1=-1$.
Beforehand, we prove the converse of Proposition~\ref{prop: antisym to sym} for even symmetric matroids.

\begin{theorem}\label{thm: even to antisymmetric}
    Let $M = (\ground{n},\cB)$ be an even symmetric matroid.
    There is unique $\cB' \subseteq \cA_n$ such that $M' = (\ground{n}, \cB \cup \cB')$ is an antisymmetric matroid.
\end{theorem}

Theorem~\ref{thm: even to antisymmetric} is easily deduced from the circuit axiom of even symmetric matroids by Booth, Moreira, and Pinto~\cite{BMP2003}. 

\begin{theorem}[{\cite[Theorem~12]{BMP2003}}]\label{thm: circuits of even symmetric}
    Let $\cC$ be a set of subsets $C$ of $E = \ground{n}$ such that $C$ is a subtransversal.
    Then $\cC$ is the set of circuits of an even symmetric matroid if and only if it satisfies \ref{item: sC1},~\ref{item: sC2},~\ref{item: isotropic}, and~\ref{item: maximal}.
\end{theorem}

Remark that
the original statement in~\cite{BMP2003} includes the additional condition~\ref{item: addition'}, 
but it can be omitted because of Lemma~\ref{lem: addition}.

\begin{proof}[\mbox{\bf Proof of Theorem~\ref{thm: even to antisymmetric}}]
    Let $\cB'$ be the set of almost-transversals $A \in \cA_n$ such that for some $x,y\in E$ with $\skewpair{x} \subseteq A$ and $\skewpair{y} \cap A = \emptyset$, both $A-x+y$ and $A-x^*+y^*$ are in $\cB$.
    Let $\cC$ be the set of circuits of the even symmetric matroid $M$.
    By Theorems~\ref{thm: circuits of even symmetric} and~\ref{thm: cryptomorphism}, $\cC$ is the set of circuits of an antisymmetric matroid, say $M'$, on $E$.

    We claim that $\cB(M') = \cB \cup \cB'$.
    It is equivalent to show that for each $B\in \cT_n\cup \cA_n$, the set $B$ is in $\cB\cup \cB'$ if and only if $B$ does not contain any $C\in \cC$.
    
    Suppose that $B\in \cB\cup \cB'$.
    If $B\in \cB \subseteq \cT_n$, then $B$ is not contained in any $C\in \cC$ because it is a basis of the even symmetric matroid $M$.
    If $B\in \cB' \subseteq \cA_n$, then $B$ is not contained in any $C\in \cC$ because $B$ contains a skew pair and $C$ is a subtransversal.
    Thus, $B$ is not a superset of any elements of~$\cC$.
    
    Conversely, suppose that $B$ does not contain any $C\in \cC$.
    If $B\in \cT_n$, then $B$ is a basis of the even symmetric matroid $M$ and equivalently $B\in \cB$.
    Thus, we may assume that $B\in \cA_n$.
    Let $\skewpair{x}$ and $\skewpair{y}$ be the skew pairs such that $\skewpair{x} \subseteq B$ and $\skewpair{y} \cap B = \emptyset$.
    By~\ref{item: maximal}, there is $C \in \cC$ such that $C \subseteq B+y$.
    Then $y \in C$ since $C\not\subseteq B$.
    As $C$ is a subtransversal, by interchanging $x$ and $x^*$ if necessary, we may assume that $x^* \notin C$.
    Then there is no $D\in \cC$ such that $D\subseteq B-x^*+y^*$, because otherwise $C\cap D^* = \{y\}$ contradicting~\ref{item: isotropic}.
    Hence $B-x^*+y^* \in \cB$.
    Assume that $x\notin C$.
    Then similarly there is no $D\in \cC$ such that $D\subseteq B-x+y^*$, and hence $B-x+y^* \in \cB$, which violates that $M$ is even.
    Therefore, $x \in C$.
    Assume that there is $D\in \cC$ such that $D\subseteq B-x+y$.
    If $x^*\in D$, then $C\cap D^* = \{x\}$, contradicting~\ref{item: isotropic}.
    Thus, $x^*\notin D$.
    Hence $D\subseteq B-x-x^*+y$ and $D \ne C$.
    Also, $y \in D$ because $D\not\subseteq B$.
    By~\ref{item: addition}, there is $D'\in \cC$ such that $D' \subseteq C+D-y \subseteq B$, a contradiction.
    Therefore, no such $D$ exists and $B-x+y \in \cB$.
    Thus, we conclude that $B \in \cB'$ because $B-x^*+y^*, B-x+y \in \cB$.
    Therefore, the claim is proved.
    
    It remains to check that if there is an antisymmetric matroid $M''$ such that $\cB(M'') \cap \cT_n = \cB$, then $\cB(M'') \cap \cA_n = \cB'$.
    It is straightforward from Lemma~\ref{lem: 3-term basis exchange}.
\end{proof}

Proposition~\ref{prop: embedding of matroids to antisym mats} can be generalized by replacing matroids with even delta-matroids, whenever $k$ has characteristic two.
It is due to Theorem~\ref{thm: even to antisymmetric}.

\begin{corollary}
    Let $k$ be a field of characteristic two.
    Then an even delta-matroid $M$ on $[n]$ admits an $[n]$-by-$[n]$ symmetric matrix $A$ over $k$ such that $\cB(M) = \{X\subseteq [n] : \det(A[X,X]) \ne 0\}$ if and only if the antisymmetric matroid associated with $\lift(M)$ is representable over~$k$.
\end{corollary}

Proposition~\ref{prop: embedding of matroids to antisym mats} cannot be further extended to delta-matroids because of Example~\ref{eg: two sym mat with same antisym mat}.
We remark that an even delta-matroid $M$ admits a symmetric matrix representation over a field $k$ if and only if either $M$ is a matroid up to twisting or $k$ has characteristic two; see Geelen~{\cite[Page~27]{Geelen1996thesis}}.

\subsection{Gaussoids}\label{sec: gaussoids}

A gaussoid is a combinatorial structure introduced by Ln\v{e}ni\v{c}ka and Mat\'{u}\v{s}~\cite{LM2007} for understanding which almost-principal submatrices of a positive definite symmetric matrix can be simultaneously singular.
We review an equivalent definition in terms of edge relations by~{\cite[Theorem~1]{BDKS2019}}.

A subset $\cG$ of $\cA_n$ is \emph{allowable} if an almost-transversal $A$ is in $\cG$ if and only if $A-p+q$ is in $\cG$, where $p$ and $q$ are skew pairs with $p\subseteq A$ and $q\cap A=\emptyset$.
An allowable subset $\cG$ of $\cA_n$ is \emph{incompatible with} a restricted G--P relation $f := \sum_{e\in S_1 \setminus S_2} (-1)^{\smaller{S_1}{e} + \smaller{S_2}{e}} x_{S_1-e} x_{S_2+e} = 0$ if there is precisely one term $x_{S_1-e}x_{S_2+e}$ in $f$ such that neither $S_1-e$ nor $S_2+e$ is in $\cG$.
Otherwise $\cG$ is \emph{compatible with} $f=0$.

Recall, in Section~\ref{sec: Preliminaries}, that every $3$-term restricted G--P relation is classified as either a square relation
\[
    x_{Sab}x_{Sa^*b^*} + x_{Sab^*}x_{Sa^*b} - x_{Saa^*}x_{Sbb^*} = 0,
\]
where $S+ab$ is a transversal, 
or an edge relation
\[
    (-1)^{\smaller{L}{a}} x_{S abc} x_{S bb^*c^*}
    +
    (-1)^{\smaller{L}{b^*}}
    x_{S abc^*} x_{S bb^*c}
    +
    (-1)^{\smaller{L}{c^*}}
    x_{S abb^*} x_{S bcc^*}
    =0,
\]
where $S+abc$ is a transversal and $L = \{a,c,b^*,c^*\}$.
In the edge relation, $Sabc$ and $Sabc^*$ are transversals, and $Sbb^*c^*$, $Sbb^*c^*$, $Sabb^*$, and $Sbcc^*$ are almost-transversals.

\begin{definition}[\cite{LM2007}; see~{\cite[Theorem~1]{BDKS2019}}]
    A subset $\cG$ of $\cA_n$ is a \emph{gaussoid} if it is allowable and is compatible with all edge relations.
\end{definition}

\begin{example}
    Let $\Sigma$ be a positive definite symmetric matrix over the real field $\bR$, and let $\Lambda := \begin{bmatrix} I \, | \, \Sigma \end{bmatrix}$.
    Then $\cG := \{B \in \cA_n : \det(\Lambda[n,B]) = 0\}$ is a gaussoid, and $\cB := \{B \in \cT_n\cup \cA_n : \det(\Lambda[n,B]) \ne 0\}$ is the set of bases of an antisymmetric matroid.
    Note that $\cT_n \subseteq \cB$ because $\Sigma$ is positive definite, and $\cG = \cA_n \setminus \cB$.
\end{example}

In general, we can obtain gaussoids from antisymmetric matroids containing all transversals.

\begin{proposition}
    Let $M = (\ground{n},\cB)$ be an antisymmetric matroid.
    If $\cT_n\subseteq \cB$, then $\cA_n \setminus \cB$ is a gaussoid.
\end{proposition}
\begin{proof}
    Let $\cG := \cA_n \setminus \cB$.
    By~\ref{item: sB2}, $\cG$ is allowable.
    By~\ref{item: exchange} and the assumption $\cT_n\subseteq \cB$, the set $\cG$ is compatible with all restricted G--P relations and thus it is a gaussoid.
\end{proof}

It is open whether, for every gaussoid $\cG$, $\cT_n \cup (\cA_n \setminus \cG)$ is the set of bases of an antisymmetric matroid.
We remark that the previous statement holds if Conjecture~1 in~\cite{BDKS2019} is true.

\section{Homotopy theorem}\label{sec: homotopy}

We show that the first homology group of a graph associated with an antisymmetric matroid is generated by short cycles, which will be used to prove the equivalence of two notions of antisymmetric matroids over tracts in Sections~\ref{sec: antisymmetric matroids over tracts}--\ref{sec: cryptomorphism}.
Our result implies the homotopy theorem for basis graphs of matroids by Maurer~\cite{Maurer1973} and that of even symmetric matroids by Wenzel~\cite{Wenzel1993}.

A \emph{transversal basis graph} $\cG_M$ of an antisymmetric matroid $M$ on $\ground{n}$ is a graph such that
\begin{itemize}
    \item its vertex set is $\cB(M) \cap \cT_n$, and 
    \item two vertices $B$ and $B'$ are adjacent if and only if (i) $|B\setminus B'| = 1$ or (ii) $|B\setminus B'| = 2$ and there is $A\in \cB(M)\cap \cA_n$ such that $|B\setminus A| = |B'\setminus A| = 1$. %
\end{itemize}
The \emph{weight} of an edge $BB'$ is defined as $\eta(BB'):= |B\setminus B'|$.

For a subgraph $G$ of $\cG_M$ and two vertices $B,B'$ in $G$, let $\dist_{G}(B,B')$ be the smallest sum $\eta(P) = \sum_{e\in E(P)} \eta(e)$ among all paths $P$ from $B$ to $B'$ in $G$.
For convenience, we write $\dist_{M}(B,B')$ for $\dist_{\cG_M}(B,B')$.

\begin{lemma}\label{lem: cGm}
    The following hold.
    \begin{enumerate}[label=\rm(\roman*)]
        \item $\dist_{M}(B,B') = |B\setminus B'|$.
        \item For every cycle $C$ in $\cG_M$, the weight $\eta(C) = \sum_{e\in E(C)} \eta(e)$ is even. 
    \end{enumerate}
\end{lemma}
\begin{proof}
    (i)
    Clearly, $\dist_{M}(B,B') \ge |B \setminus B'|$.
    We claim $\dist_{M}(B,B') \le |B \setminus B'|$ by induction on $|B\setminus B'|$, which trivially holds for $|B\setminus B'| = 0$.
    We may assume $B\ne B'$.
    By~\ref{item: exchange}, there are $e\in B\setminus B'$ and $f\in B'\setminus B$ such that $B'':=B-e+f \in \cB(M)$.
    By the induction hypothesis, $\dist_{M}(B,B') \le \dist_{M}(B,B'') + \dist_{M}(B'',B') \le 1+|B''\setminus B'| = |B\setminus B'|$.

    (ii)
    Fix a vertex $B_0$ of $\cG_M$ and for each nonnegative integer $i$, let $V_i$ be the set of vertices $B$ such that $|B\setminus B_0| = i$.
    By interchanging $j$ and $j^*$ for some $j\in\ground{n}$, we may assume that $B_0 = [n]$.
    Then $V_i = \{B\in \cB : |B \cap [n]^*| = i\}$.
    Then it is easily observed that for every edge $BB'$ in $\cG_M$ such that $|B\setminus B_0| \le |B'\setminus B_0|$,
    \begin{enumerate}[label=\rm(\roman*)]
        \item $\eta(BB') = 1$ if and only if $B \in V_i$ and $B' \in V_{i+1}$ for some $i$, and 
        \item $\eta(BB') = 2$ if and only if $B \in V_i$ and $B'\in V_i \cup V_{i+2}$ for some $i$.
    \end{enumerate}
    Thus, if a path $P$ is from $B_0$ to a vertex in $V_i$, then $\eta(P) \equiv i \pmod{2}$.
    Then it is straightforward that for every cycle $C$, its weight $\eta(C)$ is even.
\end{proof}

\begin{definition}[Graph Homology]
    For a graph $G = (V,E)$ and a linear ordering $\prec$ on $V$, let $\partial_1 : \bZ^E \to \bZ^V$ be a group homomorphism such that $\partial_1(vw) = v-w$ if $v\prec w$.
    The %
    \emph{(first) homology group} of $G$ is $H(G) := \ker(\partial_1)$.
\end{definition}

For a fixed graph $G$, its homology group $H(G)$ is unique up to isomorphism for different choices of linear orderings on the vertex set.
We will often identify a cycle $v_1v_2\ldots v_kv_1$ in $G$ with an element $\sum_{i=1}^{k} \epsilon_iv_iv_{i+1} \in H(G)$ where $\epsilon_i= \begin{cases} 1 & \text{if $v_i\prec v_{i+1}$} \\ -1 & \text{otherwise} \end{cases}$ and $v_{k+1}:= v_1$.
Note that the cycles of $G$ generate $H(G)$.

We denote the homology group of $\cG_M$ by $H_M$.
We call a cycle $C$ in $\cG_M$ is \emph{reducible} if, in $H_M$, it can be generated by the cycles of weight smaller than $\eta(C)$.
Otherwise, we say $C$ is \emph{irreducible}.
Now we are ready to state the Homotopy Theorem for transversal basis graphs.

\begin{theorem}[Homotopy Theorem] %
    \label{thm: homotopy}
    The homology group $H_M$ is generated by cycles $C$ of $\eta(C) \le 8$.
\end{theorem}

Our result implies Maurer's homotopy theorem for basis graphs of matroids~\cite{Maurer1973}. %
We account for this implication in a bit more general setting concerning even symmetric matroids.
The \emph{basis graph} of an even symmetric matroid $N$ is a graph on $\cB(N)$ such that two vertices $B$ and $B'$ are adjacent if and only if $|B\setminus B'| = 2$.
By Theorem~\ref{thm: even to antisymmetric}, there is an antisymmetric matroid $M$ such that $\cB(M) \cap \cT_n = \cB(N)$.
Clearly, the basis graph of $N$ is identical to the transversal basis graph of $M$.
Hence, Wenzel's homotopy theorem~\cite{Wenzel1995} for even symmetric matroid is deduced immediately. %

\begin{corollary}[\cite{Wenzel1995}] %
    The homology group of the basis graph of an even symmetric matroid is generated by cycles of length at most four.
\end{corollary}

The \emph{basis graph} of a matroid is a graph on the set of bases such that two vertices $B,B'$ are adjacent if and only if $|B\setminus B'|=1$.

\begin{corollary}[\cite{Maurer1973}] %
    The homology group of the basis graph of a matroid is generated by cycles of length at most four.
\end{corollary}

In order to prove Theorem~\ref{thm: homotopy}, we first observe several properties of irreducible cycles in $\cG_M$.

\begin{lemma}\label{lem: irreducible1}
    Let $C$ be an irreducible cycle in $\cG_M$.
    Then for every pair of vertices $B,B'$ in $C$, $\dist_C(B,B') = \dist_M(B,B')$.
\end{lemma}
\begin{proof}
    Suppose not.
    We take a pair of distinct vertices $B,B'$ in $C$ and a path $P$ from $B$ to $B'$ in $C$ such that
    \begin{enumerate}[label=\rm(\roman*)]
        \item $\eta(P) \le \eta(C)/2$,
        \item $\eta(P) > \dist_M(B,B')$, and
        \item subject to (i) and (ii), $\dist_M(B,B')$ is minimized.
    \end{enumerate}
    Let $P'$ be the path from $B$ to $B'$ in $C$ other than $P$, and 
    let $Q$ be a path from $B$ to $B'$ in $\cG_M$ such that $\eta(Q) = \dist_M(B,B')$.
    By (iii), no internal vertex of $Q$ is in $C$.
    Then the cycle induced by paths $P$ and $Q$ has weight $\eta(P)+\eta(Q) < 2\eta(P) \le \eta(C)$, and the cycle induced by paths $P'$ and $Q$ has weight $\eta(P')+\eta(Q) < \eta(P')+\eta(P) = \eta(C)$.
    It contradicts that $C$ is irreducible.
\end{proof}

\begin{lemma}\label{lem: irreducible2}
    Let $C$ be a cycle of weight $2\ell$ in $\cG_M$ and let $B_0\in V(C)$.
    If $C$ is irreducible, then either
    \begin{itemize}
        \item there is $B\in V(C)$ such that $\dist_M(B_0,B) = \ell$ or
        \item there is $BB' \in E(C)$ such that $\dist_M(B_0,B) = \dist_M(B_0,B') = \ell-1$ and $\eta(BB')=2$.
    \end{itemize}
\end{lemma}

\begin{proof}[\mbox{\bf Proof of Theorem~\ref{thm: homotopy}}]
    Suppose to the contrary that $\cG_M$ has irreducible cycles of weight larger than~$8$.
    Among such cycles, 
    we choose $C$ such that 
    \begin{enumerate}[label=\rm(\roman*)]
        \item its weight $\eta(C)$ is minimized, and
        \item subject to (i), the number $|E(C)\cap \eta^{-1}(1)|$ of weight $1$ edges in $C$ is minimized. %
    \end{enumerate}
    We denote by $\eta(C) = 2\ell > 8$.
    Then all cycles of weight $2\ell'$ with $4< \ell' < \ell$ are reducible. %
    We select an arbitrary vertex $B_0 \in V(C)$.
    There are two cases by Lemma~\ref{lem: irreducible2}.

    
    \noindent\textbf{Case I.}
    There is a vertex $B$ in $C$ such that $\dist_M(B_0,B) = \ell$.
    Let $B_1$ and $B_2$ be two distinct neighbors of $B$ in $C$, and let $P$ be a path in $C$ from $B_1$ to $B_2$ containing $B_0$.
    %


    \noindent\textbf{Subcase I.1.}
    $\eta(BB_1) = \eta(BB_2) = 1$.
    Then there are distinct elements $e,f\in B$ such that $B_1 = B\symdiff\skewpair{f}$ and $B_2 = B\symdiff\skewpair{e}$.
    By (ii), $B_1B_2$ is not an edge in~$\cG_M$ and hence $B+e^*-f$ is not a basis of~$M$.
    Thus, by the basis exchange~\ref{item: exchange} applied to $B_1$ and $B_2$, $B' := B\symdiff\{e,e^*,f,f^*\}$ is a basis; see Figure~\ref{fig: homotopy1}(top left).
    Then $B_1B', B_2B'$ are edges in $\cG_M$ of weight $1$ and $\dist_M(B_0,B') = \ell-2$.
    If $X$ is a vertex in $C$ such that $\dist_M(B_0,X) = \ell-2$, then $\dist_{C}(B_i,X) = 3$ for some $i\in\{1,2\}$ and by Lemma~\ref{lem: irreducible1}, $\dist_M(B_i,X) = 3$.
    Thus, $B'\notin V(C)$.
    Let $C'$ be a cycle concatenating two paths $P$ and $B_1B'B_2$.
    Then $\dist_{C'}(B_0,B') = \ell$ and thus $C'$ is reducible by Lemma~\ref{lem: irreducible1}.
    Let $C'' := BB_1B'B_2B$ be a cycle of weight $4$.
    As $C = a C' + b C''$ in $H_M$ for some $a,b\in\{1,-1\}$, it contradicts that $C$ is irreducible.
    %
    %

    %
    %
    %


    \noindent\textbf{Subcase I.2.}
    $\eta(BB_1) = 2$ and $\eta(BB_2) = 1$.
    Then there are elements $e,f,g\in B$ such that $B_1 = B\symdiff\{f,f^*,g,g^*\}$ and $B_2 = B\symdiff\skewpair{e}$.
    Because $\dist_{M}(B_1,B_2) = \dist_{C}(B_1,B_2)=3$, the elements $e,f,g$ are distinct.
    By~\ref{item: exchange} applied to $B_1$ and $B_2$, there is $h\in \{e,f^*,g^*\} = B_1\setminus B_2$ such that $B_1+e^*-h$ is a basis of~$M$.

    Suppose that $B' := B_1+e^*-e$ is a basis; see Figure~\ref{fig: homotopy1}(top middle).
    Then $\dist_{M}(B_2,B')=2$ and $\dist_{M}(B_0,B') = \ell-3$.
    If $X$ is a vertex in $C$ such that $\dist_{M}(B_0,X) = \ell-3$, then either $\dist_{C}(B_2,X) = 4$ or $\dist_{C}(B_1,X) = 5$ and by Lemma~\ref{lem: irreducible1}, $\dist_M(B_2,X) = 4$ or $\dist_M(B_5,X) = 5$.
    Thus, $B' \notin V(C)$.
    Let $Q$ be a path from $B_2$ to $B'$ of weight $\dist_M(B_2,B')=2$.
    Then by Lemma~\ref{lem: irreducible1}, a cycle induced by $P$, $Q$, and $B'B_1$ is reducible.
    Since a cycle induced by two paths $Q$ and $B'B_1BB_2$ has weight $6$, we deduce that $C$ is reducible, a contradiction.
    Therefore, $B_1+e^*-e$ is not a basis.
    Then $h\ne e$.

    By symmetry, we can assume that $h=f^*$.
    By Lemma~\ref{lem: 3-term basis exchange}, one of two transversals $(B_1+e^*-f^*)-e+f$ and $(B_1+e^*-f^*)-e+f^*$ is a basis.
    Hence $B'' := (B_1+e^*-f^*)-e+f = B\symdiff\{e,e^*,g,g^*\}$ is a basis; see Figure~\ref{fig: homotopy1}(top right).
    Then $B_1B''$ is an edge of weight $2$, $B_2B''$ is an edge of weight $1$, and $\dist_M(B_0,B'') = \ell-2$.
    Similar to before, we can deduce that $B''\notin V(C)$.
    Let $Q$ be a path from $B_2$ to $B''$ of weight $\dist_M(B_2,B'')=2$.
    Then by Lemma~\ref{lem: irreducible1}, a cycle induced by $P$, $Q$, and $B''B_1$ is reducible.
    A cycle consisting of two paths $Q$ and $B''B_1BB_2$ has weight $6$.
    This contradicts that $C$ is irreducible.


    \noindent\textbf{Subcase I.3.}
    $\eta(BB_1) = \eta(BB_2) = 2$.
    Then there are four elements $e,f,g,h\in B$ such that $B_1 = B\symdiff\{g,g^*,h,h^*\}$ and $B_2 = B\symdiff\{e,e^*,f,f^*\}$.
    Because $\dist_M(B_1,B_2) = \dist_C(B_1,B_2) = 4$, these four elements are distinct; see Figure~\ref{fig: homotopy1}(bottom).
    By the basis exchange~\ref{item: exchange}, there is $i\in \{e,f,g^*,h^*\} = B_1\setminus B_2$ such that $B_1+e^*-i$ is a basis.

    Suppose that $B':=B_1+e^*-e$ is a basis.
    Let $Q$ be a path from $B'$ to $B_2$ of weight $\dist_M(B',B_2)=3$.
    By Lemma~\ref{lem: irreducible1} applied to $C$, $B_2$ is the only vertex in both paths $P$ and $Q$.
    A cycle consisting of two paths $Q$ and $B'B_1BB_2$ has weight $8$.
    A cycle consisting of $P$, $Q$, and $B'B_1$ is reducible by Lemma~\ref{lem: irreducible1}.
    It contradicts that $C$ is irreducible.
    Therefore, $B_1+e^*-e$ is not a basis and $i \ne e$.

    Suppose that $i = f$.
    By Lemma~\ref{lem: 3-term basis exchange}, $(B_1+e^*-f)-e+f$ or $(B_1+e^*-f)-e+f^*$ is a basis.
    Hence $B'' := (B_1+e^*-f)-e+f^*$ is a basis.
    Then $\dist_M(B_0,B'') = \ell-4$ and $B''\notin V(C)$.
    Let $Q$ be a path from $B''$ to $B_2$ of weight $\dist_M(B_2,B'')=2$.
    Then $V(P)\cap V(Q) = \{B_2\}$ by Lemma~\ref{lem: irreducible1}.
    Then the union of $C$, $Q$, and $B''B_1$ has two cycles other than $C$.
    One has weight $8$ and the other is reducible by Lemma~\ref{lem: irreducible1}, a contradiction.
    Thus, $i \in \{g^*,h^*\}$.

    By symmetry, we can assume that $i = g^*$.
    By Lemma~\ref{lem: 3-term basis exchange}, $(B_1+e^*-g^*)-e+g$ or $(B_1+e^*-g^*)-e+g^*$ is a basis.
    Hence $B''' := (B_1+e^*-g^*)-e+g = B\symdiff\{e,e^*,h,h^*\}$ is a basis. 
    Then $\dist_M(B_0,B''') = \ell-2$ and $B'''\notin V(C)$.
    Similarly, Lemma~\ref{lem: irreducible1} yields a contradiction, and we skip details.

\begin{figure}
    \centering
    \begin{tikzpicture}
        \node[v] (B) at (0,0) {};
        \node[v] (B1) at (-1,-0.7) {};
        \node[v] (B2) at (1,-0.7) {};
        \draw (B1)node[above]{\scriptsize{$B_1$}}--(B)node[left]{\scriptsize{$B$}}--(B2)node[above]{\scriptsize{$B_2$}};

        \begin{scope}[xshift=-0.15cm, yshift=0.3cm]
            \node[box1s] () at (0.3*0,0.3) {};
            \node[box1s] () at (0.3*1,0.3) {};

            \node[box2s] () at (0.3*0,0) {};
            \node[box2s] () at (0.3*1,0) {};

            \node () at (0.3*0,0.3) {\scriptsize{$e$}};
            \node () at (0.3*1,0.3) {\scriptsize{$f$}};
        \end{scope}
        
        \begin{scope}[xshift=-1.7cm, yshift=-0.85cm]
            \node[box1s] () at (0.3*0,0.3) {};
            \node[box2s] () at (0.3*1,0.3) {};

            \node[box2s] () at (0.3*0,0) {};
            \node[box1s] () at (0.3*1,0) {};

            \node () at (0.3*1,0) {\scriptsize{$f^*$}};
        \end{scope}

        \begin{scope}[xshift=1.4cm, yshift=-0.85cm]
            \node[box2s] () at (0.3*0,0.3) {};
            \node[box1s] () at (0.3*1,0.3) {};

            \node[box1s] () at (0.3*0,0) {};
            \node[box2s] () at (0.3*1,0) {};

            \node () at (0.3*0,0) {\scriptsize{$e^*$}};
        \end{scope}

        \node[v] (B') at (0,-1.4) {};
        \draw (B1)--(B')node[left]{\scriptsize{$B'$}}--(B2);

        \begin{scope}[xshift=-0.15cm, yshift=-2cm]
            \node[box2s] () at (0.3*0,0.3) {};
            \node[box2s] () at (0.3*1,0.3) {};

            \node[box1s] () at (0.3*0,0) {};
            \node[box1s] () at (0.3*1,0) {};    
        \end{scope}

        \node[v] (B0) at (0,-4) {};
        \draw[dashed] (B1)to[bend right=30](B0)node[below]{\scriptsize{$B_0$}};
        \draw[dashed] (B0)to[bend right=30](B2);

        \node () at (0,-3.4) {I.1};
    \end{tikzpicture}
    \hspace{1.35cm}
    \begin{tikzpicture}
        \node[v] (B) at (0,0) {};
        \node[v] (B1) at (-1,-1.4) {};
        \node[v] (B2) at (1,-0.7) {};
        \draw[double] (B1)node[above]{\scriptsize{$B_1$}}--(B)node[left]{\scriptsize{$B$}};
        \draw (B)--(B2)node[above]{\scriptsize{$B_2$}};

        \begin{scope}[xshift=-0.3cm, yshift=0.3cm]
            \node[box1s] () at (0.3*0,0.3) {};
            \node[box1s] () at (0.3*1,0.3) {};
            \node[box1s] () at (0.3*2,0.3) {};

            \node[box2s] () at (0.3*0,0) {};
            \node[box2s] () at (0.3*1,0) {};
            \node[box2s] () at (0.3*2,0) {};

            \node () at (0.3*0,0.3) {\scriptsize{$e$}};
            \node () at (0.3*1,0.3) {\scriptsize{$f$}};
            \node () at (0.3*2,0.3) {\scriptsize{$g$}};

        \end{scope}
        
        \begin{scope}[xshift=-2.0cm, yshift=-1.55cm]
            \node[box1s] () at (0.3*0,0.3) {};
            \node[box2s] () at (0.3*1,0.3) {};
            \node[box2s] () at (0.3*2,0.3) {};

            \node[box2s] () at (0.3*0,0) {};
            \node[box1s] () at (0.3*1,0) {};
            \node[box1s] () at (0.3*2,0) {};

            \node () at (0.3*1,0) {\scriptsize{$f^*$}};
            \node () at (0.3*2,0) {\scriptsize{$g^*$}};
        \end{scope}

        \begin{scope}[xshift=1.4cm, yshift=-0.85cm]
            \node[box2s] () at (0.3*0,0.3) {};
            \node[box1s] () at (0.3*1,0.3) {};
            \node[box1s] () at (0.3*2,0.3) {};

            \node[box1s] () at (0.3*0,0) {};
            \node[box2s] () at (0.3*1,0) {};
            \node[box2s] () at (0.3*2,0) {};

            \node () at (0.3*0,0) {\scriptsize{$e^*$}};
        \end{scope}

        \node[v] (B') at (0,-2.1) {};
        \draw (B1)--(B')node[left]{\scriptsize{$B'$}};
        \draw[dashed] (B')--node[midway]{\scriptsize{2}}(B2);

        \begin{scope}[xshift=-0.3cm, yshift=-2.7cm]
            \node[box2s] () at (0.3*0,0.3) {};
            \node[box2s] () at (0.3*1,0.3) {};
            \node[box2s] () at (0.3*2,0.3) {};

            \node[box1s] () at (0.3*0,0) {};
            \node[box1s] () at (0.3*1,0) {};
            \node[box1s] () at (0.3*2,0) {};
        \end{scope}

        \node[v] (B0) at (0,-4) {};
        \draw[dashed] (B1)to[bend right=30](B0)node[below]{\scriptsize{$B_0$}};
        \draw[dashed] (B0)to[bend right=30](B2);

        \node () at (0,-3.4) {I.2};
    \end{tikzpicture}
    \hspace{1.15cm}
    \begin{tikzpicture}
        \node[v] (B) at (0,0) {};
        \node[v] (B1) at (-1,-1.4) {};
        \node[v] (B2) at (1,-0.7) {};
        \draw[double] (B1)node[above]{\scriptsize{$B_1$}}--(B)node[left]{\scriptsize{$B$}};
        \draw (B)--(B2)node[above]{\scriptsize{$B_2$}};

        \begin{scope}[xshift=-0.3cm, yshift=0.3cm]
            \node[box1s] () at (0.3*0,0.3) {};
            \node[box1s] () at (0.3*1,0.3) {};
            \node[box1s] () at (0.3*2,0.3) {};

            \node[box2s] () at (0.3*0,0) {};
            \node[box2s] () at (0.3*1,0) {};
            \node[box2s] () at (0.3*2,0) {};
        \end{scope}
        
        \begin{scope}[xshift=-2.0cm, yshift=-1.55cm]
            \node[box1s] () at (0.3*0,0.3) {};
            \node[box2s] () at (0.3*1,0.3) {};
            \node[box2s] () at (0.3*2,0.3) {};

            \node[box2s] () at (0.3*0,0) {};
            \node[box1s] () at (0.3*1,0) {};
            \node[box1s] () at (0.3*2,0) {};
        \end{scope}

        \begin{scope}[xshift=1.4cm, yshift=-0.85cm]
            \node[box2s] () at (0.3*0,0.3) {};
            \node[box1s] () at (0.3*1,0.3) {};
            \node[box1s] () at (0.3*2,0.3) {};

            \node[box1s] () at (0.3*0,0) {};
            \node[box2s] () at (0.3*1,0) {};
            \node[box2s] () at (0.3*2,0) {};
        \end{scope}

        \node[v] (B') at (0,-1.4) {};
        \draw[double] (B1)--(B')node[above]{\scriptsize{$B''$}};
        \draw (B')--(B2);

        \begin{scope}[xshift=-0.3cm, yshift=-2.0cm]
            \node[box2s] () at (0.3*0,0.3) {};
            \node[box1s] () at (0.3*1,0.3) {};
            \node[box2s] () at (0.3*2,0.3) {};

            \node[box1s] () at (0.3*0,0) {};
            \node[box2s] () at (0.3*1,0) {};
            \node[box1s] () at (0.3*2,0) {};
        \end{scope}

        \node[v] (B0) at (0,-4) {};
        \draw[dashed] (B1)to[bend right=30](B0)node[below]{\scriptsize{$B_0$}};
        \draw[dashed] (B0)to[bend right=30](B2);

        \node () at (0,-3.4) {I.2};
    \end{tikzpicture}

    \begin{tikzpicture}
        \node[v] (B) at (0,0) {};
        \node[v] (B1) at (-1.2,-1.2) {};
        \node[v] (B2) at (1.2,-1.2) {};
        \draw[double] (B1)node[above]{\scriptsize{$B_1$}}--(B)node[left]{\scriptsize{$B$}}--(B2)node[above]{\scriptsize{$B_2$}};

        \begin{scope}[xshift=-0.45cm, yshift=0.3cm]
            \node[box1s] () at (0.3*0,0.3) {};
            \node[box1s] () at (0.3*1,0.3) {};
            \node[box1s] () at (0.3*2,0.3) {};
            \node[box1s] () at (0.3*3,0.3) {};

            \node[box2s] () at (0.3*0,0) {};
            \node[box2s] () at (0.3*1,0) {};
            \node[box2s] () at (0.3*2,0) {};
            \node[box2s] () at (0.3*3,0) {};

            \node () at (0.3*0,0.3) {\scriptsize{$e$}};
            \node () at (0.3*1,0.3) {\scriptsize{$f$}};
            \node () at (0.3*2,0.3) {\scriptsize{$g$}};
            \node () at (0.3*3,0.3) {\scriptsize{$h$}};
        \end{scope}
        
        \begin{scope}[xshift=-2.5cm, yshift=-1.35cm]
            \node[box1s] () at (0.3*0,0.3) {};
            \node[box1s] () at (0.3*1,0.3) {};
            \node[box2s] () at (0.3*2,0.3) {};
            \node[box2s] () at (0.3*3,0.3) {};

            \node[box2s] () at (0.3*0,0) {};
            \node[box2s] () at (0.3*1,0) {};
            \node[box1s] () at (0.3*2,0) {};
            \node[box1s] () at (0.3*3,0) {};

            \node () at (0.3*2,0) {\scriptsize{$g^*$}};
            \node () at (0.3*3,0) {\scriptsize{$h^*$}};
        \end{scope}

        \begin{scope}[xshift=1.6cm, yshift=-1.35cm]
            \node[box2s] () at (0.3*0,0.3) {};
            \node[box2s] () at (0.3*1,0.3) {};
            \node[box1s] () at (0.3*2,0.3) {};
            \node[box1s] () at (0.3*3,0.3) {};

            \node[box1s] () at (0.3*0,0) {};
            \node[box1s] () at (0.3*1,0) {};
            \node[box2s] () at (0.3*2,0) {};
            \node[box2s] () at (0.3*3,0) {};

            \node () at (0.3*0,0) {\scriptsize{$e^*$}};
            \node () at (0.3*1,0) {\scriptsize{$f^*$}};
        \end{scope}

        \node[v] (B') at (-0.3,-1.8) {};
        \draw (B1)--(B')node[above]{\scriptsize{$B'$}};
        \draw[dashed] (B')--node[midway]{\scriptsize{3}}(B2);

        \begin{scope}[xshift=-0.75cm, yshift=-2.4cm]
            \node[box2s] () at (0.3*0,0.3) {};
            \node[box1s] () at (0.3*1,0.3) {};
            \node[box2s] () at (0.3*2,0.3) {};
            \node[box2s] () at (0.3*3,0.3) {};

            \node[box1s] () at (0.3*0,0) {};
            \node[box2s] () at (0.3*1,0) {};
            \node[box1s] () at (0.3*2,0) {};
            \node[box1s] () at (0.3*3,0) {};
        \end{scope}

        \node[v] (B0) at (0,-4) {};
        \draw[dashed] (B1)to[bend right=30](B0)node[below]{\scriptsize{$B_0$}};
        \draw[dashed] (B0)to[bend right=30](B2);
        
        \node () at (0,-3.4) {I.3};
    \end{tikzpicture}
    \hspace{0.1cm}
    \begin{tikzpicture}
        \node[v] (B) at (0,0) {};
        \node[v] (B1) at (-1.2,-1.2) {};
        \node[v] (B2) at (1.2,-1.2) {};
        \draw[double] (B1)node[above]{\scriptsize{$B_1$}}--(B)node[left]{\scriptsize{$B$}}--(B2)node[above]{\scriptsize{$B_2$}};

        \begin{scope}[xshift=-0.45cm, yshift=0.3cm]
            \node[box1s] () at (0.3*0,0.3) {};
            \node[box1s] () at (0.3*1,0.3) {};
            \node[box1s] () at (0.3*2,0.3) {};
            \node[box1s] () at (0.3*3,0.3) {};

            \node[box2s] () at (0.3*0,0) {};
            \node[box2s] () at (0.3*1,0) {};
            \node[box2s] () at (0.3*2,0) {};
            \node[box2s] () at (0.3*3,0) {};
        \end{scope}
        
        \begin{scope}[xshift=-2.5cm, yshift=-1.35cm]
            \node[box1s] () at (0.3*0,0.3) {};
            \node[box1s] () at (0.3*1,0.3) {};
            \node[box2s] () at (0.3*2,0.3) {};
            \node[box2s] () at (0.3*3,0.3) {};

            \node[box2s] () at (0.3*0,0) {};
            \node[box2s] () at (0.3*1,0) {};
            \node[box1s] () at (0.3*2,0) {};
            \node[box1s] () at (0.3*3,0) {};
        \end{scope}

        \begin{scope}[xshift=1.6cm, yshift=-1.35cm]
            \node[box2s] () at (0.3*0,0.3) {};
            \node[box2s] () at (0.3*1,0.3) {};
            \node[box1s] () at (0.3*2,0.3) {};
            \node[box1s] () at (0.3*3,0.3) {};

            \node[box1s] () at (0.3*0,0) {};
            \node[box1s] () at (0.3*1,0) {};
            \node[box2s] () at (0.3*2,0) {};
            \node[box2s] () at (0.3*3,0) {};
        \end{scope}

        \node[v] (B') at (0,-2.4) {};
        \draw[double] (B1)--(B')node[left]{\scriptsize{$B''$}};
        \draw[dashed] (B')--node[midway]{\scriptsize{2}}(B2);

        \begin{scope}[xshift=-0.45cm, yshift=-3.0cm]
            \node[box2s] () at (0.3*0,0.3) {};
            \node[box2s] () at (0.3*1,0.3) {};
            \node[box2s] () at (0.3*2,0.3) {};
            \node[box2s] () at (0.3*3,0.3) {};

            \node[box1s] () at (0.3*0,0) {};
            \node[box1s] () at (0.3*1,0) {};
            \node[box1s] () at (0.3*2,0) {};
            \node[box1s] () at (0.3*3,0) {};
        \end{scope}

        \node[v] (B0) at (0,-4) {};
        \draw[dashed] (B1)to[bend right=30](B0)node[below]{\scriptsize{$B_0$}};
        \draw[dashed] (B0)to[bend right=30](B2);

        \node () at (0,-3.4) {I.3};
    \end{tikzpicture}
    \hspace{0.1cm}
    \begin{tikzpicture}
        \node[v] (B) at (0,0) {};
        \node[v] (B1) at (-1.2,-1.2) {};
        \node[v] (B2) at (1.2,-1.2) {};
        \draw[double] (B1)node[above]{\scriptsize{$B_1$}}--(B)node[left]{\scriptsize{$B$}}--(B2)node[above]{\scriptsize{$B_2$}};

        \begin{scope}[xshift=-0.45cm, yshift=0.3cm]
            \node[box1s] () at (0.3*0,0.3) {};
            \node[box1s] () at (0.3*1,0.3) {};
            \node[box1s] () at (0.3*2,0.3) {};
            \node[box1s] () at (0.3*3,0.3) {};

            \node[box2s] () at (0.3*0,0) {};
            \node[box2s] () at (0.3*1,0) {};
            \node[box2s] () at (0.3*2,0) {};
            \node[box2s] () at (0.3*3,0) {};
        \end{scope}
        
        \begin{scope}[xshift=-2.5cm, yshift=-1.35cm]
            \node[box1s] () at (0.3*0,0.3) {};
            \node[box1s] () at (0.3*1,0.3) {};
            \node[box2s] () at (0.3*2,0.3) {};
            \node[box2s] () at (0.3*3,0.3) {};

            \node[box2s] () at (0.3*0,0) {};
            \node[box2s] () at (0.3*1,0) {};
            \node[box1s] () at (0.3*2,0) {};
            \node[box1s] () at (0.3*3,0) {};
        \end{scope}

        \begin{scope}[xshift=1.6cm, yshift=-1.35cm]
            \node[box2s] () at (0.3*0,0.3) {};
            \node[box2s] () at (0.3*1,0.3) {};
            \node[box1s] () at (0.3*2,0.3) {};
            \node[box1s] () at (0.3*3,0.3) {};

            \node[box1s] () at (0.3*0,0) {};
            \node[box1s] () at (0.3*1,0) {};
            \node[box2s] () at (0.3*2,0) {};
            \node[box2s] () at (0.3*3,0) {};
        \end{scope}

        \node[v] (B') at (0,-1.2) {};
        \draw[double] (B1)--(B')node[above]{\scriptsize{$B'''$}};
        \draw[dashed] (B')--node[midway]{\scriptsize{2}}(B2);

        \begin{scope}[xshift=-0.45cm, yshift=-1.8cm]
            \node[box2s] () at (0.3*0,0.3) {};
            \node[box1s] () at (0.3*1,0.3) {};
            \node[box1s] () at (0.3*2,0.3) {};
            \node[box2s] () at (0.3*3,0.3) {};

            \node[box1s] () at (0.3*0,0) {};
            \node[box2s] () at (0.3*1,0) {};
            \node[box2s] () at (0.3*2,0) {};
            \node[box1s] () at (0.3*3,0) {};
        \end{scope}

        \node[v] (B0) at (0,-4) {};
        \draw[dashed] (B1)to[bend right=30](B0)node[below]{\scriptsize{$B_0$}};
        \draw[dashed] (B0)to[bend right=30](B2);

        \node () at (0,-3.4) {I.3};
    \end{tikzpicture}
    \caption{Descriptions of Case I in the proof of Theorem~\ref{thm: homotopy}. Solid lines represent edges of weight $1$, double lines represent edges of weight $2$, and dashed lines represent paths in $\cG_M$.}
    \label{fig: homotopy1}
\end{figure}


    \noindent\textbf{Case II.}
    There is an edge $BB'$ in $C$ such that $\dist_M(B_0,B) = \dist(B_0,B') = \ell-1$ and $\eta(BB')=2$.
    Then $B' = B \symdiff \{e,e^*,f,f^*\}$ for some distinct elements $e,f^* \in B$.
    Let $B_1$ be the neighbor of $B$ in $C$ which is not $B'$.
    %
    

    \noindent\textbf{Subcase II.1.}
    $\eta(BB_1) = 2$.
    Then $B_1 = B\symdiff\{g,g^*,h,h^*\}$ for some distinct $g,h\in B\setminus\{e,f^*\}$; see Figure~\ref{fig: homotopy2.1}.
    Let $P$ be a path in $C$ from $B_1$ to $B'$ containing $B_0$.
    By the basis exchange~\ref{item: exchange}, there is $i\in B_1\setminus B' = \{e,f^*,g^*,h^*\}$ such that $B_1+e^*-i$ is a basis of~$M$.

    Suppose that $D:=B_1+e^*-e$ is a basis.
    Then $\dist_M(B_0,D) = \ell-4$ and $D \notin V(C)$.
    Let $Q$ be a path from $B'$ to $D$ of weight $\dist_M(B',D)=3$.
    Note that $V(P) \cap V(Q) = \{B'\}$ by Lemma~\ref{lem: irreducible1} applied to $C$.
    Then a cycle induced by $P$, $Q$, and $DB_1$ is reducible by Lemma~\ref{lem: irreducible1}, and
    a cycle concatenating two paths $Q$ and $DB_1BB'$ has length $8$.
    It contradicts that $C$ is irreducible.
    Therefore, $B_1+e^*-e$ is not a basis 
    and $i \ne e$.

    Suppose that $i = f^*$.
    By Lemma~\ref{lem: 3-term basis exchange}, $(B_1+e^*-f^*)-e+f$ or $(B_1+e^*-f^*)-e+f^*$ is a basis.
    Hence $D' := (B_1+e^*-f^*)-e+f = B' \symdiff \{g,g^*,h,h^*\}$ is a basis.
    Then $\dist_M(B_0,B'') = \ell-3$ and $D'\notin V(C)$.
    Using Lemma~\ref{lem: irreducible1}, we can similarly conclude that $C$ is reducible, a contradiction.
    Thus, $i\ne f^*$ and so $i$ is either $g^*$ or $h^*$.

    By symmetry, we may assume that $i=g^*$.
    By Lemma~\ref{lem: 3-term basis exchange}, 
    $D'' := (B_1+e^*-g^*)-e+g = B'\symdiff\{f,f^*,h,h^*\}$ is a basis.
    Then $\dist_M(B_0,D'') = \ell-2$ and $D''\notin V(C)$.
    Using Lemma~\ref{lem: irreducible1}, similarly one can deduce a contradiction.
    %


    \noindent\textbf{Subcase II.2.}
    $\eta(BB_1) = 1$.
    Let $B_2$ be the neighbor of $B'$ in $V(C)$ other than $B$.
    By Subcase II.1, we may assume that $\eta(BB_2) = 1$.
    Then for some distinct elements $g,h\in B\cap B'$, we have $B_1 = B\symdiff\skewpair{h}$ and $B_2 = B'\symdiff\skewpair{g}$; see Figure~\ref{fig: homotopy2.2}.
    Let $P$ be a path in $C$ from $B_1$ to $B_2$ containing $B_0$.
    By the basis exchange~\ref{item: exchange}, there is $i\in B_1\setminus B_2 = \{e,f^*,g,h^*\}$ such that $B_1+e^*-i$ is a basis of $M$.

    Suppose that $D:=B_1+e^*-e$ is a basis.
    Then $\dist_M(B_0,D) = \ell-3$ and $D \notin V(C)$.
    Let $Q$ be a path from $B_2$ to $D$ of weight $\dist_M(B_2,D) = 3$.
    By Lemma~\ref{lem: irreducible1}, $P$ and $Q$ only meet at a vertex $B_2$.
    A cycle consisting of two paths $DB_1BB'B_2$ and $Q$ has length $8$, and a cycle consisting of $P$, $Q$, and $DB_1$ is reducible by Lemma~\ref{lem: irreducible1}.
    It contradicts that $C$ is irreducible.
    Therefore, $B_1+e^*-e$ is not a basis and $i \ne e$.

    Suppose that $i = f^*$.
    By Lemma~\ref{lem: 3-term basis exchange}, $D' := (B_1+e^*-f^*)-e+f = B_2\symdiff\{g,g^*,h,h^*\}$ is a basis.
    Then $\dist_M(B_0,B'') = \ell-2$ and $D'\notin V(C)$.
    Similarly, 
    we can deduce a contradiction using Lemma~\ref{lem: irreducible1}.
    Thus, we can assume that $i\ne f^*$.

    Suppose that $i=h^*$.
    By Lemma~\ref{lem: 3-term basis exchange}, $D'' := (B_1+e^*-h^*)-e+h = B_2 \symdiff\{f,f^*,g,g^*\}$ is a basis.
    Then $\dist_M(B_0,D'') = \ell-2$ and $D''\notin V(C)$.
    So, we can deduce a contradiction similarly as before.

    Thus, we may assume that $i=g$.
    By Lemma~\ref{lem: 3-term basis exchange}, $D''' := (B_1+e^*-g)-e+g^*$ is a basis.
    Then $\dist_M(B_0,D''') = \ell-4$ and $D'''\notin V(C)$.
    One can deduce a contradiction similarly as before.
\end{proof}

\begin{figure}
    \centering
    \begin{tikzpicture}
        \node[v] (B) at (-1,0) {};
        \node[v] (B') at (1,0) {};
        \draw[double] (B)node[above]{\scriptsize{$B$}}--(B')node[above]{\scriptsize{$B'$}};
        \node[v] (B1) at (-1.3,-1.6) {};
        \draw[double] (B)--(B1)node[right]{\scriptsize{$B_1$}};
        
        \begin{scope}[xshift=-2.3cm, yshift=-0.15cm]
            \node[box1s] () at (0.3*0,0.3) {};
            \node[box2s] () at (0.3*1,0.3) {};
            \node[box1s] () at (0.3*2,0.3) {};
            \node[box1s] () at (0.3*3,0.3) {};

            \node[box2s] () at (0.3*0,0) {};
            \node[box1s] () at (0.3*1,0) {};
            \node[box2s] () at (0.3*2,0) {};
            \node[box2s] () at (0.3*3,0) {};

            \node () at (0.3*0,0.3) {\scriptsize{$e$}};
            \node () at (0.3*1,0) {\scriptsize{$f^*$}};
            \node () at (0.3*2,0.3) {\scriptsize{$g$}};
            \node () at (0.3*3,0.3) {\scriptsize{$h$}};
        \end{scope}

        \begin{scope}[xshift=1.4cm, yshift=-0.15cm]
            \node[box2s] () at (0.3*0,0.3) {};
            \node[box1s] () at (0.3*1,0.3) {};
            \node[box1s] () at (0.3*2,0.3) {};
            \node[box1s] () at (0.3*3,0.3) {};

            \node[box1s] () at (0.3*0,0) {};
            \node[box2s] () at (0.3*1,0) {};
            \node[box2s] () at (0.3*2,0) {};
            \node[box2s] () at (0.3*3,0) {};

            \node () at (0.3*0,0) {\scriptsize{$e^*$}};
            \node () at (0.3*1,0.3) {\scriptsize{$f$}};
        \end{scope}

        \begin{scope}[xshift=-2.6cm, yshift=-1.75cm]
            \node[box1s] () at (0.3*0,0.3) {};
            \node[box2s] () at (0.3*1,0.3) {};
            \node[box2s] () at (0.3*2,0.3) {};
            \node[box2s] () at (0.3*3,0.3) {};

            \node[box2s] () at (0.3*0,0) {};
            \node[box1s] () at (0.3*1,0) {};
            \node[box1s] () at (0.3*2,0) {};
            \node[box1s] () at (0.3*3,0) {};

            \node () at (0.3*2,0) {\scriptsize{$g^*$}};
            \node () at (0.3*3,0) {\scriptsize{$h^*$}};
        \end{scope}

        \node[v] (D) at (0,-2.4) {};
        \draw (B1)--(D)node[right]{\scriptsize{$D$}};
        \draw[dashed] (D)--node[midway]{\scriptsize{3}}(B');

        \begin{scope}[xshift=-0.45cm, yshift=-3.00cm]
            \node[box2s] () at (0.3*0,0.3) {};
            \node[box2s] () at (0.3*1,0.3) {};
            \node[box2s] () at (0.3*2,0.3) {};
            \node[box2s] () at (0.3*3,0.3) {};

            \node[box1s] () at (0.3*0,0) {};
            \node[box1s] () at (0.3*1,0) {};
            \node[box1s] () at (0.3*2,0) {};
            \node[box1s] () at (0.3*3,0) {};
        \end{scope}

        \node[v] (B0) at (0,-4) {};
        \draw[dashed] (B1)to[bend right=30](B0)node[below]{\scriptsize{$B_0$}};
        \draw[dashed] (B0)to[bend right=40](B');

        \node () at (0,-3.4) {II.1};
    \end{tikzpicture}
    \hspace{0.3cm}
    \begin{tikzpicture}
        \node[v] (B) at (-1,0) {};
        \node[v] (B') at (1,0) {};
        \draw[double] (B)node[above]{\scriptsize{$B$}}--(B')node[above]{\scriptsize{$B'$}};
        \node[v] (B1) at (-1.3,-1.6) {};
        \draw[double] (B)--(B1)node[above right]{\scriptsize{$B_1$}};
        
        \begin{scope}[xshift=-2.3cm, yshift=-0.15cm]
            \node[box1s] () at (0.3*0,0.3) {};
            \node[box2s] () at (0.3*1,0.3) {};
            \node[box1s] () at (0.3*2,0.3) {};
            \node[box1s] () at (0.3*3,0.3) {};

            \node[box2s] () at (0.3*0,0) {};
            \node[box1s] () at (0.3*1,0) {};
            \node[box2s] () at (0.3*2,0) {};
            \node[box2s] () at (0.3*3,0) {};

        \end{scope}

        \begin{scope}[xshift=1.4cm, yshift=-0.15cm]
            \node[box2s] () at (0.3*0,0.3) {};
            \node[box1s] () at (0.3*1,0.3) {};
            \node[box1s] () at (0.3*2,0.3) {};
            \node[box1s] () at (0.3*3,0.3) {};

            \node[box1s] () at (0.3*0,0) {};
            \node[box2s] () at (0.3*1,0) {};
            \node[box2s] () at (0.3*2,0) {};
            \node[box2s] () at (0.3*3,0) {};

        \end{scope}

        \begin{scope}[xshift=-2.6cm, yshift=-1.75cm]
            \node[box1s] () at (0.3*0,0.3) {};
            \node[box2s] () at (0.3*1,0.3) {};
            \node[box2s] () at (0.3*2,0.3) {};
            \node[box2s] () at (0.3*3,0.3) {};

            \node[box2s] () at (0.3*0,0) {};
            \node[box1s] () at (0.3*1,0) {};
            \node[box1s] () at (0.3*2,0) {};
            \node[box1s] () at (0.3*3,0) {};

        \end{scope}

        \node[v] (D) at (0.3,-1.6) {};
        \draw[double] (B1)--(D)node[right]{\scriptsize{$D'$}};
        \draw[dashed] (D)--node[midway]{\scriptsize{2}}(B');

        \begin{scope}[xshift=-0.15cm, yshift=-2.2cm]
            \node[box2s] () at (0.3*0,0.3) {};
            \node[box1s] () at (0.3*1,0.3) {};
            \node[box2s] () at (0.3*2,0.3) {};
            \node[box2s] () at (0.3*3,0.3) {};

            \node[box1s] () at (0.3*0,0) {};
            \node[box2s] () at (0.3*1,0) {};
            \node[box1s] () at (0.3*2,0) {};
            \node[box1s] () at (0.3*3,0) {};
        \end{scope}

        \node[v] (B0) at (0,-4) {};
        \draw[dashed] (B1)to[bend right=30](B0)node[below]{\scriptsize{$B_0$}};
        \draw[dashed] (B0)to[bend right=40](B');

        \node () at (0,-3.4) {II.1};
    \end{tikzpicture}
    \hspace{0.3cm}
    \begin{tikzpicture}
        \node[v] (B) at (-1,0) {};
        \node[v] (B') at (1,0) {};
        \draw[double] (B)node[above]{\scriptsize{$B$}}--(B')node[above]{\scriptsize{$B'$}};
        \node[v] (B1) at (-1.3,-1.6) {};
        \draw[double] (B)--(B1)node[above right]{\scriptsize{$B_1$}};
        
        \begin{scope}[xshift=-2.3cm, yshift=-0.15cm]
            \node[box1s] () at (0.3*0,0.3) {};
            \node[box2s] () at (0.3*1,0.3) {};
            \node[box1s] () at (0.3*2,0.3) {};
            \node[box1s] () at (0.3*3,0.3) {};

            \node[box2s] () at (0.3*0,0) {};
            \node[box1s] () at (0.3*1,0) {};
            \node[box2s] () at (0.3*2,0) {};
            \node[box2s] () at (0.3*3,0) {};

        \end{scope}

        \begin{scope}[xshift=1.4cm, yshift=-0.15cm]
            \node[box2s] () at (0.3*0,0.3) {};
            \node[box1s] () at (0.3*1,0.3) {};
            \node[box1s] () at (0.3*2,0.3) {};
            \node[box1s] () at (0.3*3,0.3) {};

            \node[box1s] () at (0.3*0,0) {};
            \node[box2s] () at (0.3*1,0) {};
            \node[box2s] () at (0.3*2,0) {};
            \node[box2s] () at (0.3*3,0) {};

        \end{scope}

        \begin{scope}[xshift=-2.6cm, yshift=-1.75cm]
            \node[box1s] () at (0.3*0,0.3) {};
            \node[box2s] () at (0.3*1,0.3) {};
            \node[box2s] () at (0.3*2,0.3) {};
            \node[box2s] () at (0.3*3,0.3) {};

            \node[box2s] () at (0.3*0,0) {};
            \node[box1s] () at (0.3*1,0) {};
            \node[box1s] () at (0.3*2,0) {};
            \node[box1s] () at (0.3*3,0) {};

        \end{scope}

        \node[v] (D) at (0.3,-1.6) {};
        \draw[double] (B1)--(D)node[right]{\scriptsize{$D''$}};
        \draw[dashed] (D)--node[midway]{\scriptsize{2}}(B');

        \begin{scope}[xshift=-0.15cm, yshift=-2.2cm]
            \node[box2s] () at (0.3*0,0.3) {};
            \node[box2s] () at (0.3*1,0.3) {};
            \node[box1s] () at (0.3*2,0.3) {};
            \node[box2s] () at (0.3*3,0.3) {};

            \node[box1s] () at (0.3*0,0) {};
            \node[box1s] () at (0.3*1,0) {};
            \node[box2s] () at (0.3*2,0) {};
            \node[box1s] () at (0.3*3,0) {};
        \end{scope}

        \node[v] (B0) at (0,-4) {};
        \draw[dashed] (B1)to[bend right=30](B0)node[below]{\scriptsize{$B_0$}};
        \draw[dashed] (B0)to[bend right=40](B');

        \node () at (0,-3.4) {II.1};
    \end{tikzpicture}
    \caption{Descriptions of Case II.1 in the proof of Theorem~\ref{thm: homotopy}.}
    \label{fig: homotopy2.1}
\end{figure}

\begin{figure}
    \centering
    \begin{tikzpicture}
        \node[v] (B) at (-0.9,0) {};
        \node[v] (B') at (0.9,0) {};
        \node[v] (B1) at (-1.1,-1.1) {};
        \node[v] (B2) at (1.1,-1.1) {};
        \draw[double] (B)node[above]{\scriptsize{$B$}}--(B')node[above]{\scriptsize{$B'$}};
        \draw (B)--(B1)node[right]{\scriptsize{$B_1$}};
        \draw (B')--(B2)node[left]{\scriptsize{$B_2$}};
        
        \begin{scope}[xshift=-2.2cm, yshift=-0.15cm]
            \node[box1s] () at (0.3*0,0.3) {};
            \node[box2s] () at (0.3*1,0.3) {};
            \node[box1s] () at (0.3*2,0.3) {};
            \node[box1s] () at (0.3*3,0.3) {};

            \node[box2s] () at (0.3*0,0) {};
            \node[box1s] () at (0.3*1,0) {};
            \node[box2s] () at (0.3*2,0) {};
            \node[box2s] () at (0.3*3,0) {};

            \node () at (0.3*0,0.3) {\scriptsize{$e$}};
            \node () at (0.3*1,0) {\scriptsize{$f^*$}};
            \node () at (0.3*2,0.3) {\scriptsize{$g$}};
            \node () at (0.3*3,0.3) {\scriptsize{$h$}};
        \end{scope}

        \begin{scope}[xshift=1.3cm, yshift=-0.15cm]
            \node[box2s] () at (0.3*0,0.3) {};
            \node[box1s] () at (0.3*1,0.3) {};
            \node[box1s] () at (0.3*2,0.3) {};
            \node[box1s] () at (0.3*3,0.3) {};

            \node[box1s] () at (0.3*0,0) {};
            \node[box2s] () at (0.3*1,0) {};
            \node[box2s] () at (0.3*2,0) {};
            \node[box2s] () at (0.3*3,0) {};

            \node () at (0.3*0,0) {\scriptsize{$e^*$}};
            \node () at (0.3*1,0.3) {\scriptsize{$f$}};
        \end{scope}

        \begin{scope}[xshift=-2.4cm, yshift=-1.25cm]
            \node[box1s] () at (0.3*0,0.3) {};
            \node[box2s] () at (0.3*1,0.3) {};
            \node[box1s] () at (0.3*2,0.3) {};
            \node[box2s] () at (0.3*3,0.3) {};

            \node[box2s] () at (0.3*0,0) {};
            \node[box1s] () at (0.3*1,0) {};
            \node[box2s] () at (0.3*2,0) {};
            \node[box1s] () at (0.3*3,0) {};

            \node () at (0.3*3,0) {\scriptsize{$h^*$}};
        \end{scope}

        \begin{scope}[xshift=1.5cm, yshift=-1.25cm]
            \node[box2s] () at (0.3*0,0.3) {};
            \node[box1s] () at (0.3*1,0.3) {};
            \node[box2s] () at (0.3*2,0.3) {};
            \node[box1s] () at (0.3*3,0.3) {};

            \node[box1s] () at (0.3*0,0) {};
            \node[box2s] () at (0.3*1,0) {};
            \node[box1s] () at (0.3*2,0) {};
            \node[box2s] () at (0.3*3,0) {};

            \node () at (0.3*2,0) {\scriptsize{$g^*$}};
        \end{scope}

        \node[v] (D) at (-0.3,-1.8) {};
        \begin{scope}[xshift=-0.75cm, yshift=-2.4cm]
            \node[box2s] () at (0.3*0,0.3) {};
            \node[box2s] () at (0.3*1,0.3) {};
            \node[box1s] () at (0.3*2,0.3) {};
            \node[box2s] () at (0.3*3,0.3) {};

            \node[box1s] () at (0.3*0,0) {};
            \node[box1s] () at (0.3*1,0) {};
            \node[box2s] () at (0.3*2,0) {};
            \node[box1s] () at (0.3*3,0) {};
        \end{scope}
        \draw (B1)--(D)node[above]{\scriptsize{$D$}};
        \draw[dashed] (D)--node[midway]{\scriptsize{3}}(B2);

        \node[v] (B0) at (0,-4) {};
        \draw[dashed] (B1)to[bend right=30](B0)node[below]{\scriptsize{$B_0$}};
        \draw[dashed] (B0)to[bend right=30](B2);

        \node () at (0,-3.4) {II.2};
    \end{tikzpicture}
    \hspace{0.3cm}
    \begin{tikzpicture}
        \node[v] (B) at (-0.9,0) {};
        \node[v] (B') at (0.9,0) {};
        \node[v] (B1) at (-1.1,-1.1) {};
        \node[v] (B2) at (1.1,-1.1) {};
        \draw[double] (B)node[above]{\scriptsize{$B$}}--(B')node[above]{\scriptsize{$B'$}};
        \draw (B)--(B1)node[above right]{\scriptsize{$B_1$}};
        \draw (B')--(B2)node[above left]{\scriptsize{$B_2$}};
        
        \begin{scope}[xshift=-2.2cm, yshift=-0.15cm]
            \node[box1s] () at (0.3*0,0.3) {};
            \node[box2s] () at (0.3*1,0.3) {};
            \node[box1s] () at (0.3*2,0.3) {};
            \node[box1s] () at (0.3*3,0.3) {};

            \node[box2s] () at (0.3*0,0) {};
            \node[box1s] () at (0.3*1,0) {};
            \node[box2s] () at (0.3*2,0) {};
            \node[box2s] () at (0.3*3,0) {};

        \end{scope}

        \begin{scope}[xshift=1.3cm, yshift=-0.15cm]
            \node[box2s] () at (0.3*0,0.3) {};
            \node[box1s] () at (0.3*1,0.3) {};
            \node[box1s] () at (0.3*2,0.3) {};
            \node[box1s] () at (0.3*3,0.3) {};

            \node[box1s] () at (0.3*0,0) {};
            \node[box2s] () at (0.3*1,0) {};
            \node[box2s] () at (0.3*2,0) {};
            \node[box2s] () at (0.3*3,0) {};

        \end{scope}

        \begin{scope}[xshift=-2.4cm, yshift=-1.25cm]
            \node[box1s] () at (0.3*0,0.3) {};
            \node[box2s] () at (0.3*1,0.3) {};
            \node[box1s] () at (0.3*2,0.3) {};
            \node[box2s] () at (0.3*3,0.3) {};

            \node[box2s] () at (0.3*0,0) {};
            \node[box1s] () at (0.3*1,0) {};
            \node[box2s] () at (0.3*2,0) {};
            \node[box1s] () at (0.3*3,0) {};

        \end{scope}

        \begin{scope}[xshift=1.5cm, yshift=-1.25cm]
            \node[box2s] () at (0.3*0,0.3) {};
            \node[box1s] () at (0.3*1,0.3) {};
            \node[box2s] () at (0.3*2,0.3) {};
            \node[box1s] () at (0.3*3,0.3) {};

            \node[box1s] () at (0.3*0,0) {};
            \node[box2s] () at (0.3*1,0) {};
            \node[box1s] () at (0.3*2,0) {};
            \node[box2s] () at (0.3*3,0) {};

        \end{scope}

        \node[v] (D) at (0,-1.1) {};
        \begin{scope}[xshift=-0.45cm, yshift=-1.7cm]
            \node[box2s] () at (0.3*0,0.3) {};
            \node[box1s] () at (0.3*1,0.3) {};
            \node[box1s] () at (0.3*2,0.3) {};
            \node[box2s] () at (0.3*3,0.3) {};

            \node[box1s] () at (0.3*0,0) {};
            \node[box2s] () at (0.3*1,0) {};
            \node[box2s] () at (0.3*2,0) {};
            \node[box1s] () at (0.3*3,0) {};
        \end{scope}
        \draw[double] (B1)--(D)node[above]{\scriptsize{$D'$}};
        \draw[dashed] (D)--node[midway]{\scriptsize{2}}(B2);

        \node[v] (B0) at (0,-4) {};
        \draw[dashed] (B1)to[bend right=30](B0)node[below]{\scriptsize{$B_0$}};
        \draw[dashed] (B0)to[bend right=30](B2);

        \node () at (0,-3.4) {II.2};
    \end{tikzpicture}

    \begin{tikzpicture}
        \node[v] (B) at (-0.9,0) {};
        \node[v] (B') at (0.9,0) {};
        \node[v] (B1) at (-1.1,-1.1) {};
        \node[v] (B2) at (1.1,-1.1) {};
        \draw[double] (B)node[above]{\scriptsize{$B$}}--(B')node[above]{\scriptsize{$B'$}};
        \draw (B)--(B1)node[above right]{\scriptsize{$B_1$}};
        \draw (B')--(B2)node[above left]{\scriptsize{$B_2$}};
        
        \begin{scope}[xshift=-2.2cm, yshift=-0.15cm]
            \node[box1s] () at (0.3*0,0.3) {};
            \node[box2s] () at (0.3*1,0.3) {};
            \node[box1s] () at (0.3*2,0.3) {};
            \node[box1s] () at (0.3*3,0.3) {};

            \node[box2s] () at (0.3*0,0) {};
            \node[box1s] () at (0.3*1,0) {};
            \node[box2s] () at (0.3*2,0) {};
            \node[box2s] () at (0.3*3,0) {};

        \end{scope}

        \begin{scope}[xshift=1.3cm, yshift=-0.15cm]
            \node[box2s] () at (0.3*0,0.3) {};
            \node[box1s] () at (0.3*1,0.3) {};
            \node[box1s] () at (0.3*2,0.3) {};
            \node[box1s] () at (0.3*3,0.3) {};

            \node[box1s] () at (0.3*0,0) {};
            \node[box2s] () at (0.3*1,0) {};
            \node[box2s] () at (0.3*2,0) {};
            \node[box2s] () at (0.3*3,0) {};

        \end{scope}

        \begin{scope}[xshift=-2.4cm, yshift=-1.25cm]
            \node[box1s] () at (0.3*0,0.3) {};
            \node[box2s] () at (0.3*1,0.3) {};
            \node[box1s] () at (0.3*2,0.3) {};
            \node[box2s] () at (0.3*3,0.3) {};

            \node[box2s] () at (0.3*0,0) {};
            \node[box1s] () at (0.3*1,0) {};
            \node[box2s] () at (0.3*2,0) {};
            \node[box1s] () at (0.3*3,0) {};

        \end{scope}

        \begin{scope}[xshift=1.5cm, yshift=-1.25cm]
            \node[box2s] () at (0.3*0,0.3) {};
            \node[box1s] () at (0.3*1,0.3) {};
            \node[box2s] () at (0.3*2,0.3) {};
            \node[box1s] () at (0.3*3,0.3) {};

            \node[box1s] () at (0.3*0,0) {};
            \node[box2s] () at (0.3*1,0) {};
            \node[box1s] () at (0.3*2,0) {};
            \node[box2s] () at (0.3*3,0) {};

        \end{scope}

        \node[v] (D) at (0,-1.1) {};
        \begin{scope}[xshift=-0.45cm, yshift=-1.7cm]
            \node[box2s] () at (0.3*0,0.3) {};
            \node[box2s] () at (0.3*1,0.3) {};
            \node[box1s] () at (0.3*2,0.3) {};
            \node[box1s] () at (0.3*3,0.3) {};

            \node[box1s] () at (0.3*0,0) {};
            \node[box1s] () at (0.3*1,0) {};
            \node[box2s] () at (0.3*2,0) {};
            \node[box2s] () at (0.3*3,0) {};
        \end{scope}
        \draw[double] (B1)--(D)node[above]{\scriptsize{$D''$}};
        \draw[dashed] (D)--node[midway]{\scriptsize{2}}(B2);

        \node[v] (B0) at (0,-4) {};
        \draw[dashed] (B1)to[bend right=30](B0)node[below]{\scriptsize{$B_0$}};
        \draw[dashed] (B0)to[bend right=30](B2);

        \node () at (0,-3.4) {II.2};
    \end{tikzpicture}
    \hspace{0.3cm}
    \begin{tikzpicture}
        \node[v] (B) at (-0.9,0) {};
        \node[v] (B') at (0.9,0) {};
        \node[v] (B1) at (-1.1,-1.1) {};
        \node[v] (B2) at (1.1,-1.1) {};
        \draw[double] (B)node[above]{\scriptsize{$B$}}--(B')node[above]{\scriptsize{$B'$}};
        \draw (B)--(B1)node[right]{\scriptsize{$B_1$}};
        \draw (B')--(B2)node[left]{\scriptsize{$B_2$}};
        
        \begin{scope}[xshift=-2.2cm, yshift=-0.15cm]
            \node[box1s] () at (0.3*0,0.3) {};
            \node[box2s] () at (0.3*1,0.3) {};
            \node[box1s] () at (0.3*2,0.3) {};
            \node[box1s] () at (0.3*3,0.3) {};

            \node[box2s] () at (0.3*0,0) {};
            \node[box1s] () at (0.3*1,0) {};
            \node[box2s] () at (0.3*2,0) {};
            \node[box2s] () at (0.3*3,0) {};

        \end{scope}

        \begin{scope}[xshift=1.3cm, yshift=-0.15cm]
            \node[box2s] () at (0.3*0,0.3) {};
            \node[box1s] () at (0.3*1,0.3) {};
            \node[box1s] () at (0.3*2,0.3) {};
            \node[box1s] () at (0.3*3,0.3) {};

            \node[box1s] () at (0.3*0,0) {};
            \node[box2s] () at (0.3*1,0) {};
            \node[box2s] () at (0.3*2,0) {};
            \node[box2s] () at (0.3*3,0) {};

        \end{scope}

        \begin{scope}[xshift=-2.4cm, yshift=-1.25cm]
            \node[box1s] () at (0.3*0,0.3) {};
            \node[box2s] () at (0.3*1,0.3) {};
            \node[box1s] () at (0.3*2,0.3) {};
            \node[box2s] () at (0.3*3,0.3) {};

            \node[box2s] () at (0.3*0,0) {};
            \node[box1s] () at (0.3*1,0) {};
            \node[box2s] () at (0.3*2,0) {};
            \node[box1s] () at (0.3*3,0) {};

        \end{scope}

        \begin{scope}[xshift=1.5cm, yshift=-1.25cm]
            \node[box2s] () at (0.3*0,0.3) {};
            \node[box1s] () at (0.3*1,0.3) {};
            \node[box2s] () at (0.3*2,0.3) {};
            \node[box1s] () at (0.3*3,0.3) {};

            \node[box1s] () at (0.3*0,0) {};
            \node[box2s] () at (0.3*1,0) {};
            \node[box1s] () at (0.3*2,0) {};
            \node[box2s] () at (0.3*3,0) {};

        \end{scope}

        \node[v] (D) at (0,-2.2) {};
        \begin{scope}[xshift=-0.45cm, yshift=-2.8cm]
            \node[box2s] () at (0.3*0,0.3) {};
            \node[box2s] () at (0.3*1,0.3) {};
            \node[box2s] () at (0.3*2,0.3) {};
            \node[box2s] () at (0.3*3,0.3) {};

            \node[box1s] () at (0.3*0,0) {};
            \node[box1s] () at (0.3*1,0) {};
            \node[box1s] () at (0.3*2,0) {};
            \node[box1s] () at (0.3*3,0) {};
        \end{scope}
        \draw[double] (B1)--(D)node[right]{\scriptsize{$D'''$}};
        \draw[dashed] (D)--node[midway]{\scriptsize{2}}(B2);

        \node[v] (B0) at (0,-4) {};
        \draw[dashed] (B1)to[bend right=30](B0)node[below]{\scriptsize{$B_0$}};
        \draw[dashed] (B0)to[bend right=30](B2);

        \node () at (0,-3.4) {II.2};
    \end{tikzpicture}
    \caption{Descriptions of Case II.2 in the proof of Theorem~\ref{thm: homotopy}.}
    \label{fig: homotopy2.2}
\end{figure}

\section{Antisymmetric matroids over tracts}\label{sec: antisymmetric matroids over tracts}

A \emph{tract} is a field-like structure defined by Baker and Bowler~\cite{BB2019}.
It encompasses fields, partial fields, and hyperfields, and thus matroids with coefficients in a tract generalize partial field representations of matroids~\cite{SW1996}, oriented matroids~\cite{BVSWZ1999}, valuated matroids~\cite{DW1992b}, and ordinary matroids.

We define antisymmetric matroids with coefficients in tracts.
As a byproduct, it provides a concept for Lagrangian Grassmannians over hyperfields.
We establish antisymmetric matroids with coefficients in two ways which generalize, first, a point in the projective space satisfying the restricted G--P relations in Section~\ref{sec: antisym mat over tracts via bases} and, second, a Lagrangian subspace in the symplectic vector space in Section~\ref{sec: antisym mat over tracts via circuits}.
We will show these two notions are equivalent in Section~\ref{sec: cryptomorphism}.
In Sections~\ref{sec: examples}--\ref{sec: orthgonal matroids with coefficients}, we give examples of antisymmetric matroids over tracts, including matroids over tracts introduced in~\cite{BB2019}.
In Sections~\ref{sec: symplectic dressian} and~\ref{sec: oriented gaussoids}, we compare antisymmetric matroids over the tropical hyperfield~$\bT$ and the sign hyperfield~$\bS$ with symplectic Dressians~\cite{BO2023} and oriented gaussoids~\cite{BDKS2019}, respectively.

We first review the definition of tracts and some basic properties.

\begin{definition}[Tracts]
    A \emph{tract} is a pair $F = (G, N_F)$ of an abelian group $G$, written multiplicatively, and a subset $N_F$ of the group semiring $\bN[G]$ satisfying the following axioms:
    \begin{enumerate}[label=\rm(T\arabic*)]
        \item The zero element $0$ of $\bN[G]$ is in $N_F$.
        \item The identity $1$ of $G$ is not in $N_F$.
        \item There is a unique $\epsilon \in G$ such that $1+\epsilon \in N_F$.
        \item If $g\in G$ and $\sum_{i=1}^n h_i \in N_F$, then $\sum_{i=1}^n gh_i \in N_F$.
    \end{enumerate}
    Abusing a notation, we write $F$ for the set $G \cup \{0\}$.
    We denote by $F^\times := G$.
    We call $N_F$ the {\em null set}, which can be regarded as the set of linear combinations of $G$ that sum to zero.
\end{definition}

\begin{lemma}[\cite{BB2019}]
    Let $F$ be a tract. The following holds.
    \begin{enumerate}[label=\rm(\roman*)]
        \item $\epsilon^2 = 1$.
        \item For $x,y\in F^\times$, if $x+y \in N_F$, then $y = \epsilon x$.
        \item $F^\times \cap N_F = \emptyset$.
    \end{enumerate}
\end{lemma}

Because of the previous lemma, we often write $-1$ instead of $\epsilon$.
A \emph{morphism} $f$ from a tract $F_1$ to another tract $F_2$ is a map such that $f(0) = 0$ and the restriction $f: F_1 ^\times \to F_2 ^\times$ is a group homomorphism inducing $\varphi(N_{F_1}) \subseteq N_{F_2}$.

\begin{example}\label{eg: tracts}
    We give tracts associated with fields, partial fields, and some hyperfields without precise definitions of partial fields and hyperfields; for more details, see~\cite{BB2019}.
    \begin{enumerate}[label=(\arabic*)]
        \item For a field $k$, let $N_k$ be the set of linear combinations of nonzero elements that sum to zero in~$k$.
        Then $(k^\times,N_k)$ is a tract.
        \item\label{item: partial field} Let $R$ be a commutative ring with unity and $G$ be its unit subgroup, i.e., $G \le R^\times$.
        A {\em partial field} $P$ associated with $G$ and $R$ is a tract such that $P^\times = G$ and its null set is the set of linear combinations of elements of $G$ summing to zero in $R$.
        The \emph{regular partial field} $\bU_0$ is a partial field associated with $G = \{1,-1\}$ and $R = \bZ$. For every field $k$, there is a unique tract morphism $\bU_0 \to k$.
        \item The \emph{initial tract} is $\bI := (\{\pm 1\}, \{0, 1+(-1)\})$.
        \item The \emph{Krasner hyperfield} $\bK$ is a tract such that $\bK^\times = \{1\}$ and its null set is $\bN[\bK^\times] \setminus \{1\}$.
        Note that $\epsilon = 1$ in $\bK$, and $\bK$ is the terminal object of the category of tracts.
        \item The \emph{sign hyperfield} $\bS$ is a tract such that $\bS^\times = \{\pm 1\}$ and its null set is the set of zero or sums appearing both $+1$ and $-1$ at least once.
        For every ordered field $k$, a map $k \to \bS$ sending zero to zero, all positive elements in $k$ to $1\in \bS$, and all negative elements in $k$ to $-1\in\bS$ is a tract morphism.
        \item The \emph{tropical hyperfield} $\bT$ is a tract such that $\bT^\times = \bR_{>0}$ and its null set is the set of zero or sums in which the maximum element appears at least twice. 
    \end{enumerate}
\end{example}

\subsection{Antisymmetric $F$-matroids}\label{sec: antisym mat over tracts via bases} %

We define an antisymmetric matroid with coefficients in a tract by introducing the restricted Grassmann--Pl\"{u}cker relations over tracts.

\begin{definition}\label{def: restricted G--P F}
    A {\em restricted Grassmann--Pl\"{u}cker function on $E=\ground{n}$ with coefficients in a tract $F$} is a nontrivial function $\varphi: \cT_n \cup \cA_n \to F$ that satisfies the following.
    \begin{enumerate}[label=\rm(Sym)]
    \item\label{item: almost-transversals F}
    If $A \in \cA_n$ and skew pairs $p,q$ such that $p\subseteq A$ and $q\cap A = \emptyset$, then
    \[
        \varphi(A) = (-1)^{i+j} \varphi(A-p+q)  
    \]
    where $i,j \in [n]$ such that $p=\skewpair{i}$ and $q=\skewpair{j}$.
    \end{enumerate}
    \begin{enumerate}[label=\rm(rGP)]
    \item\label{item: restricted G--P relations F} For $S\in \binom{E}{n+1}$ and $T \in \binom{E}{n-1}$ such that $S$ contains exactly one skew pair and $T$ has no skew pair, 
    \[
        \sum_{x\in S\setminus T}
        (-1)^{\smaller{S}{x} + \smaller{T}{x}}
        \varphi(S-x)
        \varphi(T+x) \in N_F.
        \tag{$\ddagger'$}
        \label{eq: restricted G--P relations F}
    \]
    \end{enumerate}
\end{definition}

Then $\cB = \{B \in \cT_n \cup \cA_n : \varphi(B) \ne 0\}$ satisfies \ref{item: sB1},~\ref{item: sB2}, and~\ref{item: exchange}, and we call a pair $(E,\cB)$ the underlying antisymmetric matroid of~$\varphi$.
Two restricted G--P functions $\varphi$ and $\varphi'$ are \emph{equivalent} if $\varphi'= c \cdot \varphi$ for some $c\in F^\times$.
An {\em antisymmetric $F$-matroids} (or an {\em antisymmetric matroid over $F$}) is an equivalence class $[\varphi]$ of restricted G--P functions with coefficients in $F$.
The \emph{$F$-Lagrangian Grassmannian} $\lag{F}{n}$ is the set of antisymmetric $F$-matroids on $\ground{n}$.
The antisymmetric $\bK$-matroids can be regarded as the antisymmetric matroids.
We call antisymmetric $\bS$-matroids {\em oriented antisymmetric matroids} 
and call antisymmetric $\bT$-matroids {\em valuated antisymmetric matroids}.
In Section~\ref{sec: examples}, we show that those encompass oriented matroids and valuated matroids, respectively.

For a restricted G--P function $\varphi$ with coefficient in a tract $F$ and a tract morphism $f: F\to F'$, the composition $f \circ \varphi$ is a restricted G--P function with coefficients in $F'$.
Therefore, there is a pushforward operation $f_*$ such that for each antisymmetric $F$-matroid $M$, $f_* M$ is an antisymmetric $F'$-matroid.
In particular, if $F' = \bK$, then $f_* M$ is identified with the underlying antisymmetric matroid of $M$.

We will frequently use the following lemma without explicit reference.
\begin{lemma}
    Let $S$ be a subtransversal of size $n-2$ and let $i,j$ be distinct elements in $[n]$ such that $S\cap \skewpair{i} = \emptyset = S\cap \skewpair{j}$.
    Then $i+j \equiv 1+\sum_{z\in \{i,i^*,j,j^*\}} \smaller{(S+\{i,i^*,j,j^*\})}{z} \pmod{2}$.
\end{lemma}

\subsection{Antisymmetric $F$-circuit sets}\label{sec: antisym mat over tracts via circuits} %

We can identify each subset $S$ of $E= \ground{n}$ with the indicator vector $\mathbf{1}_S$ in $\{0,1\}^E = \bK^E$.
Then the circuits of an antisymmetric matroid $M$ on $E$ are vectors in $\bK^E$,
and the axioms~\ref{item: isotropic} and~\ref{item: maximal} in Theorem~\ref{thm: cryptomorphism} can be rephrased as follows.
\begin{enumerate} 
    \item[\rm(Orth)]    $\sum_{i=1}^n(X(i)Y(i^*) + \epsilon X(i^*)Y(i)) \in N_{\bK}$ for all $X,Y \in \cC(M) \subseteq \bK^E$.
    \item[\rm(Max)]     For every $S \subseteq E$ such that $|S|=n+1$ and $S$ contains exactly one skew pair, there is $X\in \cC(M)$ such that $\supp(X) \subseteq S$.
\end{enumerate}
Replacing the Krasner hyperfield $\bK$ with an arbitrary tract $F$ in~\ref{item: isotropic}, we define an antisymmetric $F$-circuit set which is equivalent to an antisymmetric $F$-matroid.

\begin{definition}\label{def: prepared}
    A set $\cC$ of vectors in $F^{E}$ is \emph{prepared} if the following conditions hold.
    \begin{enumerate}[label=\rm(\roman*)]
        \item $\mathbf{0} \notin \cC$.
        \item The support of each vector in $\cC$ contains at most one skew pair.
        \item If $X\in \cC$, then $cX \in \cC$ for all $c\in F^\times$.
        \item For $X,Y\in F^E$, if $\supp(X) \subseteq \supp(Y)$ and $Y\in \cC$, then $X = cY$ for some $c\in F^\times$.
    \end{enumerate}
\end{definition}

\begin{definition}\label{def: antisymmetric F-circuit set}
    An \emph{antisymmetric $F$-circuit set} is a prepared set $\cC$ of vectors in $F^E$ satisfying the next two properties:
    \begin{enumerate}[label=\rm(Orth$'$)]
        \item\label{item: isotropic F} $\sympl{X}{Y} := \sum_{i=1}^{n} (X(i)Y(i^*) + \epsilon X(i^*)Y(i)) \in N_F$ for all $X,Y\in \cC$.
    \end{enumerate}
    \begin{enumerate}[label=\rm(Max$'$)]
        \item\label{item: maximality F} For every $S \subseteq E$ such that $|S|=n+1$ and $S$ contains exactly one skew pair, there is $X\in \cC(M)$ such that $\supp(X) \subseteq S$.
    \end{enumerate}
\end{definition}

\begin{lemma}
    If $\cC$ is an antisymmetric $F$-circuit set, then $\ul{\cC} := \{\ul{X} : X\in \cC\}$ is the set of circuits of an antisymmetric matroid. 
    %
\end{lemma}

The set of circuits of an antisymmetric matroid is identified with an antisymmetric $\bK$-circuit set.
For a field $k$ and a Lagrangian subspace $W$ of $k^E$, let $\cC$ be the set of nonzero vectors $X$ in $W$ such that $\ul{X}$ is minimal and $\ul{X}$ contains at most one skew pair.
Then $\cC$ is an antisymmetric $k$-circuit set.
Conversely, if $\cC' \subseteq k^E$ is an antisymmetric $k$-circuit set, 
then the span of $\cC'$ is isotropic by~\ref{item: isotropic F} and has dimension $n$ by~\ref{item: maximality F}.

\subsection{Basic examples}\label{sec: examples}

We first recall basic examples when $F = \bK$ or $F$ is a field.

\begin{example}
    An antisymmetric $\bK$-matroid is identified with the set of bases of an antisymmetric matroid. An antisymmetric $\bK$-circuit set is identified with the set of circuits of an antisymmetric matroid.
\end{example}

\begin{example}\label{eg: field}
    Let $k$ be a field.
    Then an antisymmetric $k$-matroid is equal to a point in the projective space of dimension $2^{n-2}(4+\binom{n}{2})$ satisfying the restricted G--P relations~\eqref{eq: r G--P}.
    For a Lagrangian subspace $W$ in $k^{\ground{n}}$, if $\cC$ is the set of vectors $X$ in $W\setminus \{\mathbf{0}\}$ such that $\ul{X}$ is minimal and $\ul{X}$ contains at most one skew pair,
    then $\cC$ is an antisymmetric $k$-circuit set.
\end{example}

We give an explicit explanation of why Theorem~\ref{thm: bij} implies Theorem~\ref{thm: parameterization}, which was mentioned in the introduction.
Let $W \in \lag{k}{n}$ and let $M_1$ be an antisymmetric $k$-circuit set obtained by collecting all nonzero vectors $X$ in $W$ such that $\ul{X}$ is minimal and $\ul{X}$ contains at most one skew pair.
Note that $W$ is equal to the span of $M_1$.
Let $M_2 := \Phi(W)$ be an antisymmetric $k$-matroid.
By Proposition~\ref{prop: Lag subsp and antisym mat}, their underlying antisymmetric matroids $\ul{M_1}$ and $\ul{M_2}$ are different, but $\cB(\ul{M_1}) = \cB(\ul{M_2})^*$.
Note that $W^* := \{X^* : X\in W\}$ is also a Lagrangian subspace in $k^E$, where $X^*\in k^E$ such that $X^*(i) = X(i^*)$ for each $i\in E$.
Therefore, the equivalence of antisymmetric $k$-circuit sets and antisymmetric $k$-matroids implies that $\Phi$ is a parameterization of $\lag{k}{n}$ and its image is set-theoretically cut out by the restricted G--P relations~\eqref{eq: r G--P}.

\subsection{Matroids with coefficients}\label{sec: matroids with coefficients}

We show how antisymmetric matroids over tracts generalize matroids over tracts.
We briefly review the theory of matroids with coefficients in tracts by Baker and Bowler~\cite{BB2019}, which extends earlier works of Semple and Whittle~\cite{SW1996} and Dress and Wenzel~\cite{DW1991,DW1992}.

\begin{definition}[\cite{BB2019}]
    Let $F$ be a tract and $0\le r\le n$ be integers.
    A \emph{Grassmann--Pl\"{u}cker function of rank $r$ on $[n]$ with coefficients in $F$} is a function $\psi: \binom{[n]}{r} \to F$ such that $\psi$ is not identically zero and satisfies the \emph{Grassmann--Pl\"{u}cker relations}:
    \begin{align*}
        \sum_{x \in S\setminus T} (-1)^{\smaller{S}{x} + \smaller{T}{x}} \psi(S-x) \psi(T+x) \in N_F
    \end{align*}
    for every $S\in \binom{[n]}{r+1}$ and $T\in \binom{[n]}{r-1}$.
    Two Grassmann--Pl\"{u}cker functions $\psi$ and $\psi'$ are \emph{equivalent} if $\psi' = c \cdot \psi$ for some $c\in F^\times$.
    An \emph{$F$-matroid} (or a \emph{matroid over $F$}) is an equivalence class $[\psi]$ of Grassmann--Pl\"{u}cker functions with coefficients in $F$.

    For a Grassmann--Pl\"{u}cker function $\psi$ of rank $r$ on $[n]$, let $\psi^\perp : \binom{[n]}{n-r} \to F$ be a function such that $\psi^\perp([n]\setminus B) := \sgn(B) \cdot \psi(B)$, where $\sgn(B)$ is the sign of the permutation mapping $i$ to the $i$-th smallest element in $B$ if $i\le r$ and to the $(i-r)$-th smallest element in $[n]\setminus B$ if $i>r$.
    Then $\psi^\perp$ is a Grassmann--Pl\"{u}cker function of rank $n-r$ on $[n]$.
    We call $M^\perp := [\psi^\perp]$ the \emph{dual} of $M = [\psi]$.
\end{definition}

The \emph{$F$-Grassmannian} $\gr{F}{r}{n}$ is the set of $F$-matroids of rank $r$ on $[n]$.
If $F=k$ is a field, then it is the ordinary Grassmannian over $k$.
Note that $\bS$-matroids and $\bT$-matroids are equal to oriented matroids~\cite{BVSWZ1999} and valuated matroids~\cite{DW1992b}, respectively, and the set $\gr{\bT}{r}{n}$ of valuated matroids is called the \emph{Dressian} $\dressian{r}{n}$.

We show that every $F$-matroid naturally induces an antisymmetric $F$-matroid, extending that the Grassmannian $\gr{k}{r}{n}$ is a subset of the Lagrangian Grassmannian $\lag{k}{n}$.

\begin{lemma}\label{lem: F-matroids are antisymmetric F-matroids}
    Let $\psi$ be a Grassmann--Pl\"{u}cker function of rank $r$ on $[n]$ with coefficients in a tract $F$.
    Then a function $\varphi: \cT_n \cup \cA_n \to F$ such that for each $B\in \cT_n \cup \cA_n$, 
    \begin{align*}
        \varphi(B) = 
        \begin{cases}
            \psi(B\cap[n]) \cdot \psi^\perp(B^*\cap [n])  & \text{if $|B\cap[n]| = r$},\\
            0               & \text{otherwise},
        \end{cases}
    \end{align*}
    is a restricted G--P function.
\end{lemma}
\begin{proof}
    We first show that $\varphi$ satisfies~\ref{item: almost-transversals F}.
    Let $A\in \cA_n$ such that $A = (B-i+j) \cup ([n]\setminus B)^*$ for some $B\in \binom{[n]}{r}$, $i\in B$, and $j\in [n]\setminus B$.
    Note that $\skewpair{i} \cap A = \emptyset$, $\skewpair{j} \subseteq A$, and $\sgn(B) \cdot \sgn(B-i+j) = (-1)^m$ where $m := 1+ \sum_{z\in\{i,i^*,j,j^*\}} \smaller{(A+\skewpair{i})}{z}$.
    Then
    \begin{align*}
        \varphi((B-i+j) \cup ([n]\setminus B)^*)
        &=
        \psi(B-i+j) \cdot \psi^\perp([n]\setminus B) \\
        &=
        (-1)^m \cdot \psi^\perp([n]\setminus (B-i+j)) \cdot \psi(B) \\
        &=
        (-1)^m \cdot \varphi(B \cup ([n] \setminus (B-i+j))^*).
    \end{align*}

    Now we claim~\ref{item: restricted G--P relations F}.
    Let $S\in\binom{\ground{n}}{n+1}$ and $T\in\binom{\ground{n}}{n-1}$ such that $S$ contains exactly one skew pair and $T$ has no skew pair. 
    Let $S_1 = S\cap [n]$, $S_2 = S\cap [n]^*$, $T_1 = T\cap [n]$, and $T_2 = T\cap [n]^*$.
    We can assume that either $|S_2| = |T_2| = n-r$ or $|S_1| = |T_1| = r$.
    In the former case, $|S_1| = r+1$ and $|T_1| = r-1$, and thus
    \begin{align*}
        &\sum_{x\in S\setminus T}
        (-1)^{\smaller{S}{x} + \smaller{T}{x}}
        \varphi(S-x)
        \varphi(T+x) \\
        &\hspace{1cm}=
        \psi^\perp(S_2) \psi^\perp(T_2)
        \sum_{x\in S_1\setminus T_1}
        (-1)^{\smaller{S_1}{x} + \smaller{T_1}{x}}
        \psi(S_1-x) \psi(T_1+x)
        \in N_F.
    \end{align*}
    The latter case holds similarly.
\end{proof}

\begin{theorem}\label{thm: F-gr to F-lag}
    Let $F$ be a tract. There is an injective map $\gr{F}{r}{n} \to \lag{F}{n}$ such that the following diagram commutes,
    \begin{center}
        \begin{tikzcd}
            \gr{F}{r}{n} \ar[r] \ar[d] & \lag{F}{n} \ar[d] \\
            \gr{\bK}{r}{n} \ar[r] & \lag{\bK}{n} \\
        \end{tikzcd}
    \end{center}
    \vspace{-0.7cm}
    where the vertical arrows mean taking underlying matroids or underlying antisymmetric matroids.
\end{theorem}

By Theorem~\ref{thm: F-gr to F-lag}, every oriented matroid is an oriented antisymmetric matroid.
Also, every valuated matroid is a valuated antisymmetric matroid, equivalently, 
the Dressian $\dressian{r}{n} = \gr{\bT}{r}{n}$ is a subset of the Lagrangian Grassmannian $\lag{\bT}{n}$ over the tropical hyperfield.

\begin{remark}
    The $F$-matroids have several cryptomorphic definitions~\cite{BB2019,Anderson2019} in terms of circuits and vectors.
    It is straightforward to show that a \emph{dual pair of $F$-signature of a matroid~$M$}~\cite{BB2019} induces an antisymmetric $F$-circuit set such that the circuit set of the underlying antisymmetric matroid is exactly $\cC(\antlift(M))$.
    It extends the diagrams in Section~\ref{sec: antisymmetric matroids circuits} more directly.
\end{remark}

\subsection{Even symmetric matroids with coefficients}\label{sec: orthgonal matroids with coefficients}

In~\cite{JK2023}, the author and Jin introduced \emph{orthogonal matroids over tracts} as a generalization of both matroids over tracts and the Lagrangian orthogonal Grassmannian $\mathrm{OGr}_k(n,2n)$, where an orthogonal matroid is another name for an even symmetric matroid.
Three equivalent definitions of orthogonal $F$-matroids were provided in~\cite{JK2023}, and we review one of them here.
An \emph{orthogonal $F$-matroid} is a set $\cC \subseteq F^{\ground{n}}$ such that $\ul{\cC} = \{\ul{X} : X\in\cC\}$ is the set of circuits of an orthogonal matroid, i.e., an even symmetric matroid, and 
\begin{enumerate}
    \item[\rm(Orth$''$)] $\ortho{X}{Y} := \sum_{i=1}^n(X(i)Y(i^*) + X(i^*)Y(i)) \in N_F$ for every $X,Y \in \cC$.
\end{enumerate}
Recall that by Theorems~\ref{thm: circuits of even symmetric} and~\ref{thm: cryptomorphism}, the circuit set of an orthogonal matroid is the circuit set of an antisymmetric matroid.
Therefore, we deduce the next proposition under assuming Theorem~\ref{thm: bij}, the cryptomorphism between antisymmetric $F$-matroids and antisymmetric $F$-circuit sets.
It generalizes the fact that over a field of characteristic two, every skew-symmetric matrix (with zero diagonals) is a symmetric matrix.

\begin{proposition}
    Let $F$ be a tract with $-1=1$.
    Then an orthogonal $F$-matroid is an antisymmetric $F$-matroid.
\end{proposition}

\begin{example}
    The tropical hyperfield $\bT$ is a tract with $-1=1$.
    Thus, every valuated orthogonal matroid in~\cite{JK2023}, which is equivalent to a valuated delta-matroid in~\cite{Wenzel1993b} and a tropical Wick vector in~\cite{Rincon2012}, is a valuated antisymmetric matroid.
\end{example}

\subsection{Symplectic Dressian and isotropic tropical linear spaces}\label{sec: symplectic dressian}

Balla and Olarte~\cite{BO2023} introduced the symplectic Dressian $\mathrm{SpDr}(r,2n)$ as a tropical counterpart to the symplectic Grassmannian, which is based on the work of De Concini~\cite{DeConcini1979} showing that the symplectic Grassmannian is cut out by the Grassmann--Pl\"{u}cker relations together with certain linear relations.
We show that each tropical symplectic Pl\"{u}cker vector in $\mathrm{SpDr}(n,2n)$ naturally induces a valuated antisymmetric matroid by restricting the domain of the vector to transversals and almost-transversals.
They also investigated isotropic tropical linear spaces, which were first introduced by Rinc\'{o}n~\cite{Rincon2012} to study valuated even delta-matroids.
We show that if we collect certain minimal vectors in an isotropic tropical linear space associated with a rank-$n$ valuated matroid on $2n$ elements, then it forms an antisymmetric $\bT$-circuit set.
It is analogous to a relation between Lagrangian subspaces over a field $k$ and antisymmetric $k$-circuit sets shown in Example~\ref{eg: field}

\emph{Valuated matroids} were first introduced by Dress and Wenzel~\cite{DW1992b}.
They are identified with $\bT$-matroids in the sense of~\cite{BB2019}, which was reviewed in Section~\ref{sec: matroids with coefficients}.
For integers $0\le r\le n$, the \emph{Dressian} $\dressian{r}{n}$ is the set of valuated matroids of rank $r$ on $[n]$.
A \emph{tropical linear space} is a subset of $\bT^{n}$ defined as
\begin{align*}
    L_{\mu} = \left\{ x \in \bT^n : \sum_{i\in S} \mu_{S-i} x_i \in N_{\bT} \text{ for every $S\in \binom{[n]}{r+1}$} \right\}
\end{align*}
for some valuated matroid $[\mu] \in \dressian{r}{n}$.
In the remainder of this section, we regard the ground set of valuated matroids in $\dressian{r}{2n}$ as $\ground{n}$.

\begin{definition}[\cite{BO2023}]
    A \emph{tropical symplectic Pl\"{u}cker vector} is a $\bT$-matroid $[\mu] \in \dressian{r}{2n}$ satisfying the \emph{tropical symplectic relations}:
    \begin{align*}
        \sum_{i\in [n]\setminus(S\cup S^*)} \mu(S+\skewpair{i}) \in N_\bT
    \end{align*}
    for every subset $S$ of $E$ with $|S|=r-2$.
    The \emph{symplectic Dressian} $\mathrm{SpDr}(r,2n)$ is the set of tropical symplectic Pl\"{u}cker vectors in $\dressian{r}{2n}$.
\end{definition}

\begin{definition}[\cite{Rincon2012,BO2023}]
    A tropical linear space $L \subseteq \bT^{\ground{n}}$ is \emph{isotropic} if for all $X,Y\in L$,
    \begin{align*}
        \alpha(X,Y) = \sum_{i=1}^n ( X(i)Y(i^*) + X(i^*)Y(i) ) \in N_{\bT}.
    \end{align*}
\end{definition}

\begin{proposition}[\cite{BO2023}]\label{prop: Balla and Olarte}
    For $n\le 3$ and $[\mu] \in \dressian{n}{2n}$, the valuated matroid $[\mu]$ is a symplectic Pl\"{u}cker vector if and only if the tropical linear space $L_{\mu}$ is isotropic.
\end{proposition}

Proposition~\ref{prop: Balla and Olarte} does not hold for $n\ge 4$ by~{\cite[Examples~B and~C]{BO2023}}.
We show that symplectic Pl\"{u}cker vectors and isotropic tropical linear spaces naturally produce valuated antisymmetric matroids and antisymmetric $\bT$-circuit sets, respectively.
We note that valuated antisymmetric matroids and antisymmetric $\bT$-circuit sets are equivalent by Theorem~\ref{thm: bij} applied to $F=\bT$.

\begin{proposition}
    If $[\mu] \in \mathrm{SpDr}(n,2n)$ is a tropical symplectic Pl\"{u}cker vector, then $[\mu|_{\cT_n\cup \cA_n}]$ is an antisymmetric $\bT$-matroid. 
\end{proposition}
\begin{proof}
    Let $\mu' := \mu|_{\cT_n\cup \cA_n}$, i.e., $\mu' : \cT_n \cup \cA_n \to \bT$ such that $\mu'(B) = \mu(B)$ for all $B\in\cT_n\cup \cA_n$.
    Since $[\mu]\in \dressian{n}{2n}$, $\mu'$ satisfies~\ref{item: restricted G--P relations F}.
    Since $\mu$ satisfies the tropical symplectic relations applied to subtransversals $S$ of size $n-2$, we have $\mu'(S + p) = \mu'(S + q)$ where $p,q$ are distinct skew pairs in $(\ground{n}) - S$. 
    Therefore, $\mu'$ satisfies~\ref{item: almost-transversals F}.
    Thus, it suffices to show that $\mu'$ is nontrivial, i.e., $\mu'(B) \ne 0$ for some $B\in \cT_n \cup \cA_n$.

    Since $[\mu]$ is a valuated matroid, there is an $n$-element subset $X$ of $\ground{n}$ such that $\mu(X) \ne 0$.
    We choose $X$ minimizing $|X\cap X^*|$.
    We may assume that $\frac{1}{2}|X\cap X^*| \ge 2$.
    Let $p$ be a skew pair contained in $X$ and let $S:= X-p \in \binom{\ground{n}}{n-2}$.
    By the tropical symplectic relations applied to $S$, we obtain a skew pair $q \ne p$ such that $\mu(S-p+q) \ne 0$.
    Let $x\in p$ and $y\in q$, and then applying the tropical $3$-term Pl\"{u}cker relation to $S+y$ and $(S-p+q)-y$, we deduce that $\mu(S+y-x) \ne 0$ or $\mu(S+y-x^*)\ne 0$.
    It contradicts the minimality.
\end{proof}

\begin{proposition}
    Let $L_{\mu} \subseteq \bT^{\ground{n}}$ be a tropical linear space associated with a rank-$n$ valuated matroid $[\mu] \in \dressian{n}{2n}$. 
    Let $\cC$ be the set of vectors $X\in L_{\mu} \setminus \{\mathbf{0}\}$ such that $\ul{X}$ is minimal and contains at most one skew pair.
    Then if $L_{\mu}$ is isotropic, then $\cC$ is an antisymmetric $\bT$-circuit set.
\end{proposition}
\begin{proof}
    Clearly, $\cC$ is prepared and satisfies~\ref{item: isotropic F}.
    Because $[\mu]$ has rank $n$, for every $(n+1)$-element subset $S$ of $\ground{n}$, there is a vector $X\in L_{\mu}$ such that $\ul{X} \subseteq S$.
    Hence $\cC$ satisfies~\ref{item: maximality F}.
\end{proof}

\subsection{Oriented gaussoids}\label{sec: oriented gaussoids}

We show that the oriented gaussoids~{\cite[Section~5]{BDKS2019}} are a subclass of the oriented antisymmetric matroids.
We review the definition of oriented gaussoids first.
Let $\Sigma$ be an $n$-by-$n$ symmetric matrix that is real and positive definite.
For $L\subseteq[n]$ and distinct $i,j,k\in [n]\setminus L$, the following equation holds:
\begin{align*}
    \det(\Sigma[L,L]) \det(\Sigma[L+ij,L+ik]) 
    &-
    \det(\Sigma[L+i,L+i]) \det(\Sigma[L+j,L+k]) \\
    &+
    \det(\Sigma[L+i,L+j]) \det(\Sigma[L+i,L+k]) = 0.
\end{align*}
We remark that this equation corresponds to the edge relations for $\begin{bmatrix} I_n \,|\, \Sigma \end{bmatrix}$.
For $I\subseteq [n]$, let $p_I$ be an unknown representing the $I\times I$ principal minor of $\Sigma$.
For $K\subseteq [n]$ and distinct $i,j\in [n]\setminus K$, let $a_{ij|K}$ be an unknown representing the $(K+i)\times (K+j)$ almost-principal minor of $\Sigma$.
Then $a_{ij|K} = a_{ji|K}$.
Let $\cP\cA$ be the set of all such unkowns $p_I$ and $a_{ij|K}$.
An oriented gaussoid is defined as follows, which can be regarded as a counterpart of a positive definite symmetric matrix over the sign hyperfield~$\bS$.

\begin{definition}[\cite{BDKS2019}]
    An \emph{oriented gaussoid} is a map $\varphi: \cP\cA \to \bS$ such that
    \begin{enumerate}[label=\rm(\roman*)]
        \item $\varphi(p_I) = 1$ for every $I\subseteq [n]$, and
        \item for $L\subseteq[n]$ and distinct $i,j,k\in [n]\setminus L$, 
        \begin{align*}
            p_L a_{jk|L+i} - p_{L+i} a_{jk|L} + p_{ij|L} a_{ik|L} \in N_{\bS}.
        \end{align*}
    \end{enumerate}
\end{definition}

\begin{example}
    Let $\Sigma := \begin{bmatrix}
        1       & 1/2   & 1/4 \\
        1/2     & 1     & 1/4 \\
        1/4     & 1/4   & 1   \\
    \end{bmatrix}$ be a real symmetric matrix.
    Then it is positive definite.
    $\det(\Sigma[12,13]) = 1/8$, $\det(\Sigma[12,23]) = -1/8$, and $\det(\Sigma[13,23]) = 7/16$.
    Note that there is a tract morphism $\bR \to \bS$ such that all positives go to $+1$ and all negatives go to $-1$.
    This implies that a map $\varphi:\cP\cA \to \bS$ such that $\varphi^{-1}(+1) = \cP\cA \setminus \{a_{13|2}\}$ and $\varphi^{-1}(-1) = \{a_{13|2}\}$ is an oriented gaussoid.
\end{example}

Recall that $\det(\Sigma[X,Y]) = (-1)^{\sum_{i\in X} (i+\smallereq{X}{i})} \det(\Lambda[n, [n]-X+Y^*])$, where $\Lambda$ is an $n\times (\ground{n})$ matrix $\begin{bmatrix} I_n \,|\, \Sigma \end{bmatrix}$.
Therefore, we can define an oriented gaussoid alternatively in terms of transversals and almost-transversals.

\begin{definition}\label{def: oriented gaussoid alt}
    An \emph{oriented gaussoid} is a map $\varphi: \cT_n \cup \cA_n \to \bS$ such that
    \begin{enumerate}[label=\rm(\roman*)]
        \item $\varphi([n]-X+X^*) = (-1)^{\sum_{i\in X} (i+\smallereq{X}{i})}$ for every $X\subseteq [n]$, and
        \item\label{item: og edge} for every transversal $S\cup \{a,b,c\}$ with $S\cap \{a,b,c\} = \emptyset$,
        \begin{align*}
            (-1)^{\smaller{L}{a}} \varphi(S \cup abc) \varphi(S \cup bb^*c^*)
            &+
            (-1)^{\smaller{L}{b^*}}
            \varphi(S \cup abc^*) \varphi(S \cup bb^*c) \\
            &+
            (-1)^{\smaller{L}{c^*}}
            \varphi(S \cup abb^*) \varphi(S \cup bcc^*)
            \in N_{\bS}
        \end{align*}
        where $L = \{a,c,b^*,c^*\}$.
    \end{enumerate}
\end{definition}

Note that the formula in Definition~\ref{def: oriented gaussoid alt}\ref{item: og edge} is the edge relations over $\bS$.
Hence $\varphi^{-1}(0)$ is a gaussoid for every oriented gaussoid $\varphi:\cT_n\cup \cA_n \to \bS$.
Also, the following proposition is straightforward.

\begin{proposition}\label{prop: ori ant mat to ori gaussoid}
    Let $[\varphi]$ be an oriented antisymmetric matroid on $\ground{n}$.
    If $\varphi([n]-X+X^*) = (-1)^{\sum_{i\in X} (i+\smallereq{X}{i})}$ for every $X\subseteq [n]$, then $\varphi$ is an oriented gaussoid.
\end{proposition}

An oriented gaussoid $\varphi:\cP\cA \to \bS$ is \emph{realizable} if there is a positive definite real symmetric matrix~$\Sigma$ such that $\varphi(a_{ij|K}) \in \{0,\pm1\} = \bS$ equals the sign of $\det(\Sigma[K+i,K+j])$ for all $K\subseteq [n]$ and distinct $i,j\in [n]\setminus K$.
A \emph{positive gaussoid} is an oriented gaussoid $\varphi:\cP\cA \to \bS$ such that $\varphi(a_{ij|K}) \in \{0,1\}$ for every $a_{ij|K}$.
Equivalently, a positive gaussoid is an oriented gaussoid $\varphi:\cT_n\cup \cA_n \to \bS$ such that $\varphi([n]-X+Y^*) = 0$ or $(-1)^{\sum_{i\in X} (i+\smallereq{X}{i})}$ for all $[n]-X+Y^*\in \cA_n$.

Ardila, Rinc\'{o}n, and Williams~\cite{ARW2017} showed that every positively orientable matroid is representable over the real field~$\bR$.
Boege et al.~\cite{BDKS2019} showed an analogous result for positive gaussoids.

\begin{theorem}[{\cite[Theorem~4]{BDKS2019}}]
    Every positive gaussoid is realizable.
\end{theorem}

By the previous theorem together with Proposition~\ref{prop: ori ant mat to ori gaussoid}, we deduce the following.

\begin{corollary}\label{cor: realizable positively oriented antisym mat}
    Let $M = [\varphi]$ be an oriented antisymmetric matroid.
    Suppose that 
    \begin{enumerate}[label=\rm(\roman*)]
        \item $\varphi([n]-X+X^*) = (-1)^{\sum_{i\in X} (i+\smallereq{X}{i})}$ for every $X\subseteq [n]$, and 
        \item $\varphi([n]-X+Y^*) = 0$ or $(-1)^{\sum_{i\in X} (i+\smallereq{X}{i})}$ for all $[n]-X+Y^*\in \cA_n$.
    \end{enumerate}
    Then there is an antisymmetric matroid~$N$ over~$\bR$ such that $f_* N = M$, where $f:\bR \to \bS$ is a tract morphism sending all positives to $1$ and all negatives to $-1$.
\end{corollary}

It is interesting to show whether the condition~(i) can be weakened to that $\varphi([n]-X+X^*) = 0$ or $(-1)^{\sum_{i\in X} (i+\smallereq{X}{i})}$ for every $X\subseteq [n]$.

\section{Cryptomorphism}\label{sec: cryptomorphism}

In Section~\ref{sec: antisymmetric matroids over tracts}, we introduced two concepts generalizing Lagrangian subspaces in the standard symplectic vector space.
First, an antisymmetric $F$-matroid is defined by the restricted G--P relations over tracts.
Second, an antisymmetric $F$-circuit set is defined as a maximal set of vectors that are orthogonal to each other, subject to a certain condition on their supports.
We show that those two notions are equivalent.

\begin{theorem-bij}
    There is a natural bijection between antisymmetric $F$-matroids and antisymmetric $F$-circuit sets.
\end{theorem-bij}

A proof of Theorem~\ref{thm: bij} is provided in Section~\ref{sec: equiv}.
Theorem~\ref{thm: bij} generalizes not only the cryptomorphism on antisymmetric matroids but also the parametrization of the Lagrangian Grassmannian into the projective space of dimension $2^{n-2}(n+\binom{n}{2})$, as explained in Section~\ref{sec: examples}.
In Section~\ref{sec: B2C}, we construct an antisymmetric $F$-circuit set from an antisymmetric $F$-matroid.
In Section~\ref{sec: C2B}, we oppositely build an antisymmetric $F$-matroid from an antisymmetric $F$-circuit set.
Those two constructions are evidently inverse to one another, and hence we deduce Theorem~\ref{thm: bij}.

\subsection{Constructing an antisymmetric $F$-circuit set}\label{sec: B2C}
In this subsection, we let $\varphi : \cT_n\cup \cA_n \to F$ be a restricted G--P function on $E:= \ground{n}$ with coefficients in a tract $F$.
We denote the underlying antisymmetric matroid of $\varphi$ by $M = (E,\cB)$.
The goal is to construct an antisymmetric $F$-circuit set from $\varphi$.

Let $S\subseteq E$ be a subset of size $n+1$ which contains exactly one skew pair, say $\skewpair{x}$.
Suppose that $S-x$ or $S-x^*$ is a basis.
Let $X_S \in F^E$ be a vector defined as follows:
\begin{itemize}
    \item $\supp(X_S)\subseteq S$ and 
    \item $X_S(y) = (-1)^{\chi(y)+\smaller{S}{y}} \varphi(S-y)$ for each $y\in S$.
\end{itemize}
Let $\cC$ be the set of all $c X_S$ such that $c\in F^\times$ and $S=B+x^*$ with $B\in \cB\cap \cT_n$ and $x\in B$.

\begin{lemma}\label{lem: supports of Xs}
    $\ul{\cC} = \{\ul{X} : X\in \cC\}$ is the set of circuits of $M$.
\end{lemma}
\begin{proof}
    Let $B$ be a transversal basis of $M$ and let $e\in B^*$.
    Then by Lemma~\ref{lem: fundamental circuit}, the support of $X_{B+e}$ is the fundamental circuit with respect to $B$ and $e$.
    By Lemma~\ref{lem: basis from circuit}, every circuit is a fundamental circuit with respect to some transversal basis and element, and thus $\ul{\cC} = \cC(M)$.
\end{proof}

\begin{theorem}\label{thm: B2C}
    $\cC$ is an antisymmetric $F$-circuit set.
\end{theorem}
\begin{proof}
    It is clear that $\cC$ satisfies Definition~\ref{def: prepared}(i)--(iii).
    By Lemma~\ref{lem: supports of Xs}, it satisfies~\ref{item: maximality F}.
    Next we examine~\ref{item: isotropic F}, i.e., $\sympl{X}{Y} \in N_F$ for all $X,Y \in \cC$.

    Fix $X,Y\in \cC$.
    Let $S,T$ be subsets of size $n+1$ in $E$ such that
    $\ul{X} \subseteq S$, $\ul{Y} \subseteq T$, and for some $x,y\in E$, 
    $S-x^*$ and $T-y^*$ are transversal bases.
    Applying the restricted G--P relation \ref{item: restricted G--P relations F} to $S$ and $T':=T\setminus \skewpair{y}$, we have
    \[
        \sum_{z \in S\setminus T'} (-1)^{\smaller{S}{z}+\smaller{T'}{z}} \cdot \varphi(S-z) \cdot \varphi(T'+z) \in N_F.
    \]
    Note that $S\setminus T' = (S\cap \skewpair{y}) \cup (S\setminus T)$.
    Let $m := \smaller{T'}{y} + \smaller{T'}{y^*}$. 
    Then for each $z\in\skewpair{y}$, we have
    $\smaller{T}{z^*} = \smaller{T'}{z^*} + \chi(z^*) \equiv \smaller{T'}{z} + m + \chi(z^*) \pmod{2}$ and thus $Y(z^*) = (-1)^{\smaller{T}{z^*}+\chi(z^*)} \varphi(T-z^*) = (-1)^{\smaller{T'}{z}+m} \varphi(T'+z)$.

    For $z \in S\setminus T$, we have $1 + \sum_{w\in\{y,y^*,z,z^*\}} \smaller{(T+z)}{w} \equiv \smaller{T'}{z} + \smaller{T}{z^*} + \chi(z^*) + m \pmod{2}$.
    Then by~\ref{item: almost-transversals F}, $Y(z^*) = (-1)^{\chi(z^*)+\smaller{T}{z^*}} \varphi(T-z^*) = (-1)^{\smaller{T'}{z}+m}\varphi(T'+z)$. 
    Therefore,
    \begin{align*}
        \sympl{X}{Y}
        &=
        \sum_{z\in S\setminus T'}
        (-1)^{\chi(z)}X(z)Y(z^*) \\
        &=
        (-1)^{m}
        \sum_{z \in S\setminus T'} (-1)^{\smaller{S}{z}+\smaller{T'}{z}} \cdot \varphi(S-z) \cdot \varphi(T'+z) \in N_F.
    \end{align*}

    Finally, we check Definition~\ref{def: prepared}(iv).
    Let $X,Y \in \cC$ such that $\ul{X} \subseteq \ul{Y}$.
    By Lemma~\ref{lem: supports of Xs}, $\ul{X} = \ul{Y}$.
    Choose $e\in \ul{X}$ such that $\ul{X}-e$ is a subtransversal.
    For each $f\in \ul{X}-e$, $M$ has a circuit $C$ such that $\ul{X} \cap C^* = \{e,f\}$ by Lemma~\ref{lem: circuits with two intersection}, and let $Z\in\cC$ such that $\ul{Z} = C$.
    Then by~\ref{item: isotropic F}, $\frac{X(f)}{X(e)} = (-1)^{1+\chi(e)+\chi(f)} \frac{Z(e^*)}{Z(f^*)} = \frac{Y(f)}{Y(e)}$.
    Since $f$ is arbitrary, we conclude that $X = cY$ for some $c\in F^\times$.
\end{proof}

\subsection{Constructing an antisymmetric $F$-matroid}\label{sec: C2B}
Let $\cC$ be an antisymmetric $F$-circuit set on $E = \ground{n}$, and let $M$ be its underlying antisymmetric matroid, i.e., $\cC(M) = \ul{\cC}$.
We will construct a restricted G--P function by approaching a reverse step of Section~\ref{sec: B2C}.

\begin{definition}\label{def: gamma}
    Let $B_1$ and $B_2$ be bases of $M$ such that $S:= B_1\cup B_2$ has exactly one skew pair.
    We denote by $\{x\}=S\setminus B_1$ and $\{y\}=S\setminus B_2$, and we define 
    \[
        \gamma(B_1,B_2) := (-1)^{\chi(x)+\chi(y)+\smaller{S}{x}+\smaller{S}{y}} \frac{X(y)}{X(x)},
    \]
    where $X$ is a vector in $\cC$ such that $\ul{X} \subseteq S$.
\end{definition}

\begin{definition}\label{def: basis graph}
    The \emph{basis graph} of $M$ is a graph $G_M$ on $\cB(M)$ such that two vertices $B$ and $B'$ are adjacent if and only if $|B\setminus B'|=1$ and at least one of $B$ and $B'$ is a transversal.
\end{definition}

\begin{lemma}
    The basis graph $G_M$ is connected.
\end{lemma}
\begin{proof}
    It is an immediate consequence of Lemma~\ref{lem: 3-term basis exchange} and~\ref{lem: cGm}.
\end{proof}

One natural candidate of a restricted G--P function $\varphi : \cT_n\cup \cA_n \to F$ can be constructed as follows.
\begin{enumerate}[label=\rm(\roman*)]
    \item Fix a basis $B_0\in \cB(M)$ and let $\varphi(B_0) = 1 \in F^\times$.
    \item For each $B\in \cB(M)$, let $\varphi(B) = \prod_{i=0}^{k-1}\gamma(B_i,B_{i+1}) \in F^\times$ where $B_0B_1\ldots B_k$ be a path from $B_0$ to $B_k:=B$ in the basis graph $G_M$.
    \item For $B\in (\cT_n\cup \cA_n)\setminus \cB(M)$, let $\varphi(B) = 0$.
\end{enumerate}
By similar proof of Theorem~\ref{thm: B2C}, we can show that $\varphi$ is a restricted G--P function on $E$ with coefficients in $F$ (Theorem~\ref{thm: C2B}).
However, the hardest part is to show that $\varphi$ is well defined, i.e., for different paths $P = B_0\ldots B_k$ and $P' = B_0'\ldots B_\ell'$ from $B_0=B_0'$ to $B_k=B_{\ell'}=B$, we should prove that $\prod_{i=0}^{k-1}\gamma(B_i,B_{i+1}) = \prod_{i=0}^{\ell-1}\gamma(B_i',B_{i+1}')$.
Henceforth, we devoted most of the subsection to show that $\varphi$ is well defined using the Homotopy Theorem (Theorem~\ref{thm: homotopy}).

\begin{lemma}\label{lem: gamma 2cycle}
    $\gamma(B_1,B_2) = \gamma(B_2,B_1)^{-1}$
    for each $B_1B_2 \in E(G_M)$.
\end{lemma}

\begin{lemma}\label{lem: gamma 2almost}
    Let $B$ be a transversal basis, and let $A$ and $A'$ be distinct almost-transversal bases such that $A=B+x^*-y$ and $A'=B-x+y^*$ for some $x,y\in B$.
    Then
    \[
        \gamma(B,A) = (-1)^m \gamma(B,A'),  
    \]
    where $m:= 1+\sum_{z\in\{x,x^*,y,y^*\}} \smaller{(B+x^*+y^*)}{z}$. %
\end{lemma}
\begin{proof}
    Let $S = B+x^*$, $T = B+y^*$, and $U = B+x^*+y^*$. 
    Let $X,Y$ be vectors in $\cC$ such that $\supp(X)\subseteq S$ and $\supp(Y)\subseteq T$.
    Then $\sympl{X}{Y} = (-1)^{\chi(x^*)}X(x^*)Y(x) + (-1)^{\chi(y)}X(y)Y(y^*) \in N_F$.
    Note that $\smaller{S}{x^*}+\smaller{T}{y^*} = \smaller{U}{x^*}+\smaller{U}{y^*} -1$, $\smaller{S}{y} = \smaller{U}{y} - \chi(y)$, and $\smaller{T}{x} = \smaller{U}{x} - \chi(x)$.
    Therefore,
    \[
        \gamma(B,A)\gamma(B,A')^{-1}
        =
        (-1)^{\smaller{S}{x^*}+\smaller{T}{y^*}+\smaller{S}{y}+\smaller{T}{x}} \frac{X(y)Y(y^*)}{X(x^*)Y(x)}
        =
        (-1)^m.
    \]
\end{proof}

\begin{lemma}\label{lem: gamma 4cycle}
    Let $B_1B_2B_3B_4B_1$ be a $4$-cycle in $G_M$.
    Then $\prod_{i=1}^{4}\gamma(B_i,B_{i+1}) = 1$, where $B_5:=B_1$.
\end{lemma}
\begin{proof}
    If $B_i$ is an almost-transversal, then $B_{i-1}$ and $B_{i+1}$ are transversals.
    Hence, by relabelling, we may assume that $B_1$ and $B_3$ are transversals.
    As $B_1B_2B_3B_4B_1$ is a $4$-cycle in $G_M$, $|B_2\setminus B_1| = |B_4\setminus B_1| = 1$ and $|B_3\setminus B_1| = 2$.
    Then $B_3 = B_1 \symdiff \{x,x^*,y,y^*\}$ for some $x,y \in B$.


    \noindent\textbf{Case I.}
    Both $B_2$ and $B_4$ are almost-transversals.
    Then by symmetry, we may assume that $B_2 = B_1+x^*-y = B_3+x-y^*$ and $B_4 = B_1-x+y^* = B_3-x^*+y$.
    Therefore, $\gamma(B_1,B_2)\gamma(B_2,B_3) = \gamma(B_1,B_4)\gamma(B_4,B_3)$ by Lemma~\ref{lem: gamma 2almost}.


    \noindent\textbf{Case II.}
    $B_2$ is a transversal and $B_4$ is an almost-transversal.
    Then by symmetry, we may assume that $B_2 = B_1\symdiff\skewpair{x}$.
    Then $B_4$ is either $B_1-y+x^*$ or $B_1-x+y^*$.
    By Case I, we can assume that $B_4 = B_1-y+x^*$.
    We denote by $B_4' := B_1-x+y^*$.

    Let $S:= B_1+x^*$ and $T:= B_3+y$. %
    Let $X,Y\in \cC$ be vectors such that $\ul{X}\subseteq S$ and $\ul{Y}\subseteq T$.
    Then 
    \begin{align*}
        \gamma(B_1,B_2)\gamma(B_2,B_3)
        &=
        (-1)^{\smaller{S}{x}+\smaller{S}{x^*}+\smaller{T}{y}+\smaller{T}{y^*}} \frac{X(x)}{X(x^*)} \frac{Y(y)}{Y(y^*)}, \\
        \gamma(B_1,B_4)\gamma(B_4',B_3)
        &=
        (-1)^{\smaller{S}{y}+\smaller{S}{x^*}+\smaller{T}{y}+\smaller{T}{x^*}} \frac{X(y)}{X(x^*)} \frac{Y(y)}{Y(x^*)}.
    \end{align*}
    As $\sympl{X}{Y}\in N_F$, we have $\frac{X(x)}{Y(y^*)} = (-1)^{\chi(x)+\chi(y)+1} \frac{X(y)}{Y(x^*)}$.
    Hence $\gamma(B_1,B_2)\gamma(B_2,B_3) \gamma(B_3,B_4')\gamma(B_4,B_1) = (-1)^{1+\chi(x)+\chi(y)+m}$ where $m:= \smaller{S}{x} + \smaller{S}{y} + \smaller{T}{y^*} + \smaller{T}{x^*}$.
    If $U:= S+y^* = T+x$, then $m = \sum_{e\in\{x,x^*,y,y^*\}} \smaller{U}{e} - (1+\chi(x^*)+\chi(y))$.
    Therefore, by Lemma~\ref{lem: gamma 2almost}, we obtain the desired equality.


    \noindent\textbf{Case III.}
    Both $B_2$ and $B_4$ are transversals.
    Then by symmetry, $B_2 = B_1\symdiff\skewpair{x} = B_3\symdiff\skewpair{y}$ and $B_4 = B_1\symdiff\skewpair{y} = B_3\symdiff\skewpair{x}$.
    If $B_1-x+y^*$ is a basis, then applying Case II twice, we can deduce the desired equality.
    Therefore, we can assume that $B_1-x+y^*$ is not a basis.

    Let $X,Y,Z,W\in \cC$ be vectors such that $\ul{X}\subseteq B_1\cup B_2$, $\ul{Y}\subseteq B_2\cup B_3$, $\ul{Z}\subseteq B_1\cup B_4$, and $\ul{W}\subseteq B_4\cup B_3$.
    Because neither $B_1-x+y^*$ nor $B_1+x^*-y$ is a basis, $X(y)=Y(x^*)=Z(x)=W(y^*) = 0$.
    Then $\ul{W} \subseteq (B_4\cup B_3) - y^* \subseteq B_1 \cup B_2$ and thus $\ul{X} = \ul{W}$ by Lemma~\ref{lem: fundamental circuit}.
    Similarly, $\ul{Y} = \ul{Z}$.
    Hence $\frac{X(x)}{X(x^*)} = \frac{W(x)}{W(x^*)}$ and $\frac{Y(y)}{Y(y^*)} = \frac{Z(y)}{Z(y^*)}$.
    Therefore, $\gamma(B_1,B_2)\gamma(B_2,B_3) = \gamma(B_1,B_4)\gamma(B_4,B_3)$.
\end{proof}

\begin{definition}
    For two transversal bases $B_1$ and $B_2$ of $M$ such that $|B_1\setminus B_2| = 2$,
    let
    \[
        \gamma(B_1,B_2) := \gamma(B_1,B)\gamma(B,B_2)
    \]
    where $B$ is an arbitrary basis such that $BB_1, BB_2 \in E(G_M)$. 
    It is well defined by Lemma~\ref{lem: gamma 4cycle}.
\end{definition}

Recall that the transversal basis graph $\cG_M$ is a graph on $\cB(M)\cap \cT_n$ together with weights $\eta(BB') = |B\setminus B'| \in \{1,2\}$ on its edges $BB'$. %
We say two cycles $C_1$ and $C_2$ in $\cG_M$ are \emph{$4$-homotopic}, denoted by $C_1 \simeq C_2$, if there is a sequence of cycles $D_1,\ldots,D_k$ such that 
\begin{itemize}
    \item $D_1 = C_1$,
    \item $D_k = C_2$, and
    \item each $D_{i+1}$ is obtained from $D_i$ either 
    \begin{enumerate}[label=\rm(\alph*)]
        \item by replacing two edges $B_1B_2$ and $B_2B_3$ of weight $1$ with an edge $B_1B_3$ of weight $2$, or 
        \item by replacing an edge $B_1B_3$ of weight $2$ with two edges $B_1B_2$ and $B_2B_3$ of weight $1$.
    \end{enumerate}
\end{itemize}

\begin{lemma}\label{lem: irreducible 6cycle}
    Let $C$ be an irreducible cycle of weight $6$ in $\cG_M$.
    Then $C \simeq B_1B_2 \ldots B_kB_1$ such that 
    \begin{enumerate}[label=\rm(\roman*)]
        \item $k=3$ and $\eta(B_1B_2) = \eta(B_2B_3) = \eta(B_3B_1) = 2$, or
        \item $k=4$ and $\eta(B_1B_2) = \eta(B_3B_4) = 1$ and $\eta(B_2B_3) = \eta(B_4B_1) = 2$.
    \end{enumerate}
\end{lemma}
\begin{proof}
    Let $C'$ be an irreducible cycle $4$-homotopic to $C$, which maximizes the number $|E(C') \cap \eta^{-1}(2)|$ of edges $e \in E(C')$ such that $\eta(e)=2$.

    \begin{claim}
        For each $B \in V(C')$, there are no three consecutive vertices $D_1,D_2,D_3$ in $C'$ such that $\eta(D_1D_2) = \eta(D_2D_3) = 1$, $\dist_M(B,D_1) = \dist_M(B,D_3) = 2$, and $\dist_M(B,D_2) = 3$.
    \end{claim}
    \begin{proof}
        Suppose that such vertices $D_1, D_2, D_3$ exist.
        If $D_1D_3\in E(\cG_M)$, then it contradicts our choice of $C'$.
        Thus, $D_1D_3\notin E(\cG_M)$ and hence by Lemma~\ref{lem: 3-term basis exchange}, $M$ has a transversal basis $D' \ne D_2$ such that $|D'\setminus D_i| = 1$ for $i\in\{1,3\}$.
        Then $|D'\setminus B| = 1$.
        Hence $C'$ is generated by three cycles of weight $4$, a contradiction.
    \end{proof}

    By the claim, we can easily deduce that $C'=B_1\ldots B_kB_1$ satisfies either (i) or (ii).
\end{proof}

\begin{lemma}\label{lem: gamma 6cycle}
    Let $C = B_1 B_2 \ldots B_k B_1$ be a cycle of weight $6$ in $\cG_M$.
    Then $\prod_{i=1}^{k} \gamma(B_i,B_{i+1}) = 1$, where $B_{k+1} := B_1$.
\end{lemma}
\begin{proof}
    We may assume that $C$ is irreducible by Lemma~\ref{lem: gamma 4cycle}.
    Then by Lemma~\ref{lem: irreducible 6cycle} and rotational symmetry, we may assume that either
    \begin{enumerate}[label=\rm(\roman*)]
        \item $k=3$ and $\eta(B_iB_{i+1})=2$ for each $1\le i\le 3$, or
        \item $k=4$ and $\eta(B_1B_2) = \eta(B_3B_4) = 1$ and $\eta(B_2B_3) = \eta(B_4B_1) = 2$.
    \end{enumerate}

    In the case (i), $B_i = T\symdiff\skewpair{x_i}$ for some transversal $T$ and elements $x_1,x_2,x_3\in T$.
    Then $A_i := T+x_i^*-x_{i+1}$ with $1\le i\le 3$ are bases, where we read the subscripts modulo $3$, and $B_1A_1B_2A_2B_3A_3B_1$ is a $6$-cycle in $G_M$; see Figure~\ref{fig: 6-cycle}(left).
    Since $C$ is irreducible, $T$ is not a basis of~$M$.

    Let $X_1,X_2,X_3,Y_1,Y_2,Y_3 \in \cC$ be vectors such that $\ul{X_i} \subseteq B_i \cup A_i = T + x_i$ and $\ul{Y_i} \subseteq A_i \cup B_{i+1}$.
    Because $T$ is not a basis, by Lemma~\ref{lem: fundamental circuit}, $\ul{X_i} \subseteq T$ for each $1\le i\le 3$ and thus $\ul{X_1} = \ul{X_2} = \ul{X_3}$.
    By multiplying elements in $F^\times$, we can assume that $X_1 = X_2 = X_3 =: X$.

    Because $\sympl{X}{Y_i} \in N_F$, we have $\frac{Y_i(x_{i}^*)}{Y_i(x_{i+1}^*)} = (-1)^{1+\chi(x_i)+\chi(x_{i+1})} \frac{X(x_{i+1})}{X(x_{i})}$.
    Note that $\sum_{i=1}^{3} \big( \smaller{(T+x_i^*)}{x_i} + \smaller{(T+x_i^*)}{x_{i+1}} + \smaller{(T\symdiff\skewpair{x_{i+1}}+x_i^*)}{x_{i+1}^*} + \smaller{(T\symdiff\skewpair{x_{i+1}}+x_i^*)}{x_i^*} \big) \equiv 1 \pmod{2}$.
    Therefore,
    \begin{align*}
        \gamma(B_1,B_2) \gamma(B_2,B_3) \gamma(B_3,B_1)
        =
        -
        \frac{X(x_2)}{X(x_1)} \frac{Y_1(x_1^*)}{Y_1(x_2^*)}
        \frac{X(x_3)}{X(x_2)} \frac{Y_2(x_2^*)}{Y_2(x_3^*)}
        \frac{X(x_1)}{X(x_3)} \frac{Y_3(x_3^*)}{Y_3(x_1^*)}
        =
        1.
    \end{align*}

    Now, we prove the case (ii).
    For some $x_1,x_2,x_3 \in B_1$, we have $B_2 = B_1\symdiff\skewpair{x_1}$, $B_3 = B_1\symdiff\{x_1,x_1^*,x_2,x_2^*,x_3,x_3^*\}$, and $B_4 = B_1\symdiff\{x_2,x_2^*,x_3,x_3^*\}$; see Figure~\ref{fig: 6-cycle}(right).

    Suppose that $D := B_1 \symdiff\{x_1,x_1^*,x_2,x_2^*\}$ is a basis.
    Then $B_1 \symdiff\skewpair{x_2}$ is not a basis since $C$ is irreducible.
    By Lemma~\ref{lem: 3-term basis exchange}, $DB_1$ and $DB_4$ are edges of weight $2$ in $\cG_M$.
    Then by Lemma~\ref{lem: gamma 4cycle} and the case (i), we have $\prod_{i=1}^{k} \gamma(B_i,B_{i+1}) = 1$.
    Thus, we can assume that $B_1 \symdiff\{x_1,x_1^*,x_2,x_2^*\}$ is not a basis.
    Similarly, we can assume that none of $B_1 \symdiff\{x_1,x_1^*,x_3,x_3^*\}$, $B_1\symdiff\skewpair{x_2}$, and $B_1\symdiff\skewpair{x_3}$ is a basis.
    Then for each $(i,j)\in [3]^2 \setminus \{(1,1),(2,3),(3,2)\}$, neither $B_1-x_i+x_j^*$ nor $B_3+x_i-x_j^*$ is a basis.

    Let $X_1,X_2,X_3,Y_1,Y_2,Y_3 \in \cC$ be vectors such that $\ul{X_1} \subseteq B_1+x_1$, $\ul{X_2} \subseteq B_2+x_2^*$, $\ul{X_3} \subseteq B_3+x_2$, $\ul{Y_1} \subseteq B_3+x_1^*$, $\ul{Y_2} \subseteq B_4+x_2$, and $\ul{Y_3} \subseteq B_1+x_2^*$. 
    Because of the previous observations on non-bases of $M$, we have that 
    \begin{align*}
        \sympl{X_1}{Y_1} 
        &= 
        (-1)^{\chi(x_1)}X_1(x_1)Y_1(x_1^*) + (-1)^{\chi(x_1^*)}X_1(x_1^*)Y_1(x_1)
        \in N_F, \\
        \sympl{X_2}{Y_2} 
        &= 
        (-1)^{\chi(x_2^*)}X_2(x_2^*)Y_2(x_2) + (-1)^{\chi(x_3)}X_2(x_3)Y_2(x_3^*)
        \in N_F, \\
        \sympl{X_3}{Y_3} 
        &= 
        (-1)^{\chi(x_2)}X_3(x_2)Y_3(x_2^*) + (-1)^{\chi(x_3^*)}X_3(x_3^*)Y_3(x_3)
        \in N_F.
    \end{align*}
    Therefore, $\gamma(B_1,B_2)\gamma(B_2,B_3) = \gamma(B_1,B_4)\gamma(B_4,B_3)$.
\end{proof}

\begin{figure}
    \centering
    \begin{tikzpicture}
        \node[vo] (v1) at (-90+60*0:1.3/2) {};
        \node[v] (v2) at (-90+60*1:1.3) {};
        \node[vo] (v3) at (-90+60*2:1.3/2) {};
        \node[v] (v4) at (-90+60*3:1.3) {};
        \node[vo] (v5) at (-90+60*4:1.3/2) {};
        \node[v] (v6) at (-90+60*5:1.3) {};

        \draw (v2)--(v4)--(v6)--(v2);

        \begin{scope}[xshift=-1.8cm, yshift=-1.15cm]
            \node[box1s] () at (-0.3+0.3*0, 0.15) {};
            \node[box2s] () at (-0.3+0.3*1, 0.15) {};
            \node[box2s] () at (-0.3+0.3*2, 0.15) {};

            \node[box2s] () at (-0.3+0.3*0, -0.15) {};
            \node[box1s] () at (-0.3+0.3*1, -0.15) {};
            \node[box1s] () at (-0.3+0.3*2, -0.15) {};

            \node () at (0,-0.6) {\scriptsize{$B_1$}};

            \node () at (-0.3+0.3*0, 0.15) {\scriptsize{$x_1^*$}};
            \node () at (-0.3+0.3*0, -0.15) {\scriptsize{$x_1$}};
            \node () at (-0.3+0.3*1, -0.15) {\scriptsize{$x_2$}};
            \node () at (-0.3+0.3*2, -0.15) {\scriptsize{$x_3$}};
        \end{scope}

        \begin{scope}[xshift=-1.35cm, yshift=0.325cm]
            \node[box1s] () at (-0.3+0.3*0, 0.15) {};
            \node[box2s] () at (-0.3+0.3*1, 0.15) {};
            \node[box2s] () at (-0.3+0.3*2, 0.15) {};

            \node[box1s] () at (-0.3+0.3*0, -0.15) {};
            \node[box2s] () at (-0.3+0.3*1, -0.15) {};
            \node[box1s] () at (-0.3+0.3*2, -0.15) {};

            \node () at (0,0.55) {\scriptsize{$A_1$}};
        \end{scope}

        \begin{scope}[xshift=0cm, yshift=1.85cm]
            \node[box2s] () at (-0.3+0.3*0, 0.15) {};
            \node[box1s] () at (-0.3+0.3*1, 0.15) {};
            \node[box2s] () at (-0.3+0.3*2, 0.15) {};

            \node[box1s] () at (-0.3+0.3*0, -0.15) {};
            \node[box2s] () at (-0.3+0.3*1, -0.15) {};
            \node[box1s] () at (-0.3+0.3*2, -0.15) {};

            \node () at (0,0.55) {\scriptsize{$B_2$}};

            \node () at (-0.3+0.3*1, 0.15) {\scriptsize{$x_2^*$}};
        \end{scope}

        \begin{scope}[xshift=1.35cm, yshift=0.325cm]
            \node[box2s] () at (-0.3+0.3*0, 0.15) {};
            \node[box1s] () at (-0.3+0.3*1, 0.15) {};
            \node[box2s] () at (-0.3+0.3*2, 0.15) {};

            \node[box1s] () at (-0.3+0.3*0, -0.15) {};
            \node[box1s] () at (-0.3+0.3*1, -0.15) {};
            \node[box2s] () at (-0.3+0.3*2, -0.15) {};

            \node () at (0,0.55) {\scriptsize{$A_2$}};
        \end{scope}

        \begin{scope}[xshift=1.8cm, yshift=-1.15cm]
            \node[box2s] () at (-0.3+0.3*0, 0.15) {};
            \node[box2s] () at (-0.3+0.3*1, 0.15) {};
            \node[box1s] () at (-0.3+0.3*2, 0.15) {};

            \node[box1s] () at (-0.3+0.3*0, -0.15) {};
            \node[box1s] () at (-0.3+0.3*1, -0.15) {};
            \node[box2s] () at (-0.3+0.3*2, -0.15) {};

            \node () at (0,-0.6) {\scriptsize{$B_3$}};

            \node () at (-0.3+0.3*2, 0.15) {\scriptsize{$x_3^*$}};
        \end{scope}

        \begin{scope}[xshift=0cm, yshift=-1.15cm]
            \node[box2s] () at (-0.3+0.3*0, 0.15) {};
            \node[box2s] () at (-0.3+0.3*1, 0.15) {};
            \node[box1s] () at (-0.3+0.3*2, 0.15) {};

            \node[box2s] () at (-0.3+0.3*0, -0.15) {};
            \node[box1s] () at (-0.3+0.3*1, -0.15) {};
            \node[box1s] () at (-0.3+0.3*2, -0.15) {};

            \node () at (0,-0.6) {\scriptsize{$A_3$}};
        \end{scope}

        \node () at (0,-2.35) {Case (i)};
    \end{tikzpicture}
    \hspace{2cm}
    \begin{tikzpicture}
        \node[v] (v1) at (-90+60*0:1.3) {};
        \node[vo] (v2) at (-90+60*1:1.3/2) {};
        \node[v] (v3) at (-90+60*2:1.3) {};
        \node[v] (v4) at (-90+60*3:1.3) {};
        \node[vo] (v5) at (-90+60*4:1.3/2) {};
        \node[v] (v6) at (-90+60*5:1.3) {};

        \draw (v1)--(v3)--(v4)--(v6)--(v1);

        \begin{scope}[xshift=-0.3cm, yshift=1.7cm]
            \node[box1s] () at (0.3*0,0.3) {};
            \node[box1s] () at (0.3*1,0.3) {};
            \node[box1s] () at (0.3*2,0.3) {};

            \node[box2s] () at (0.3*0,0) {};
            \node[box2s] () at (0.3*1,0) {};
            \node[box2s] () at (0.3*2,0) {};
            
            \node () at (1.05,0.15) {\scriptsize{$B_3$}};
        \end{scope}

        \begin{scope}[xshift=1.55cm, yshift=0.52cm]
            \node[box2s] () at (0.3*0,0.3) {};
            \node[box1s] () at (0.3*1,0.3) {};
            \node[box1s] () at (0.3*2,0.3) {};

            \node[box1s] () at (0.3*0,0) {};
            \node[box2s] () at (0.3*1,0) {};
            \node[box2s] () at (0.3*2,0) {};

            \node () at (0.3*1,0.3) {\scriptsize{$x_2^*$}};
            \node () at (0.3*2,0.3) {\scriptsize{$x_3^*$}};

            \node () at (1.05,0.15) {\scriptsize{$B_4$}};
        \end{scope}

        \begin{scope}[xshift=-2.18cm, yshift=-0.82cm]
            \node[box1s] () at (0.3*0,0.3) {};
            \node[box2s] () at (0.3*1,0.3) {};
            \node[box2s] () at (0.3*2,0.3) {};

            \node[box2s] () at (0.3*0,0) {};
            \node[box1s] () at (0.3*1,0) {};
            \node[box1s] () at (0.3*2,0) {};

            \node () at (0.3*0,0.3) {\scriptsize{$x_1^*$}};

            \node () at (-0.45,0.15) {\scriptsize{$B_2$}};
        \end{scope}

        \begin{scope}[xshift=-0.3cm, yshift=-2.0cm]
            \node[box2s] () at (0.3*0,0.3) {};
            \node[box2s] () at (0.3*1,0.3) {};
            \node[box2s] () at (0.3*2,0.3) {};

            \node[box1s] () at (0.3*0,0) {};
            \node[box1s] () at (0.3*1,0) {};
            \node[box1s] () at (0.3*2,0) {};

            \node () at (0.3*0,0) {\scriptsize{$x_1$}};
            \node () at (0.3*1,0) {\scriptsize{$x_2$}};
            \node () at (0.3*2,0) {\scriptsize{$x_3$}};

            \node () at (-0.45,0.15) {\scriptsize{$B_1$}};
        \end{scope}

        \node () at (0,-2.8) {Case (ii)};
    \end{tikzpicture}
    \caption{Two descriptions of the cycle $C$ of weight $6$ in $\cG_M$ in the proof of Lemma~\ref{lem: gamma 6cycle}.
    The cycle $C$ can be identified with a $6$-cycle in $G_M$, where solid dots represent transversal bases and hollow dots represent almost-transversal bases of $M$.}
    \label{fig: 6-cycle}
\end{figure}

\begin{lemma}\label{lem: irreducible 8cycle}
    Let $C$ be an irreducible cycle of weight $8$ in $\cG_M$.
    Then $C \simeq B_1B_2B_3B_4B_1$ for some bases $B_1,\ldots,B_4$ such that
    \begin{enumerate}[label=\rm(\roman*)]
        \item $B_2 = U+x_1^*+x_2^*+x_3+x_4$, $B_3 = U+x_1^*+x_2^*+x_3^*+x_4^*$, and $B_4 = U+x_1+x_2+x_3^*+x_4^*$, where $x_1,\ldots,x_4$ are distinct elements of $B_1$ and $U:= B_1 \setminus \{x_1,x_2,x_3,x_4\}$, and
        \item $V(\cG_M) \cap \{U\cup X : X \subseteq \{x_1,\ldots,x_4,x_1^*,\ldots,x_4^*\}\} = \{B_1,\ldots,B_4\}$.
    \end{enumerate}
\end{lemma}
\begin{proof}
    Let $C'$ be an irreducible cycle $4$-homotopic to $C$, which maximizes the number of edges $e\in E(C')$ such that $\eta(e)=2$.
    We will show that $C'$ satisfies the desired properties (i) and (ii).

    \begin{claim}\label{cl: 8irreducible 1}
        For each $B\in V(C')$, there are no three consecutive vertices $D_1,D_2,D_3$ in $C'$ such that $\eta(D_1D_2) = \eta(D_2D_3) = 1$, $\dist_M(B,D_1) = \dist_M(B,D_3) = 3$, and $\dist_M(B,D_2) = 4$.
    \end{claim}
    \begin{proof}
        Suppose that such vertices $D_1,D_2,D_3$ exist.
        By our choice of $C'$, $D_1D_3$ is not an edge in $\cG_M$.
        Thus, there is a basis $D' \ne D_2$ such that $\dist_M(D_1,D') = \dist_M(D_3,D') = 1$.
        Then $\dist_M(B,D') = 2$.
        Let $C''$ be a cycle obtained from $C'$ by replacing a vertex $D_2$ with $D'$.
        Then $C''$ is $4$-homotopic to $C'$ and it is reducible by Lemma~\ref{lem: irreducible1}, contradicting that $C'$ is irreducible.
    \end{proof}

    \begin{claim}\label{cl: 8irreducible 2}
        There are no edges $D_1D_2, D_3D_4$ of $C'$ such that $\eta(D_1D_2) = 1$, $\eta(D_3D_4) = 2 = \dist_M(D_4,D_1)$, and $\dist_M(D_2,D_3) = 3$.
    \end{claim}
    \begin{proof}
        Suppose such edges exist.
        Then for some $y_1,\ldots,y_4 \in D_1$, we have $D_2 = D_1 \symdiff\skewpair{y_3}$, $D_3 = D_1\symdiff\{y_1,y_1^*,\ldots,y_4,y_4^*\}$, and $D_4 = D_1\symdiff\{y_1,y_1^*,y_2,y_2^*\}$.
        Let $S:= D_1+y_3^*$ and $T:= D_4+y_3^*$.
        Then by Lemma~\ref{lem: fundamental circuit}, there are circuits $c\subseteq S$ and $c' \subseteq T$.

        Suppose that $D' := T-y_3 = D_4 \symdiff\skewpair{y_3}$ is a basis.
        Then $|D'\setminus D_2| = |\{x_1^*,x_2^*\}| =2$.
        Let $O$ be a cycle obtained from $C'$ by replacing an edge $D_3D_4$ with a path $D_3D'D_4$.
        Then $O$ is $4$-homotopic to $C'$ and it is reducible by Lemma~\ref{lem: irreducible1}.
        It contradicts that $C'$ is irreducible.
        Thus, $T-y_3$ is not a basis.

        Suppose that $D'' := (T-y_3) \symdiff\skewpair{y_1} = D_4\symdiff\{y_1,y_1^*,y_3,y_3^*\}$ is a basis.
        Then $|D'\setminus D_2| = 1$ and $|D'\setminus D_i| = 2$ for $i\in\{1,3,4\}$.
        Hence $C'$ is generated by four cycles of weight at most six, a contradiction.
        Thus, $(T-y_3) \symdiff\skewpair{y_1}$ is not a basis.

        By Lemma~\ref{lem: 3-term basis exchange}, $T-y_1^*$ is not a basis.
        Similarly, $T-y_2^*$ is not a basis.
        Therefore, $c' \subseteq T - \{y_1^*,y_2^*,y_3\} \subseteq S$ by Lemma~\ref{lem: fundamental circuit}.
        Note that $y_3 \in c \subseteq S$ and so $c \ne c'$.
        It contradicts Lemma~\ref{lem: fundamental circuit}.
    \end{proof}

    By Claims~\ref{cl: 8irreducible 1} and~\ref{cl: 8irreducible 2}, $C' = B_1B_2B_3B_4B_1$ such that $\eta(B_iB_{i+1})=2$ for all $i$, where $B_5:=B_1$.
    Then one can denote by $B_2 = U+x_1^*+x_2^*+x_3+x_4$, $B_3 = U+x_1^*+x_2^*+x_3^*+x_4^*$, and $B_4 = U+x_1+x_2+x_3^*+x_4^*$ for some elements $x_1,\ldots,x_4$ of $B_1$ and $U:= B_1 \setminus \{x_1,x_2,x_3,x_4\}$.
    As $C'$ is irreducible, $B_1 \symdiff\{x_i,x_i^*,x_j,x_j^*\}$ is not a basis for each $1\le i\le 2$ and $3\le j\le 4$.

    \begin{claim}
        $B_1 \symdiff\skewpair{x_1}$ is not a basis.
    \end{claim}
    \begin{proof}
        Suppose to the contrary that $B_1 \symdiff\skewpair{x_1}$ is a basis.
        Since $C'$ is irreducible, $B_4 \symdiff\skewpair{x_1}$ is not a basis.
        Let $S = B_1+x_1^*$ and $T = B_4+x_1^*$.
        Let $c$ and $c'$ be circuits of $M$ such that $c \subseteq S$ and $c' \subseteq T$.
        Then $x_1 \in c$ because $S-x_1 = B_1 \symdiff\skewpair{x_1}$ is a basis.
        Since $T-x_1 = B_4 \symdiff\skewpair{x_1}$ is not a basis, $x_1 \notin c'$.
        Hence $c \ne c'$.
        Because $(T-x_1)\symdiff \skewpair{x_3} = B_1 \symdiff\{x_1,x_1^*,x_4,x_4^*\}$ is not a basis, $T-x_3^*$ is not a basis by Lemma~\ref{lem: 3-term basis exchange}.
        Similarly, $T-x_4^*$ is not a basis.
        Hence $c' \subseteq T - x_3^* - x_4^* \subseteq S$, contradicting Lemma~\ref{lem: fundamental circuit}.
    \end{proof}

    Similarly, none of $B_i \symdiff\skewpair{x_j}$ with $i,j\in[4]$ is a basis.
    Therefore, (ii) holds.
\end{proof}

\begin{lemma}\label{lem: gamma 8cycle}
    Let $C = B_1 B_2 \ldots B_k B_1$ be a cycle of weight $8$ in $\cG_M$.
    Then $\sum_{i=1}^{k} \gamma(B_i,B_{i+1}) = 1$, where $B_{k+1} := B_1$.
\end{lemma}
\begin{proof}
    We may assume that $C$ is irreducible by Lemmas~\ref{lem: gamma 4cycle} and~\ref{lem: gamma 6cycle}.
    By Lemma~\ref{lem: irreducible 8cycle}, we can assume that $C = B_1B_2B_3B_4B_1$ and it satisfies the following:
    \begin{enumerate}[label=\rm(\roman*)]
        \item $B_2 = U+x_1^*+x_2^*+x_3+x_4$, $B_3 = U+x_1^*+x_2^*+x_3^*+x_4^*$, and $B_4 = U+x_1+x_2+x_3^*+x_4^*$, where $x_1,x_2,x_3,x_4$ are distinct elements in $B_1$ and $U:= B_1 \setminus \{x_1,x_2,x_3,x_4\}$, and 
        \item $V(\cG_M) \cap \{U\cup X : X \subseteq \{x_1,\ldots,x_4,x_1^*,\ldots,x_4^*\}\} = \{B_1,\ldots,B_4\}$.
    \end{enumerate}
    Let $S_1 = B_1+x_1^*$, $S_2 = B_2+x_1$, $S_3 = B_2+x_3^*$, and $S_4 = B_3+x_3$.
    Let $T_1 = B_1+x_3^*$, $T_2 = B_4+x_3$, $T_3 = B_4+x_1^*$, and $T_4 = B_3+x_1$; see Figure~\ref{fig: 8-cycle}.
    Then $S_1 \symdiff T_3 = S_2 \symdiff T_4 = \{x_3,x_3^*,x_4,x_4^*\}$ and $S_3 \symdiff T_1 = S_4 \symdiff T_2 = \{x_1,x_1^*,x_2,x_2^*\}$.
    Let $X_i$ and $Y_i$ be vectors in~$\cC$ such that $\ul{X_i} \subseteq S_i$ and $\ul{Y_i} \subseteq T_i$.
    By~(ii), neither $S_1-x_3$ nor $S_1-x_4$ is a basis.
    Then $\ul{X}_1 \subseteq T_3$ and thus by Lemma~\ref{lem: fundamental circuit}, $\ul{X}_1 = \ul{Y}_3$.
    Similarly, $\supp(X_{i})=\supp(Y_{i+2})$ for each $1\le i\le 4$, where the subscripts are read modulo $4$.
    Thus, %
    $X_i=c_iY_{i+2}$ for some $c_i \in F^\times$.
    Therefore, 
    for some $m \in \{0,1\}$,
    \begin{align*}
        \gamma(B_1,B_2)\gamma(B_2,B_3)
        &=
        (-1)^{m}
        \frac{X_1(x_2)}{X_1(x_1^*)}
        \frac{X_2(x_1)}{X_2(x_2^*)} 
        \frac{X_3(x_4)}{X_3(x_3^*)} 
        \frac{X_4(x_3)}{X_4(x_4^*)} \\
        &=
        (-1)^{m}
        \frac{Y_1(x_4)}{Y_1(x_3^*)}
        \frac{Y_2(x_3)}{Y_2(x_4^*)} 
        \frac{Y_3(x_2)}{Y_3(x_1^*)} 
        \frac{Y_4(x_1)}{Y_4(x_2^*)}
        =
        \gamma(B_1,B_4)\gamma(B_4,B_3). 
    \end{align*}
\end{proof}

\begin{figure}
    \centering
    \begin{tikzpicture}
        \node[v] (v1) at (-90+90*0 : 1.5) {};
        \node[v] (v2) at (-90+90*1 : 1.5) {};
        \node[v] (v3) at (-90+90*2 : 1.5) {};
        \node[v] (v4) at (-90+90*3 : 1.5) {};

        \draw (v1)--(v2)--(v3)--(v4)--(v1);
        \begin{scope}[xshift=0cm, yshift=-2.0cm]
            \node[box2s] () at (-0.45+0.3*0, 0.15) {};
            \node[box2s] () at (-0.45+0.3*1, 0.15) {};
            \node[box2s] () at (-0.45+0.3*2, 0.15) {};
            \node[box2s] () at (-0.45+0.3*3, 0.15) {};

            \node[box1s] () at (-0.45+0.3*0, -0.15) {};
            \node[box1s] () at (-0.45+0.3*1, -0.15) {};
            \node[box1s] () at (-0.45+0.3*2, -0.15) {};
            \node[box1s] () at (-0.45+0.3*3, -0.15) {};

            \node () at (0,-0.52) {\scriptsize{$B_1$}};

            \node () at (-0.45+0.3*0, -0.15) {\scriptsize{$x_1$}};
            \node () at (-0.45+0.3*1, -0.15) {\scriptsize{$x_2$}};
            \node () at (-0.45+0.3*2, -0.15) {\scriptsize{$x_3$}};
            \node () at (-0.45+0.3*3, -0.15) {\scriptsize{$x_4$}};
        \end{scope}

        \begin{scope}[xshift=-2.3cm, yshift=0cm]
            \node[box1s] () at (-0.45+0.3*0, 0.15) {};
            \node[box1s] () at (-0.45+0.3*1, 0.15) {};
            \node[box2s] () at (-0.45+0.3*2, 0.15) {};
            \node[box2s] () at (-0.45+0.3*3, 0.15) {};

            \node[box2s] () at (-0.45+0.3*0, -0.15) {};
            \node[box2s] () at (-0.45+0.3*1, -0.15) {};
            \node[box1s] () at (-0.45+0.3*2, -0.15) {};
            \node[box1s] () at (-0.45+0.3*3, -0.15) {};

            \node () at (0,-0.52) {\scriptsize{$B_2$}};

            \node () at (-0.45+0.3*0, 0.15) {\scriptsize{$x_1^*$}};
            \node () at (-0.45+0.3*1, 0.15) {\scriptsize{$x_2^*$}};
        \end{scope}

        \begin{scope}[xshift=2.3cm, yshift=0cm]
            \node[box2s] () at (-0.45+0.3*0, 0.15) {};
            \node[box2s] () at (-0.45+0.3*1, 0.15) {};
            \node[box1s] () at (-0.45+0.3*2, 0.15) {};
            \node[box1s] () at (-0.45+0.3*3, 0.15) {};

            \node[box1s] () at (-0.45+0.3*0, -0.15) {};
            \node[box1s] () at (-0.45+0.3*1, -0.15) {};
            \node[box2s] () at (-0.45+0.3*2, -0.15) {};
            \node[box2s] () at (-0.45+0.3*3, -0.15) {};

            \node () at (0,-0.52) {\scriptsize{$B_4$}};

            \node () at (-0.45+0.3*2, 0.15) {\scriptsize{$x_3^*$}};
            \node () at (-0.45+0.3*3, 0.15) {\scriptsize{$x_4^*$}};
        \end{scope}

        \begin{scope}[xshift=0cm, yshift=2.0cm]
            \node[box1s] () at (-0.45+0.3*0, 0.15) {};
            \node[box1s] () at (-0.45+0.3*1, 0.15) {};
            \node[box1s] () at (-0.45+0.3*2, 0.15) {};
            \node[box1s] () at (-0.45+0.3*3, 0.15) {};

            \node[box2s] () at (-0.45+0.3*0, -0.15) {};
            \node[box2s] () at (-0.45+0.3*1, -0.15) {};
            \node[box2s] () at (-0.45+0.3*2, -0.15) {};
            \node[box2s] () at (-0.45+0.3*3, -0.15) {};

            \node () at (0,0.52) {\scriptsize{$B_3$}};
        \end{scope}

        \begin{scope}[xshift=-4.3cm, yshift=-1.8cm]
            \node[box1s] () at (-0.45+0.3*0, 0.15) {};
            \node[box2s] () at (-0.45+0.3*1, 0.15) {};
            \node[box2s] () at (-0.45+0.3*2, 0.15) {};
            \node[box2s] () at (-0.45+0.3*3, 0.15) {};

            \node[box1s] () at (-0.45+0.3*0, -0.15) {};
            \node[box1s] () at (-0.45+0.3*1, -0.15) {};
            \node[box1s] () at (-0.45+0.3*2, -0.15) {};
            \node[box1s] () at (-0.45+0.3*3, -0.15) {};

            \node () at (-0.88,0) {\scriptsize{$S_1$}};
        \end{scope}

        \begin{scope}[xshift=-4.3cm, yshift=-0.8cm]
            \node[box1s] () at (-0.45+0.3*0, 0.15) {};
            \node[box1s] () at (-0.45+0.3*1, 0.15) {};
            \node[box2s] () at (-0.45+0.3*2, 0.15) {};
            \node[box2s] () at (-0.45+0.3*3, 0.15) {};

            \node[box1s] () at (-0.45+0.3*0, -0.15) {};
            \node[box2s] () at (-0.45+0.3*1, -0.15) {};
            \node[box1s] () at (-0.45+0.3*2, -0.15) {};
            \node[box1s] () at (-0.45+0.3*3, -0.15) {};

            \node () at (-0.88,0) {\scriptsize{$S_2$}};
        \end{scope}

        \begin{scope}[xshift=-4.3cm, yshift=0.8cm]
            \node[box1s] () at (-0.45+0.3*0, 0.15) {};
            \node[box1s] () at (-0.45+0.3*1, 0.15) {};
            \node[box1s] () at (-0.45+0.3*2, 0.15) {};
            \node[box2s] () at (-0.45+0.3*3, 0.15) {};

            \node[box2s] () at (-0.45+0.3*0, -0.15) {};
            \node[box2s] () at (-0.45+0.3*1, -0.15) {};
            \node[box1s] () at (-0.45+0.3*2, -0.15) {};
            \node[box1s] () at (-0.45+0.3*3, -0.15) {};

            \node () at (-0.88,0) {\scriptsize{$S_3$}};
        \end{scope}

        \begin{scope}[xshift=-4.3cm, yshift=1.8cm]
            \node[box1s] () at (-0.45+0.3*0, 0.15) {};
            \node[box1s] () at (-0.45+0.3*1, 0.15) {};
            \node[box1s] () at (-0.45+0.3*2, 0.15) {};
            \node[box1s] () at (-0.45+0.3*3, 0.15) {};

            \node[box2s] () at (-0.45+0.3*0, -0.15) {};
            \node[box2s] () at (-0.45+0.3*1, -0.15) {};
            \node[box1s] () at (-0.45+0.3*2, -0.15) {};
            \node[box2s] () at (-0.45+0.3*3, -0.15) {};

            \node () at (-0.88,0) {\scriptsize{$S_4$}};
        \end{scope}

        \begin{scope}[xshift=4.3cm, yshift=-1.8cm]
            \node[box2s] () at (-0.45+0.3*0, 0.15) {};
            \node[box2s] () at (-0.45+0.3*1, 0.15) {};
            \node[box1s] () at (-0.45+0.3*2, 0.15) {};
            \node[box2s] () at (-0.45+0.3*3, 0.15) {};

            \node[box1s] () at (-0.45+0.3*0, -0.15) {};
            \node[box1s] () at (-0.45+0.3*1, -0.15) {};
            \node[box1s] () at (-0.45+0.3*2, -0.15) {};
            \node[box1s] () at (-0.45+0.3*3, -0.15) {};

            \node () at (0.88,0) {\scriptsize{$T_1$}};
        \end{scope}

        \begin{scope}[xshift=4.3cm, yshift=-0.8cm]
            \node[box2s] () at (-0.45+0.3*0, 0.15) {};
            \node[box2s] () at (-0.45+0.3*1, 0.15) {};
            \node[box1s] () at (-0.45+0.3*2, 0.15) {};
            \node[box1s] () at (-0.45+0.3*3, 0.15) {};

            \node[box1s] () at (-0.45+0.3*0, -0.15) {};
            \node[box1s] () at (-0.45+0.3*1, -0.15) {};
            \node[box1s] () at (-0.45+0.3*2, -0.15) {};
            \node[box2s] () at (-0.45+0.3*3, -0.15) {};

            \node () at (0.88,0) {\scriptsize{$T_2$}};
        \end{scope}

        \begin{scope}[xshift=4.3cm, yshift=0.8cm]
            \node[box1s] () at (-0.45+0.3*0, 0.15) {};
            \node[box2s] () at (-0.45+0.3*1, 0.15) {};
            \node[box1s] () at (-0.45+0.3*2, 0.15) {};
            \node[box1s] () at (-0.45+0.3*3, 0.15) {};

            \node[box1s] () at (-0.45+0.3*0, -0.15) {};
            \node[box1s] () at (-0.45+0.3*1, -0.15) {};
            \node[box2s] () at (-0.45+0.3*2, -0.15) {};
            \node[box2s] () at (-0.45+0.3*3, -0.15) {};

            \node () at (0.88,0) {\scriptsize{$T_3$}};
        \end{scope}

        \begin{scope}[xshift=4.3cm, yshift=1.8cm]
            \node[box1s] () at (-0.45+0.3*0, 0.15) {};
            \node[box1s] () at (-0.45+0.3*1, 0.15) {};
            \node[box1s] () at (-0.45+0.3*2, 0.15) {};
            \node[box1s] () at (-0.45+0.3*3, 0.15) {};

            \node[box1s] () at (-0.45+0.3*0, -0.15) {};
            \node[box2s] () at (-0.45+0.3*1, -0.15) {};
            \node[box2s] () at (-0.45+0.3*2, -0.15) {};
            \node[box2s] () at (-0.45+0.3*3, -0.15) {};

            \node () at (0.88,0) {\scriptsize{$T_4$}};
        \end{scope}

    \end{tikzpicture}
    \caption{An illustration of a cycle $C = B_1B_2B_3B_4B_1$ of weight $8$ in $\cG_M$ together with $S_i$ and $T_i$ with $i\in[4]$ in the proof of Lemma~\ref{lem: gamma 8cycle}.}
    \label{fig: 8-cycle}
\end{figure}

\begin{proposition}\label{prop: well-defined}
    For two bases $B$ and $B'$ of $M$, let $P_1$ and $P_2$ be two paths in the basis graph $G_M$ from $B$ to $B'$.
    Then $\gamma(P_1) = \gamma(P_2)$.
\end{proposition}
\begin{proof}
    Let $C$ be a cycle in $G_M$ consisting of $P_1$ and $P_2$, and let $C'$ be the corresponding cycle in~$\cG_M$.
    By Homotopy Theorem~\ref{thm: homotopy} and Lemmas~\ref{lem: gamma 4cycle},~\ref{lem: gamma 6cycle},~\ref{lem: gamma 8cycle}, $\gamma(C') = 1$.
    Then $\gamma(C) = \gamma(C')= 1$ by Lemma~\ref{lem: gamma 4cycle}.
    Thus, $\gamma(P_1) = \gamma(P_2)$ by Lemma~\ref{lem: gamma 2cycle}.
\end{proof}

By Proposition~\ref{prop: well-defined}, the function $\varphi: \cT_n \cup \cA_n \to F$ described below Definition~\ref{def: basis graph} is well defined.
We finally show that $\varphi$ satisfies the restricted Grassmann--Pl\"{u}cker relations~\eqref{eq: restricted G--P relations F}.

\begin{theorem}\label{thm: C2B}
    $\varphi$ is a restricted G--P function.
\end{theorem}
\begin{proof}
    By Lemma~\ref{lem: gamma 2almost}, $\varphi$ satisfies~\ref{item: almost-transversals F}.

    Let $S \in \binom{E}{n+1}$ and $T\in \binom{E}{n-1}$ be sets such that ${S}$ contains exactly one skew pair, say $\skewpair{x}$, and ${T}$ has no skew pair.
    Let $\skewpair{y}$ be the unique skew pair not contained in $T$, and let $T':=T+\skewpair{y}.$
    We claim that $\varphi$ satisfies~\eqref{eq: restricted G--P relations F}. 
    We can assume that for some $z\in S\setminus T$, both ${S}-z$ and ${T}+z$ are bases of $M$.
    Then by Lemma~\ref{lem: 3-term basis exchange}, $S-x$ or $S-x^*$ is a transversal basis of $M$.
    Also, $T+y$ or $T+y^*$ is a transversal basis.
    Hence there are $X,Y\in\cC$ such that $\ul{X}\subseteq S$ and $\ul{Y}\subseteq T'$.

    By symmetry, we can assume that $S-x$ and $T+y$ are bases.
    Then $\varphi(S-x)$, $\varphi(T+y)$, $X(x)$, and $Y(y^*)$ are nonzero in $F$.
    For each $z\in S\setminus T$, we have $\frac{\varphi(S-z)}{\varphi(S-x)} = \gamma(S-x,S-z) = (-1)^{\chi(x)+\chi(z)+\smaller{S}{x}+\smaller{S}{z}}\frac{X(z)}{X(x)}$.
    We also have
    \[
        \frac{\varphi(T+y^*)}{\varphi(T+y)}
        =
        \gamma(T+y,T+y^*)
        =
        (-1)^{\chi(y)+\chi(y^*)+\smaller{T'}{y}+\smaller{T'}{y^*}}
        \frac{Y(y)}{Y(y^*)}
        =
        (-1)^{\smaller{T}{y}+\smaller{T}{y^*}}
        \frac{Y(y)}{Y(y^*)}.
    \]
    For $z\in S\setminus T'$, let $U_z := T+y+z$ and let $Y_z$ be a vector in $\cC$ such that $\ul{Y_z} \subseteq U_z$.
    Since $\sympl{Y}{Y_z} = (-1)^{\chi(y^*)}Y(y^*)Y_z(y) + (-1)^{\chi(z^*)}Y(z^*)Y_z(z) \in N_F$, we have
    \[
        \frac{\varphi(T+z)}{\varphi(T+y)}
        =
        (-1)^{\chi(y)+\chi(z)+\smaller{U_z}{y}+\smaller{U_z}{z}}
        \frac{Y_z(y)}{Y_z(z)}
        =
        (-1)^{\smaller{T}{y}+\smaller{T}{z}}
        \frac{Y(z^*)}{Y(y^*)}.
    \]
    Let $c := (-1)^{\chi(x)+\smaller{S}{x}+\smaller{T}{y}} \frac{\varphi(S-x) \varphi(T+y)}{X(x)Y(y^*)} \in F^\times$.
    Note that $X(z) = 0$ if $z\in E\setminus S$, and $Y(z^*) = 0$ if $z \in T$.
    Therefore,
    \begin{align*}
        \sum_{z\in S\setminus T}
        (-1)^{\smaller{S}{z}+\smaller{T}{z}}
        \varphi(S-z)
        \varphi(T+z)
        &=
        c \sum_{z\in S\setminus T}
        (-1)^{\chi(z)} X(z)Y(z^*) \\
        &=
        c \sum_{z\in E}
        (-1)^{\chi(z)} X(z)Y(z^*)
        \in N_F.
    \end{align*}
\end{proof}

\subsection{Equivalence}\label{sec: equiv}

In Section~\ref{sec: B2C}, we constructed an antisymmetric $F$-circuit set from an antisymmetric $F$-matroid.
Conversely, we built an antisymmetric $F$-matroid from an antisymmetric $F$-circuit set in Section~\ref{sec: C2B}.
By definition, these two constructions are the reverse step of each other, and thus we deduce Theorem~\ref{thm: bij} as follows.

\begin{proof}[\mbox{\bf Proof of Theorem~\ref{thm: bij}}]
    Let $M = [\varphi]$ be an antisymmetric $F$-matroid on $E=\ground{n}$, and let $\cC$ be the antisymmetric $F$-circuit set constructed from $\varphi$ in the sense of Section~\ref{sec: B2C}.
    Let $\varphi'$ be a restricted G--P function constructed from $\cC$ in the sense of Section~\ref{sec: C2B}.
    Then the underlying matroids of $\varphi$, $\varphi'$, and $\cC$ are the same.
    Let $B_1$ and $B_2$ be bases such that $|B_1\setminus B_2|=1$, and let $X \in\cC$ be a vector whose support $\ul{X}$ is a subset of $S:=B_1\cup B_2$.
    We denote by $\{x\} = S\setminus B_1$ and $\{y\} = S\setminus B_1$.
    Then 
    \begin{align*}
        \frac{\varphi(B_2)}{\varphi(B_1)}
        =
        (-1)^{\chi(x) + \chi(y) + \smaller{S}{x} + \smaller{S}{y}}
        \frac{X(y)}{X(x)}
        =
        \frac{\varphi'(B_2)}{\varphi'(B_1)}.
    \end{align*}
    Therefore, $M = [\varphi']$.

    Let $\cC'$ be the antisymmetric $F$-circuit set constructed from $\varphi'$.
    Then similarly we deduce that $\cC=\cC'$.
    Thus, there is a natural bijection between antisymmetric $F$-matroids and antisymmetric $F$-circuit sets.
\end{proof}

\section{Analogue of Tutte's theorem}\label{sec: tutte}

The following theorem is a fundamental result in matroid theory and linear algebra.
To our best knowledge, it was first proved by Tutte~\cite{Tutte1958} in terms of circuits and chain-group representations of matroids. %
We refer to~\cite{BJ2023} for the following statement and its proof.
A \emph{3-term Grassmann--Pl\"{u}cker relation} is a Grammann-Pl\"{u}cker relation~\eqref{eq: G--P} satisfying $|S\setminus T|=3$.

\begin{theorem}%
    \label{thm: tutte}
    For a field $k$, let $p \in \bP^{\binom{n}{r}-1}(k)$ be a point.
    Then the following are equivalent.
    \begin{enumerate}[label=\rm(\roman*)]
        \item $p$ satisfies all Grassmann--Pl\"{u}cker relations.
        \item $p$ satisfies all $3$-term Grassmann--Pl\"{u}cker relations and the support of $p$ forms a matroid.
        \item There is an $r\times n$ matrix $A$ over $k$ such that $p_B = \det(A[r,B])$ for all $B\in \binom{[n]}{r}$.
    \end{enumerate}
\end{theorem}

Now we show Theorem~\ref{thm: tutte for Lag}, an analog of Tutte's theorem for antisymmetric matroids and Lagrangian Grassmannians.
It is restated as Theorem~\ref{thm: strong and weak}.

\begin{definition} %
    \label{def: weak}
    A \emph{weak restricted G--P function on $E=\ground{n}$ with coefficients in a tract $F$} is a nontrivial function $\varphi: \cT_n \cup \cA_n \to F$ such that 
    the support of $\varphi$ form the set of bases of an antisymmetric matroid on $E$ and
    $\varphi$ satisfies~\ref{item: almost-transversals F} and the following weaker replacement of~\ref{item: restricted G--P relations F}:
    \begin{enumerate}[label=\rm(rGP$'$)]
    \item\label{item: weak} 
    For $S\in \binom{E}{n+1}$ and $T \in \binom{E}{n-1}$ such that $S$ contains exactly one skew pair and $T$ has no skew pair, if $|S\setminus T| \le 4$, then 
    \[
        \sum_{x\in S\setminus T}
        (-1)^{\smaller{S}{x} + \smaller{T}{x}}
        \varphi(S-x)
        \varphi(T+x) \in N_F.
    \]
    \end{enumerate}
    A \emph{weak antisymmetric $F$-matroid} is an equivalence class of weak restricted G--P functions.
\end{definition}

\begin{theorem}\label{thm: strong and weak}
    For a field $k$, let $\varphi \in \bP^{2^{n-2}(4+\binom{n}{2})-1}(k)$.
    Then the following are equivalent.
    \begin{enumerate}[label=\rm(\roman*)]
        \item $\varphi$ is a restricted G--P function. %
        \item $\varphi$ is a weak restricted G--P function. %
        \item There is an $n\times (\ground{n})$ matrix $A$ over $k$ such that the row-space of $A$ is in $\lag{k}{n}$ and $\varphi(B) = \det(A[n,B])$ for all $B\in \cT_n \cup \cA_n$.
    \end{enumerate}
\end{theorem}

Note that in Theorem~\ref{thm: strong and weak}(iii), the row-space of $A=\begin{bmatrix} A_1 \,|\, A_2 \end{bmatrix}$ is Lagrangian if and only if $A_1A_2^t$ is symmetric.
We show two lemmas before proving Theorem~\ref{thm: strong and weak}.

\begin{lemma}\label{lem: basis reduction}
    Let $M$ be an antisymmetric matroid on $\ground{n}$ such that $[n]$ is a basis.
    Let $X,Y \subseteq [n]$ be sets such that $|X|=|Y| \ge 2$ and $|X\setminus Y| \le 1$.
    If $[n]-X+Y^*$ is a basis, then there is $Z \subseteq X \cap Y$ such that $|X\setminus Z|  = |Y\setminus Z| \in \{1,2\}$ and $[n]-Z+Z^*$ is a basis.
\end{lemma}
\begin{proof}
    Suppose that $X=Y$.
    By~\ref{item: exchange}, $[n]-(X-e)+(Y-f)^*$ is a basis for some $e,f\in X$.
    We may assume that $e\ne f$.
    Then by Lemma~\ref{lem: 3-term basis exchange}, $[n]-(X-e)+(X-e)^*$ or $[n]-(X-e-f)+(X-e-f)^*$ is a basis.
    Therefore, we may assume that $X\ne Y$.
    We denote by $\{x\} = X-Y$ and $\{y\} = Y-X$.
    By~\ref{item: exchange}, $[n]-(X-x)+(Y-g)^*$ is a basis for some $g\in Y$.
    We may assume that $g\ne y$.
    Then by Lemma~\ref{lem: 3-term basis exchange}, $[n]-(X-x)+(Y-y)^*$ or $[n]-(X-x-g)+(Y-y-g)^*$ is a basis.
\end{proof}

\begin{lemma}\label{lem: twisting}
    For $S\subseteq [n]$, let $\Psi_s : F^{\ground{n}} \to F^{\ground{n}}$ be a linear map such that for each $i\in[n]$,
    \begin{align*}
        \mathbf{e}_i \mapsto
        \begin{cases}
            \mathbf{e}_i & \text{if $i\notin S$}, \\
            \mathbf{e}_{i^*} & \text{otherwise}, \\
        \end{cases}
        \quad\text{and}\quad
        \mathbf{e}_{i^*} \mapsto
        \begin{cases}
            \mathbf{e}_{i^*} & \text{if $i\notin S$}, \\
            -\mathbf{e}_{i} & \text{otherwise}. \\
        \end{cases}
    \end{align*}
    Then $\Psi_S$ induces a bijection from $\lag{k}{n}$ to itself such that for every $W\in \lag{k}{n}$, a set $B\in \cT_n \cup \cA_n$ is a basis of $M(W)$ if and only if $T\symdiff(S\cup S^*)$ is a basis of $M(\Psi_S(W))$.
\end{lemma}

\begin{proof}[\mbox{\bf Proof of Theorem~\ref{thm: strong and weak}}]
    Obviously, (i) implies (ii).
    By Theorem~\ref{thm: parameterization}, (i) and (iii) are equivalent.
    Hence it remains to prove that (ii) implies (iii).

    Let $\varphi$ be a weak restricted G--P function on $E := \ground{n}$ and let $M$ be its underlying antisymmetric matroid.
    By Lemma~\ref{lem: twisting},
    we may assume that $B_0 := [n]$ is a basis of $M$.
    Let $a_{ij} = (-1)^{n-i} \varphi({B_0-i+j^*}) / \varphi({B_0})$ for all $i,j\in[n]$. %
    By~\ref{item: almost-transversals F}, we have $a_{ij} = a_{ji}$ and thus
    $\Sigma := (a_{ij})_{1\le i, j\le n}$ is a symmetric matrix.

    We claim that
    \[
        \det(\Sigma[X,Y])
        =
        (-1)^{tn+\binom{t}{2}+\sum_{x\in X} x}
        \varphi(B_0 - X + Y^*)
    \]
    for all $X,Y\subseteq [n]$ such that $|X|=|Y| =: t$ and $|X\setminus Y|\le 1$.
    Note that the row-space of an $n\times E$ matrix $A:= \begin{bmatrix} I_n \,|\, \Sigma \end{bmatrix}$ is Lagrangian, and $\det(A[n,B_0-X+Y^*]) = (-1)^{tn+\binom{t}{2}+\sum_{x\in X} x} \det(\Sigma[X,Y])$.
    Thus, the claim suffices to conclude (iii).

    We prove the claim by induction on $|X|$.
    By our choice of $\Sigma$, we may assume that $|X| \ge 2$.

    \smallskip

    \noindent\textbf{Case I.}
    $X=Y$.
    By relabelling, we can assume that $X=\{1,2,\ldots,t\}$. %
    Let $m:= tn + \binom{t}{2} + \binom{t+1}{2}$.

    Suppose that none of $B_0-(X-1)+(X-i)^*$ with $i\in [m]$ is a basis.
    Then by~\ref{item: exchange}, $B_0-X+X^*$ is not a basis.
    Also, $\det(A[X-1,X-i]) = 0$ for all $i\in X$ by the induction hypothesis, and thus $\det(A[X,X])=0= \varphi({B_0-X+X^*}) / \varphi({B_0})$. %
    Therefore, we may assume that $B_0-(X-1)+(X-j)^*$ is a basis for some $j\in [m]$.
    By Lemma~\ref{lem: basis reduction}, there is $S \subseteq [t]-\{1,j\}$ such that $|\{2,\ldots,t\} \setminus S| \in \{1,2\}$ and $[n]-S+S^*$ is a basis.
    By relabelling, we may assume that $S$ is either $[t]\setminus\{1,2\}$ or $[t]\setminus\{1,2,3\}$.

    \smallskip

    \noindent\textbf{Subcase I.1.} 
    $S = [t]\setminus\{1,2\}$.
    Applying the $3$-term restricted G--P relation to $\{1,1^*,2^*\}+(B_0-X+S^*)$ and $\{2\}+(B_0-X+S^*)$,
    we have
    \begin{align*}
        \varphi(B_0-S+S^*)
        \varphi(B_0-X+X^*)
        &+
        \varphi(B_0-(S+1)+(S+1)^*)
        \varphi(B_0-(S+2)+(S+2)^*) \\
        &-
        \varphi(B_0-(S+1)+(S+2)^*)
        \varphi(B_0-(S+2)+(S+1)^*)
        =0.
    \end{align*}
    By the induction hypothesis, for $i,j\in[n]$, we have 
    \begin{align*}
        \frac{\varphi(B_0-S+S^*)}{\varphi(B_0)}
        &=
        (-1)^{m-2n-(2k-3)-3}{\det(\Sigma[S,S])}, \\
        \frac{\varphi(B_0-(S+i)+(S+j)^*)}{\varphi(B_0)}
        &=
        (-1)^{m-n-(k-1)-i}{\det(\Sigma[S+i,S+j])}.
    \end{align*}
    Then by the generalized Laplace expansion,
    \begin{align*}
        \frac{\varphi(B_0-X+X^*)}{\varphi(B_0)}
        &=
        \frac{(-1)^{m}}{\det(\Sigma[S,S])}
        \big(
            \det(\Sigma[S+1,S+1])
            \det(\Sigma[S+2,S+2]) \\
        &\hspace{2.35cm}-
            \det(\Sigma[S+1,S+2])
            \det(\Sigma[S+2,S+1])
        \big) \\
        &=
        (-1)^{m}
        \det(\Sigma[X,X]).
    \end{align*}

    \noindent\textbf{Subcase I.2.}
    $S = [t] \setminus \{1,2,3\}$.
    By the induction hypothesis, 
    for proper subsets $I,J$ of $[3]$ such that $|I|=|J|$ and $|I\setminus J| \le 1$, we have $\varphi(B_0-(S+I)+(S+J)^*)/\varphi(B_0) = (-1)^{m-\eta(I)}\det(A[S+I,S+J])$, where $\eta(I):= (n+t)(3-|I|) + \binom{|I|+1}{2}+(6-\sum_{i\in I}i)$.
    Then by the $4$-term restriced G--P relation applied to $\{1,1^*,2^*,3^*\}+(B_0-X+S^*)$ and $\{2,3\}+(B_0-X+S^*)$ and the generalized Laplace expansion, we have
    \begin{align*}
        \frac{\varphi(B_0-X+X^*)}{\varphi(B_0)}
        &=
        \frac{(-1)^{m}}{\det(\Sigma[S,S])}
        \big(
            \det(\Sigma[S+1,S+1])
            \det(\Sigma[S+2+3,S+2+3]) \\
        &\hspace{2.35cm}-
            \det(\Sigma[S+1,S+2])
            \det(\Sigma[S+2+3,S+1+3]) \\
        &\hspace{2.35cm}+
            \det(\Sigma[S+1,S+3])
            \det(\Sigma[S+2+3,S+1+2])
        \big) \\
        &=
        (-1)^{m}
        \det(\Sigma[X,X]).
    \end{align*}

    \noindent\textbf{Case II.}
    $|X\setminus Y| = 1$.
    By relabelling, we may assume that $X = [t]\setminus\{1\}$ and $Y = [t]\setminus\{2\}$.
    Let $m:= (t-1)n + \binom{t-1}{2} + \binom{t}{2}-1$.

    Suppose that none of $[n]-(X-i)+(Y-1)^*$ with $i\in X$ is a basis.
    By~\ref{item: exchange}, $B_0-X+Y^*$ is not a basis.
    For each $i\in X$, $\det(A[X-i,Y-1]) = 0$  by the induction hypothesis.
    Hence $\det(A[X,Y]) = 0 = \varphi(B_0-X+Y^*)/\varphi(B_0)$.
    Therefore, we may assume that $B_0-(X-j)+(Y-1)^*$ is a basis for some $j\in X$.
    By Lemma~\ref{lem: basis reduction}, there is $S \subseteq [t]-\{1,2,j\}$ such that $|(X-j) \setminus S| \in \{1,2\}$ and $[n]-S+S^*$ is a basis.
    By relabelling, we may assume that $S$ is either $[t]\setminus\{1,2,3\}$ or $[t]\setminus\{1,2,3,4\}$.

    \smallskip

    \noindent\textbf{Subcase II.1.} $S = [t]\setminus\{1,2,3\}$.
    By the induction hypothesis, we have $\varphi(B_0-S+S^*) = (-1)^{m-2n-(2t-3)-5} \det(\Sigma[S,S])$ and $\varphi(B_0-(S+i)+(S+j)^*) = (-1)^{m-n-(t-1)+i} \det(\Sigma[S+i,S+j])$ for each $i,j \in[3]$.
    Then by the $3$-term restricted G--P relation applied to $\{1,2,1^*,3^*\}+(B_0-[t]+S^*)$ and $\{1,3\}+(B_0-[t]+S^*)$ and the generalized Laplace expansion, we deduce that 
    \begin{align*}
        \frac{\varphi(B_0-X+Y^*)}{\varphi(B_0)}
        &=
        \frac{(-1)^{m}}{\det(\Sigma[S,S])}
        \big(
            \det(\Sigma[S+2,S+1])
            \det(\Sigma[S+3,S+3]) \\
        &\hspace{2.35cm}+
            \det(\Sigma[S+2,S+3])
            \det(\Sigma[S+3,S+1])
        \big) \\
        &=
        (-1)^{m}
        \det(\Sigma[X,Y]).
    \end{align*}

    \noindent\textbf{Subcase II.2.} $S = [t]\setminus\{1,2,3,4\}$.
    By the induction hypothesis, we have $\varphi(B_0-S+S^*) = (-1)^{m-3n-(3t-6)-9} \det(\Sigma[S,S])$.
    Also, $\varphi(B_0-(S+2)+(S+i)^*) = (-1)^{m-2n-(2t-3)+7} \det(\Sigma[S+2,S+i])$ for each $i\in[4]$, and $\varphi(B_0-(S+3+4)+(S+i+j)^*) = (-1)^{m-n-(t-1)+2} \det(\Sigma[S+3+4,S+i+j])$ for $i\in[4]$ and $j\in\{3,4\}$.
    Then by the $4$-term restricted G--P relation applied to $\{1,2,1^*,3^*,4^*\}+(B_0-[t]+S^*)$ and $\{1,3,4\}+(B_0-[t]+S^*)$ and the generalized Laplace expansion, we have
    \begin{align*}
        \frac{\varphi(B_0-X+Y^*)}{\varphi(B_0)}
        &=
        \frac{(-1)^{m}}{\det(\Sigma[S,S])}
        \big(
            \det(\Sigma[S+2,S+1])
            \det(\Sigma[S+3+4,S+3+4]) \\
        &\hspace{2.35cm}+
            \det(\Sigma[S+2,S+3])
            \det(\Sigma[S+3+4,S+1+4]) \\
        &\hspace{2.35cm}+
            \det(\Sigma[S+2,S+3])
            \det(\Sigma[S+3+4,S+1+3])
        \big) \\
        &=
        (-1)^{m}
        \det(\Sigma[X,Y]).
    \end{align*}
\end{proof}

Note that Theorem~\ref{thm: strong and weak} fails if we weaken the definition of weak restricted G--P functions by removing the condition that the support of $\varphi$ forms an antisymmetric matroid.
For example, let $F$ be an arbitrary tract and $\varphi:\cT_4 \cup \cA_4 \to F$ be a function such that $\supp(\varphi) = \{[4],[4]^*\}$.
Then $\supp(\varphi)$ is not the set of bases of an antisymmetric on $\ground{4}$ and thus $\varphi$ is not a restricted G--P function.
However, $\varphi$ satisfies all $3$- and $4$-term restricted G--P relations.
Similarly, if we moderated the definition of weak restricted G--P functions by replacing the $3$-/$4$-term restricted G--P relations with the $3$-term restricted G--P relations, then Theorem~\ref{thm: strong and weak} does not hold anymore; see Example~\ref{eg: 4term}.
We remark that Tutte's theorem was extended for \emph{perfect tracts} including all partial fields and the sign, tropical, Krasner hyperfields by Baker and Bowler~\cite{BB2019}.
Also, an analogue of Tutte's theorem holds for even symmetric matroids and Lagrangian orthogonal Grassmannians over partial fields~\cite{BJ2023} and the tropical hyperfield~$\bT$~\cite{Rincon2012}.

\begin{example}\label{eg: 4term}
    Let $M = (\ground{4},\cB)$ be an antisymmetric matroid such that
    \[
        \cB := 
        \{[4]\} \cup 
        \{ i^*j^*k\ell : ijkl = [4] \} \cup 
        \Big\{ ii^*jk : ijk \in \binom{[4]}{3} \Big\} \cup
        \Big\{ ii^*jk^* : ijk \in\binom{[4]}{3} \Big\}.
    \]
    Let $\varphi: \cT_n \cup \cA_n \to \bF_2$ be a function whose support is $\supp(\varphi)=\cB$.
    Then $\varphi$ satisfies all $3$-term restricted G--P relations, but it does not satisfy $4$-term restricted G--P relations.
    More precisely, it does not satisfy a $4$-term restricted G--P relation~\eqref{eq: restricted G--P relations F} applied to $S=\{1,2,2^*,3^*,4^*\}$ and $T=\{1^*2,3\}$. Both $S-x$ and $T+x$ are bases for each $x\in \{2^*,3^*,4^*\}$, and $S-1$ is not a basis of $M$. %
    Hence $\sum_{x\in S-T} \varphi(S-x)\varphi(T+x) = 3 \ne 0$.
\end{example}

\section{Concluding remarks}\label{sec: concluding remarks}

Finally, we propose several questions motivated by matroid theory.

\begin{itemize}
    \item Can we find more cryptomorphic definitions of antisymmetric matroids?
    The author is especially interested in establishing the hyperplane or flat axiom for antisymmetric matroids.

    \item 
    By Proposition~\ref{prop: repre and minors}, the class of antisymmetric matroids representable over a field $k$ is closed under taking minors.
    Thus, it is natural to ask for the excluded minors for antisymmetric matroids representable over a certain field.
    The easiest case would be the binary field $\bF_2$.
    We note that an antisymmetric matroid $(\ground{2}, \cT_2 \cup \cA_2)$ is non-binary.
    
    \item 
    The $3$-term Grassmannian-Pl\"{u}cker relations provide a weaker notion of matroids over tracts, called \emph{weak $F$-matroids}~\cite{BB2019}.
    Remarkably, weak $F$-matroids are $F$-matroids in most interesting tracts $F$ such as all partial fields, the sign hyperfield $\bS$, and the tropical hyperfield $\bT$.
    Thus, it would be interesting to ask whether weak antisymmetric $F$-matroids are antisymmetric $F$-matroids if $F=\bS$ or $\bT$.
    Theorem~\ref{thm: strong and weak} proves this when $F$ is a field, and the proof can be easily generalized for partial fields.
    We write its outline below:
    
    Let $\varphi$ be a weak restricted G--P function over a partial field~$P$, and we regard $P$ as a pair of a commutative ring~$R$ with unity and a unit subgroup~$G$; see Example~\ref{eg: tracts}\ref{item: partial field}.
    As the same way described in the proof of Theorem~\ref{thm: strong and weak}, we construct a symmetric matrix $\Sigma = (a_{ij})_{1\le i,j\le n}$ over $R$ such that each entry $a_{ij}$ is in~$G$, and then the main claim still holds, i.e., $\det(\Sigma[X,Y]) = (-1)^{tn+\binom{t}{2}+\sum_{x\in X} x} \varphi(B_0 - X + Y^*)$ for all $X,Y\subseteq [n]$ such that $B_0-X+Y^*\in \cT_n\cup \cA_n$.
    Then the Laplace expansions for $\Sigma$ imply that $\varphi$ is a restricted G--P function.

    \item 
    We reveal the connections between matroids and antisymmetric matroids in Sections~\ref{sec: antisymmetric matroids circuits}--\ref{sec: representability} and~\ref{sec: matroids with coefficients}.
    We showed that every matroid can be identified with an antisymmetric matroid extending its representability.
    We ask a kind of inverse question:
    For an antisymmetric matroid $M$ on $E=\ground{n}$, is there a matroid $N$ on $E$ such that $\cB(M) = \cB(N) \cap (\cT_n \cup \cA_n)$?
    This question is true for antisymmetric matroids representable over a field, and it is motivated by \emph{enveloping matroids} in~{\cite[Section~3]{BGW2003}}.
    Remark that the same question is open for even symmetric matroids.

    \item 
    Theorem~\ref{thm: even to antisymmetric} shows that every even symmetric matroid can be obtained from an antisymmetric matroid by discarding almost-transversal bases.
    It would be interesting to extend this result for general symmetric matroids or to find a counterexample.

    \item 
    Extending any results on oriented matroids to oriented antisymmetric matroids, such as the topological representation theorem, would be appealing.

    \item
    Corollary~\ref{cor: realizable positively oriented antisym mat} shows that oriented antisymmetric matroids with certain sign conditions are realizable, as a consequence of a result on oriented gaussoid~\cite{BDKS2019}.
    We ask whether the condition~(i) in Corollary~\ref{cor: realizable positively oriented antisym mat} can be weakened to that $\varphi([n]-X+X^*) = 0$ or $(-1)^{\sum_{i\in X} (i+\smallereq{X}{i})}$ for every $X\subseteq [n]$.
    Note that, if true, then the corresponding oriented antisymmetric matroid is realizable by a positive semidefinite symmetric matrix.

\end{itemize}

\section*{Funding}
This work was supported by the Institute for Basic Science (IBS-R029-C1).

\section*{Acknowledgements}
The author thanks Changxin Ding, Chris Eur, Matthew Baker, Oliver Lorscheid, June Huh, and Sang-il
Oum for their valuable comments. In particular, Changxin Ding pointed out minor errors and helped
improve readability. The author also thanks the referees for their helpful comments.

\providecommand{\bysame}{\leavevmode\hbox to3em{\hrulefill}\thinspace}
\providecommand{\MR}{\relax\ifhmode\unskip\space\fi MR }
\providecommand{\MRhref}[2]{%
  \href{http://www.ams.org/mathscinet-getitem?mr=#1}{#2}
}
\providecommand{\href}[2]{#2}

\end{document}